\RequirePackage{fix-cm}
\documentclass[smallcondensed,numbook,envcountsect]{svjour3}% 

\smartqed  % flush right qed marks, e.g. at end of proof
\usepackage{soul}
\usepackage[normalem]{ulem}
\usepackage{graphicx}%
\usepackage{multirow}%
\usepackage{amsmath,amssymb,amsfonts}%
\usepackage{mathrsfs}%
\usepackage[dvipsnames]{xcolor}
\usepackage{textcomp}%
\usepackage{manyfoot}%
\usepackage{booktabs}%
\usepackage{algorithm}%
\usepackage{algorithmicx}%
\usepackage{algpseudocode}%
\usepackage{listings}%
\usepackage{tabularx}
\usepackage{todonotes}
\usepackage{orcidlink}
\usepackage{tikz}
\usetikzlibrary{arrows.meta, positioning, shapes}
\definecolor{vibrantYellow}{RGB}{255,215,0}
\definecolor{brightOrange}{RGB}{255,165,0}
\definecolor{limeGreen}{RGB}{50,205,50}
\definecolor{crimsonRed}{RGB}{220, 20, 60}
\definecolor{royalBlue}{RGB}{65, 105, 225}
\definecolor{darkOrchid}{RGB}{153, 50, 204}
\definecolor{Teal}{RGB}{0, 128, 128}
\definecolor{tomatoRed}{RGB}{255, 99, 71}
\definecolor{chartreuse}{RGB}{127, 255, 0}

\usepackage{cancel}

\usepackage[title]{appendix}%
\usepackage[labelfont=bf]{caption}
\usepackage[labelfont=bf]{subcaption}
\usepackage{mathtools}
\usepackage{nicefrac, xfrac}

\usepackage[bb=dsserif]{mathalpha}
\usepackage{makecell}

\usepackage{pifont}% http://ctan.org/pkg/pifont
\newcommand{\n}{{\rho}}
\newcommand{\m}{{\eta}}
\DeclareMathOperator{\sgn}{sgn}
\usepackage[utf8]{inputenc}
\usepackage{ascii}
\usepackage{newunicodechar}
\newunicodechar{☺}{\SOH}
\newunicodechar{☻}{\STX}
\newunicodechar{♥}{\ETX}
\newunicodechar{♣}{\ENQ}
\newunicodechar{♠}{\ACK}
\newunicodechar{○}{\HT}

\makeatletter
\newsavebox\myboxA
\newsavebox\myboxB
\newlength\mylenA
\newcommand*\xoverline[2][0.74]{%
    \sbox{\myboxA}{$\m@th#2$}%
    \setbox\myboxB\null% Phantom box
    \ht\myboxB=\ht\myboxA%
    \dp\myboxB=\dp\myboxA%
    \wd\myboxB=#1\wd\myboxA% Scale phantom
    \sbox\myboxB{$\m@th\overline{\copy\myboxB}$}%  Overlined phantom
    \setlength\mylenA{\the\wd\myboxA}%   calc width diff
    \addtolength\mylenA{-\the\wd\myboxB}%
    \ifdim\wd\myboxB<\wd\myboxA%
       \rlap{\hskip 0.25\mylenA\usebox\myboxB}{\usebox\myboxA}%
    \else
        \hskip -0.25\mylenA\rlap{\usebox\myboxA}{\hskip 0.25\mylenA\usebox\myboxB}%
    \fi}
\makeatother
\usepackage{xparse}
\NewDocumentCommand{\Iij}{ O{i} O{j} }{\ensuremath{Q_{[#1,#2]}}}
\NewDocumentCommand{\Iji}{ O{i} O{j} }{\ensuremath{Q_{[#2,#1]}}}
\NewDocumentCommand{\alphaij}{ O{i} O{j} }{\ensuremath{\alpha_{[#1,#2]}}}
\NewDocumentCommand{\alphaji}{ O{i} O{j} }{\ensuremath{\alpha_{[#2,#1]}}}
\NewDocumentCommand{\TigLmu}{ O{k} O{\ensuremath{\gamma L}} O{\ensuremath{\kappa}} }{\ensuremath{T_{#1} (#2, #3)}}
\NewDocumentCommand{\Ei}{ O{k} O{x} }{\ensuremath{E_{#1} (#2)}}
\NewDocumentCommand{\Tinm}{ O{k} }{\ensuremath{T_{#1}(\n,\m)}}
\NewDocumentCommand{\gbari}{ O{k} O{\ensuremath{(\kappa})} }{\ensuremath{\xoverline{\gamma L}_{#1}{#2}}}
\NewDocumentCommand{\Nbar}{ }{\ensuremath{\bar{N}(\gamma L, \kappa)}}
\NewDocumentCommand{\PN}{ }{\ensuremath{P_N(\gamma L, \gamma \mu)}}
\NewDocumentCommand{\gstar}{ O{\ensuremath{\kappa}} }{\ensuremath{(\gamma L)^*{#1}}}
\NewDocumentCommand{\pigLgmu}{ O{\ensuremath{\gamma}} }{\ensuremath{p({#1}L, {#1}\mu)}}
\NewDocumentCommand{\idproof}{ O{\ensuremath{k}} }{\ensuremath{#1}}
\NewDocumentCommand{\thetagLgmu}{ O{\ensuremath{\gamma}} }{\ensuremath{\theta(#1 L, #1 \mu)}}
\NewDocumentCommand{\cgLgmu}{ O{\ensuremath{\gamma}} }{\ensuremath{c(#1 L, #1 \mu)}}
\NewDocumentCommand{\lhs}{ }{l.h.s.\xspace}
\NewDocumentCommand{\rhs}{ }{r.h.s.\xspace}

\newcommand\Tstrut{\rule{0pt}{3ex}}       % "top" strut
\newcommand\Bstrut{\rule[-0.5ex]{0pt}{0pt}} % "bottom" strut
\newcommand{\TBstrut}{\Tstrut\Bstrut} % top&bottom struts
\usepackage{enumitem}
\setlist[itemize]{align=parleft,left=0pt..1em}
\usepackage{comment}

\graphicspath{{./pics/}}
\DeclareMathOperator*\minimize{minimize}%   minimize	
\DeclareMathOperator*\argmax{arg\,max}  %   argmax
\DeclareMathOperator*\argmin{arg\,min}  %	argmin
\AtBeginEnvironment{appendices}{\crefalias{section}{appendix}}

\usepackage{varioref}
\usepackage{hyperref}
\usepackage[noabbrev, capitalise, nameinlink]{cleveref}

\makeatletter
\let\cl@chapter\undefined
\makeatletter
\numberwithin{theorem}{section}% Number theorems within sections: <section>.<theorem>
\makeatletter
\let\c@corollary\relax% Corollary
\let\c@proposition\relax% Proposition
\let\c@lemma\relax% Lemma
\let\c@definition\relax% Definition
\let\c@remark\relax% Remark
\let\c@conjecture\relax% Conjecture
\makeatother
\usepackage{aliascnt}
\newaliascnt{corollary}{theorem}% Corollary numbering: <section>.<theorem>
\newaliascnt{proposition}{theorem}% Proposition numbering: <section>.<theorem>
\newaliascnt{lemma}{theorem}% Lemma numbering: <section>.<theorem>
\newaliascnt{definition}{theorem}% Definition numbering: <section>.<theorem>
\newaliascnt{remark}{theorem}% Remark numbering: <section>.<theorem>
\newtheorem{conjecture}{Conjecture}
\hypersetup{%
      bookmarksnumbered,
      breaklinks,
      hypertexnames = true,
      plainpages    = false,
      colorlinks = true,
      urlcolor   = blue,
      linkcolor  = MidnightBlue,
      filecolor  = Mulberry,      % color of file links
      citecolor  = BrickRed,
      menubordercolor = BrickRed
    }

\DeclareMathOperator*{\sSum}{\textstyle\sum}
\usepackage{relsize}% http://ctan.org/pkg/relsize
\newif\ifshowoldtext
\showoldtextfalse     % To DELETE the old text completely
\newcommand{\oldtext}[1]{%
  \ifshowoldtext
    {\color{gray}{#1}}% The '%' prevents unwanted spaces
  \fi
}
\newif\ifhighlight
\highlightfalse     % To show the text in the document's current color
\newcommand{\shownewtext}[1]{%
  \ifhighlight
    {\color{BrickRed}{#1}}% If the switch is true, color the text red.
  \else
    #1% If the switch is false, just output the text as is.
  \fi
}
\newif\ifhighlightfinal
\highlightfinalfalse     % To show the text in the document's current color
\newcommand{\shownewtextfinal}[1]{%
  \ifhighlightfinal
    {\color{red}{#1}}% If the switch is true, color the text red.
  \else
    #1% If the switch is false, just output the text as is.
  \fi
}
\makeatletter
\renewcommand\thefigure{\@arabic\c@figure}
\renewcommand\thetable{\@arabic\c@table}
\@removefromreset{figure}{section}
\@removefromreset{table}{section}
\makeatother
\newcommand{\crefdefpart}[2]{%
  \hyperref[#2]{\namecref{#1}~\labelcref*{#1}\ref*{#2}}%
}
\begin{document}
\title{Exact worst-case convergence rates of gradient descent: a complete analysis for all constant stepsizes over nonconvex and convex functions
% A complete analysis of the gradient descent: \\
% the exact worst-case rates for convex, nonconvex and strongly convex functions for any constant stepsize
%\thanks{Grants or other notes
%about the article that should go on the front page should be
%placed here. General acknowledgments should be placed at the end of the article.}
}%
%\subtitle{Do you have a subtitle?\\ If so, write it here}

\titlerunning{Exact performance bounds of gradient descent}        % if too long for running head

\author{Teodor Rotaru \and
        Fran\c cois  Glineur \and
        Panagiotis Patrinos}

%\authorrunning{Short form of author list} % if too long for running head

\institute{Teodor Rotaru \,\orcidlink{0000-0003-2039-6228} \and Panos Patrinos \,\orcidlink{0000-0003-4824-7697} \at
              Department of Electrical Engineering (ESAT-STADIUS), KU Leuven, Kasteelpark Arenberg 10, Leuven, 3001, Belgium \\
              \email{\emph{teodor.rotaru@kuleuven.be} $\cdot$ panos.patrinos@esat.kuleuven.be}
           \and
           Teodor Rotaru \,\orcidlink{0000-0003-2039-6228} \and Fran\c cois  Glineur \,\orcidlink{0000-0002-5890-1093} \at
              Department of Mathematical Engineering (ICTEAM-INMA), UCLouvain, Av. Georges Lemaître 4, Louvain-la-Neuve, 1348, Belgium \\
              \email{\emph{teodor.rotaru@kuleuven.be} $\cdot$ francois.glineur@uclouvain.be}
}%
% \vspace{-2em}
\date{\vspace{-1em}Received: date / Accepted: date}
% \date{Received: 16 August 2024 / Accepted: 29 November 2025}
% The correct dates will be entered by the editor
\maketitle
\vspace{-.1em}
\begin{abstract}
We consider gradient descent with constant stepsizes and derive exact worst-case convergence rates on the minimum gradient norm of the iterates. Our analysis covers all possible stepsizes and arbitrary upper/lower bounds on the curvature of the objective function, thus including convex, strongly convex and weakly convex (hypoconvex) objective functions. 

Among the challenging parts of the analysis,
we note the necessity to exploit dependencies between non-consecutive iterates. 
While this complicates the proofs to some extent, it enables us to achieve an exact full-range analysis of gradient descent for any constant stepsize (covering, in particular, normalized stepsizes greater than one), whereas the literature contained only conjectured rates of this type. %
In the nonconvex case, allowing arbitrary bounds on upper and lower curvatures extends existing partial results that are valid only for gradient Lipschitz functions (i.e., where lower and upper bounds on curvature are equal), leading to improved rates for weakly convex functions. 

From our exact worst-case performance bounds, we deduce the optimal constant stepsize for gradient descent. Leveraging our analysis, we also introduce a new variant of gradient descent based on a unique, fixed sequence of variable stepsizes, demonstrating its superiority in the worst-case over any constant stepsize schedule.%
% \PACS{PACS code1 \and PACS code2 \and more}
\smallskip\noindent\textbf{Funding} This work was supported by the framework of the Global PhD Partnership KU Leuven -- UCLouvain. 
\end{abstract}

%%%%%%%%%%%% TO REMOVE AT THE END %
\clearpage
\tableofcontents
\clearpage 
%%%%%%%%%%%% TO REMOVE AT THE END %
%
%
\section{Introduction}\label{sec:Introduction}
Tight convergence rates for the canonical first-order optimization algorithm, i.e., the gradient descent \shownewtext{method}, applied to smooth functions are demonstrated. %
We analyse functions with arbitrary upper and lower curvatures, covering weakly convex, convex and strongly convex cases. For all constant stepsizes ensuring descent, we derive tight bounds and establish convergence rates parameterized by curvature bounds, stepsize, and iteration count. 
The key tool enabling to obtain a complete characterization is the analysis of iterations' interconnection not only for consecutive steps, but also for the ones lying at distance two, unlike classical convergence analysis, e.g., \cite{Nesterov_cvx_2018,Boyd_cvx_opt,bertsekas2016nonlinear}. 
% \shownewtextfinal{Our proof strategy differs from existing work on performance estimation problems (PEP) (e.g., \cite{Kim_Fesler_2021_OGM-G}), which connects consecutive iterates to the final one. In contrast, our approach allows for an analysis decoupled from the last iterate, using standard descent lemmas to establish performance bounds. Results on full range of stepsizes were previously at most conjectured for the apparently simpler scenario of convex functions.}%

Standard PEP formulations typically rely on non-consecutive iterates to characterize worst-case performance. However, incorporating all resulting dependencies generally leads to substantially harder proofs; consequently, in practice, relatively few analyses leverage more than consecutive iterates. %

Our proof strategy slightly differs from existing work on PEP (e.g., \cite{Kim_Fesler_2021_OGM-G}), which connects consecutive iterates to the final one, as well as from the recent line of work on stepsize-based acceleration (e.g., via silver stepsizes), which exploits only a subset of inter-iterate connections. %
In contrast, our approach enables an analysis that is decoupled from the final iterate, relying instead on standard descent lemmas to establish performance bounds. To the best of our knowledge, results covering the full range of stepsizes had previously been, at most, conjectured, even in the apparently simpler scenario of convex objectives.%

\subsection{Notations and Definitions}%
We consider the unconstrained optimization problem %
$$\minimize_{x \in \mathbb{R}^d} f(x),$$%
where $f$ belongs to the class of functions $\mathcal{F}_{\mu, L}$, as defined below.%
\begin{definition}\label{def:upper_lower_curvature} 
Let $L \in (0, \infty)$ and $\mu\in (-\infty,L)$. Then a function $f\colon\mathbb{R}^d\rightarrow \mathbb{R}$ belongs to the class of functions $\mathcal{F}_{\mu,L}$ if it has:%
\begin{enumerate}[label=(\roman*)]%
    \item upper curvature $L$, i.e., $ \tfrac{L}{2} \|\cdot\|^2 - f$ is convex;%
    \item lower curvature $\mu$, i.e., $ f - \tfrac{\mu}{2} \|\cdot\|^2$ is convex. %
\end{enumerate}%
\end{definition}%
Intuitively, if $f \in \mathcal{C}^2$, then the eigenvalues of its Hessian matrix belong to some interval $[\mu,L]$. 
Depending on the sign of $\mu$, $f$ is categorized as: (i) weakly convex (or hypoconvex) for $\mu<0$, (ii) convex for $\mu=0$ or (iii) strongly convex for $\mu>0$. 
The trivial case $\mu = L$ is not explicitly addressed; it leads to a simple quadratic problem, whose Hessian with all eigenvalues equal to $L$ and for which the answers are easy to derive (for a one-dimensional example, see \cref{prop:tightness_linear_regime}). 
In our framework, we allow $\mu < -L$, implying that the Lipschitz constant of the gradients is $\max\{-\mu,L\}$. 
The class of smooth weakly convex functions includes the smooth not-necessarily convex ones under the condition $\mu=-L$. 
We sometimes represent the curvature ratio $\kappa{}\coloneqq{}\frac{\mu}{L} \in (-\infty, 1)$, which, if $\mu \geq 0$, signifies the inverse condition number.%

\noindent {\textbf{Gradient descent.}} 
Starting from some initial point $x_0 \in \mathbb{R}^d$, we consider $N$ iterations of gradient descent with fixed stepsizes $\gamma_i \in (0, \tfrac{2}{L})$, defined as:
\begin{align}\label{eq:GM_it_fixed_steps}\tag{GD}
    x_{i+1} = x_i - \gamma_i \nabla f(x_i), \quad \forall i = 0,\dots,N-1.
\end{align}%
The constraint on the stepsize values ensures a decrease in the function value after a single iteration. We assume access to a first-order oracle denoted by the set of triplets $\mathcal{T}{}\coloneqq{}\{(x_i,g_i,f_i)\}_{i \in \mathcal{I}} \subseteq \mathbb{R}^d\times\mathbb{R}^d\times\mathbb{R}$, where $\mathcal{I}$ is an index set, providing information about the iterates, their gradients and the corresponding function values. We assume the function $f$ to be bounded from below and denote $f_* {}\coloneqq{} \inf_x f(x) > -\infty$. %
In the nonconvex case, where first-order methods cannot guarantee convergence to a global minimizer, analysis focuses on the stationarity measure $ \min_{0 {} \leq {} i {} \leq {} N} \{\|\nabla f(x_i)\|^2\}$, typically bounded in terms of the initial optimality global gap $f(x_0)-f_*$, or, when $f_*$ is unknown, by the observable decrease $f(x_0)-f(x_N)$. 
We denote $[x]_{+} {}\coloneqq{} \max\{x, 0\}$ for any $x \in \mathbb{R}$.
\begin{remark}
A tight (exact) performance bound is established in two steps: %
(i) providing an \textit{upper bound} of the convergence measure and %
(ii) constructing a problem instance (a \textit{lower bound}) achieving it. %
For convenience, the performance measure is often called \textit{(convergence) rate}. This term refers not only to the asymptotic expression, such as $ \mathcal{O}(\cdot)$, but also to its associated, and often exact, constant.%
\end{remark}
\subsection{Prior work}%
The conventional analysis of gradient descent relies on a constant stepsize schedule and provides non-tight rates (for example, see the textbooks \cite{Boyd_cvx_opt,bubeck2014convex,bertsekas2016nonlinear,Beck_first_order_methods_opt,Nesterov_cvx_2018}). Drori and Teboulle's seminal paper on performance estimation problem (PEP) \cite{drori_performance_2014} propelled the quest for precise convergence rates in first-order optimization methods, due to its systematic approach to analyze worst-case behaviors. The PEP frames the task of identifying an algorithm's most unfavorable behavior within a specific problem class as an optimization problem. 
Using relaxation techniques, they formulate the problem as a convex semidefinite program. Taylor et al. \cite{Taylor_2017_SIAM_Composite_convex,taylor_smooth_2017} advanced the framework by specifying conditions on when the relaxation is tight. Notably, they introduced interpolation inequalities as a necessary and sufficient tool for deriving exact performance bounds. The conditions they introduced for $\mathcal{F}_{\mu,L}$ functions are central to our derivations.

Convergence rates of first-order methods applied to smooth convex functions have been extensively examined using PEP. For gradient descent, in \cite{drori_performance_2014} tight rates are proved for smooth convex functions with stepsizes $\gamma L \in (0,1]$ and conjectured for $\gamma L \in (1,2)$; exact rates for smooth strongly convex functions are conjectured for $\gamma L \in (0,2)$ \cite{taylor_smooth_2017} and proved when employing line-search \cite{deKlerk_Taylor_line_search_2017}. Optimized first-order methods can be derived within the PEP framework, e.g., in \cite{Kim_Fesler_2021_OGM-G} for an optimal algorithm to decrease the gradient norm of smooth convex functions. A recent breakthrough by Teboulle \cite{Teboulle2022} establishes tight convergence rates for smooth convex functions using stepsizes $\gamma L \in (1, \frac{3}{2}]$.%

For smooth nonconvex functions, commonly encountered in fields like machine learning, there is a notable gap in comprehending convergence rates. Carmon et al. establish lower bounds for finding stationary points in \cite[Theorem 2]{Carmon2020_lower_bounds_grad_nrm}, demonstrating that gradient descent is worst-case optimal for functions with only a Lipschitz gradient. In a related study, Cartis et al. present tight lower bounds for steepest descent \cite{Cartis_complexity_steepest_descent}. Within the PEP framework, Taylor examines the convergence rates of gradient descent applied to smooth, not-necessarily convex functions, i.e., $\mu = -L$, with a constant stepsize $\frac{1}{L}$ \cite[p. 190]{PhD_AT_2017}. Drori and Shamir extend this result to the case of variable, yet predetermined, stepsizes below $\tfrac{1}{L}$ \cite[Corollary 1]{drori2021complexity} (using scaling arguments, this extension could also be derived from the result in \cite[p. 190]{PhD_AT_2017}). Abbaszadehpeivasti et al. further improve upon these findings \cite{abbaszadehpeivasti2021GM_smooth}, achieving state-of-the-art for stepsizes up to $\tfrac{\sqrt{3}}{L}$. Typically, these aforementioned results solely rely exploiting smoothness, neglecting the weak convexity parameter $\mu$. Recent interest in weakly convex optimization has emerged, e.g., \cite{Davis_Drusvyatskiy_weakly_convex,Solodov_weakly_cvx_2023}.%

Given our focus on constant stepsize schedule, our analysis is confined to the interval $\gamma L \in (0,2)$. However, recent advancements address longer stepsizes and achieve improved asymptotic rates. Interested readers are directed to the works of Altschuler and \shownewtextfinal{Parrilo} \cite{altschuler2023acceleration_cvx,altschuler2023acceleration_strg_cvx}, Das Gupta et al. \cite{DasGupta2024} and Grimmer et al. \cite{grimmer2023accelerated,grimmer2024provably,grimmer2024accelerated_compact}. Notably, similar to our approach, proving their rates requires connecting more than just consecutive iterations, deviating from classical analyses. Comparable demonstrations are provided in \cite{Kim_Fesler_2021_OGM-G,zamani2023exact}.

Nesterov highlighted the significance of minimizing the gradient norm \cite{Nesterov_how_to_make_gradients_small}, asserting that $\|\nabla f(x)\|$ might serve better than other performance metrics for minimization purposes. Evaluating the gradient norm for initially bounded functions is a common approach in analyzing gradient-based methods for smooth nonconvex optimization, seen in works such as \cite{drori2021complexity} for stochastic gradient descent and in \cite{Kim_Fesler_2021_OGM-G} for the OGM-G algorithm optimizing the gradient norm for smooth convex functions. This convergence measure is suitable for applications where both the initial optimality gap and the gradient norm across iterations can be directly measured, unlike metrics involving the optimizer. We maintain consistency by utilizing this performance metric for both convex and strongly convex functions. %
%%%%%%%%%%%%%%%%%%%%%%%%%%%%%%%%%%%%%%%%%%
%%%%%%%%%%%%%%%%%%%%%%%%%%%%%%%%%%%%%%%%%%
%
\subsection{Paper organization and main contributions}\label{subsec:Contributions}%
Our contributions can be summarized as follows:%
\begin{enumerate}
    \item We provide tight performance bounds for gradient descent with constant stepsizes $\gamma L \in (0,2)$, proving their tightness through worst-case examples. Our analysis applies to 
    weakly convex (\cref{thm:wc_GM_hypo_ct_all}), 
    strongly convex (\cref{thm:strg_cvx_rate}), and 
    convex functions (\cref{corollary:cvx_rate}), with a comprehensive summary in \cref{table:summary_bounds_and_tightness}.%
    \item A performance bound for not-necessarily convex functions using variable stepsizes $\gamma_i L \in (0,\gbari[1]]$\footnote{The quantities $\gbari[k]$ denote normalized stepsize thresholds that depend on the curvature ratio $\kappa$. They are monotonically increasing in the index $k \geq 0$ and are formally defined in \cref{def:stepsize_thresholds}.} (\cref{thm:wc_GM_hypo}), which we prove is tight for stepsizes below 1 and conjecture for the remaining range. This result generalizes the constant stepsize case for weakly convex functions (\cref{thm:wc_GM_hypo_ct_all}).%
    \item Leveraging our tight performance bounds, we propose optimized stepsize rules for gradient descent, namely: %
    (i) best constant stepsize for 
    convex, strongly convex and weakly convex functions (\cref{thm:optimal_ct_stepsize_and_rate_convex,thm:optimal_ct_stepsize_and_rate_strongly_convex,prop:gamma_star}) for a given number of iterations; %
    (ii) a dynamic sequence of stepsizes, independent of the iteration count $N$, achieving a superior guarantee compared to any constant stepsize policy (\cref{prop:tight_rate_all_kappas_increasing_sequence}).%
\end{enumerate}%
The aforementioned works typically address only the convex case ($\mu=0$), the strongly convex case ($\mu \in (0,L)$), the Lipschitz smooth case ($\mu=-L$), or the setting without lower curvature information ($\mu = -\infty$) (see, e.g., \cite[Section 1.2.3]{Nesterov_cvx_2018}). Our framework unifies these regimes and yields a continuous spectrum of performance guarantees interpolating between them. 
We recall that the performance metric for (strongly) convex functions usually differs, e.g., in \cite[Conjecture 1]{drori_performance_2014}, \cite[Conjecture 3]{taylor_smooth_2017} and \cite[Theorem 2.1]{Taylor_Jota_PGM_rates_proofs}.%

\begin{figure}%
    \centering%
    \includegraphics[width=\linewidth]{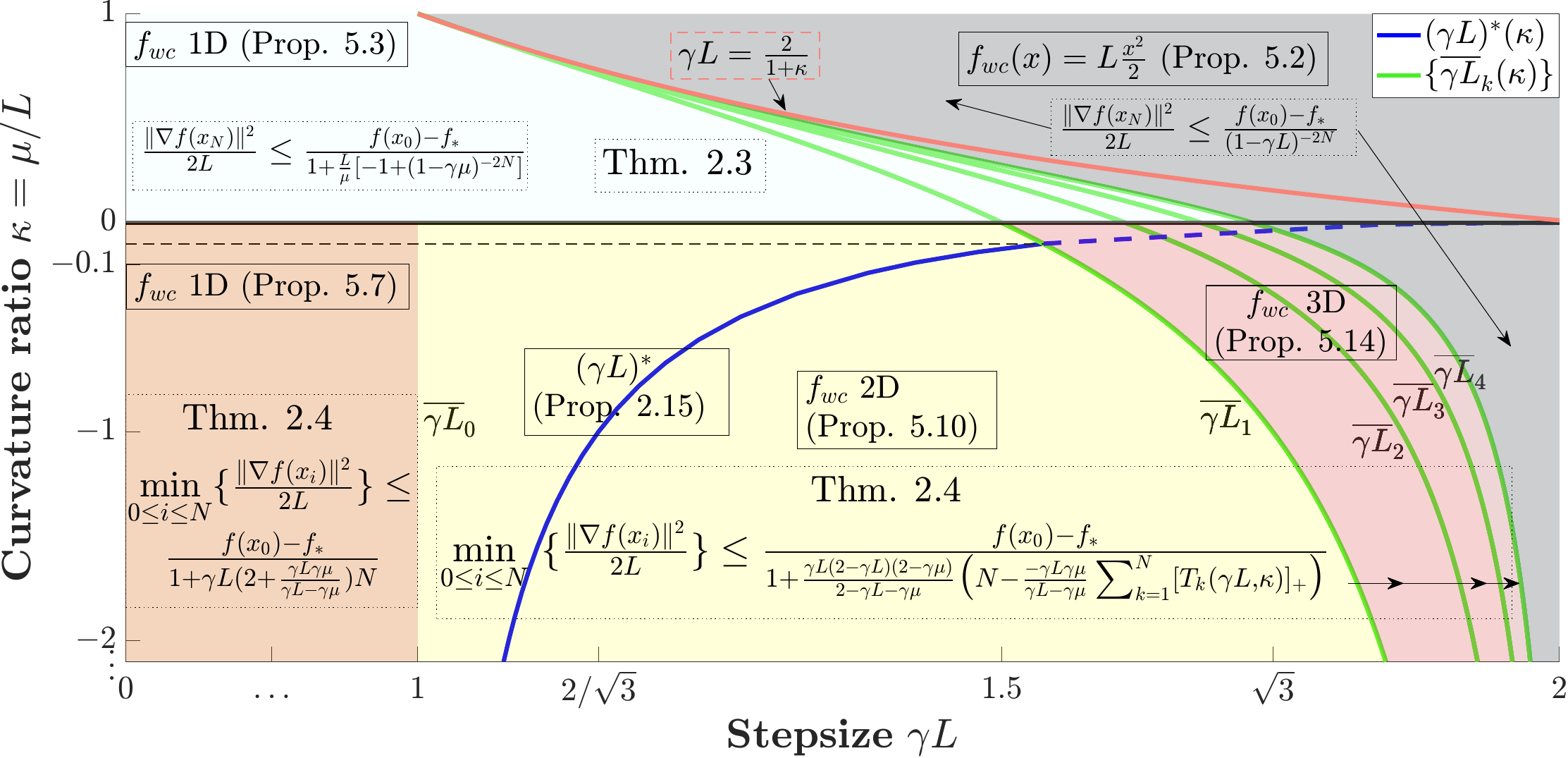}
    \caption{%
    Regimes outlined in \cref{corollary:cvx_rate} (smooth convex), \cref{thm:strg_cvx_rate} (smooth strongly convex) and \cref{thm:wc_GM_hypo_ct_all} (smooth weakly convex) with respect to the constant stepsize $\gamma L$ and curvature ratio $\kappa$, for $N=4$ iterations. State-of-the-art \textit{proven} tight convergence rates are limited to the cases $\kappa{=}-\infty$ \cite{Nesterov_cvx_2018}, $\kappa{=}-1$ \cite{abbaszadehpeivasti2021GM_smooth} and $\kappa = 0$ \cite{drori_performance_2014,Teboulle2022} and stepsizes lower than $\gbari[1]$. %
    Above this threshold, a series of $N$ regimes parameterized by $k=1,\dots,N$ exist and they are sublinear for non-strongly convex functions. 
    The optimal constant stepsize $\gstar[(\kappa)]$ for weakly convex functions is independent of $N$ for $\kappa \lesssim -0.1$, otherwise it is only asymptotically valid (as marked with {\color{blue}{\textit{blue}}} dashed line). 
    We use inequalities connecting three consecutive iterations to prove the regimes corresponding to stepsizes $\gamma L \in (\gbari[1][], \frac{2}{1+{[\kappa]_{+}}})$. The worst-case function ($f_{wc}$) dimensionality and the tightness' proof references are illustrated.}%
    \label{fig:Full_picture_all_regimes}
\end{figure}

The gist in proving tight rates on the full interval of stepsizes is its partition into a collection of intervals delimited by stepsize thresholds dependent on the curvature ratio $\kappa$, each of them requiring distinct analysis. \cref{fig:Full_picture_all_regimes} illustrates this partitioning for curvature ratio and (normalized) stepsizes, along with their corresponding rates for $N=4$ iterations. %

The remainder of the paper is organized as follows. \cref{sec:Results_wc_analysis} contains the formal statements of our main theoretical findings. \cref{subsec:Preliminaries} provides preliminary results required for the proofs. The complete proofs for the upper convergence bounds are detailed in \cref{sec:Proofs}, while \cref{sec:tightness} presents the worst-case examples that establish their tightness. %
To facilitate reproducibility, we provide a supplementary GitHub repository%
{}\footnote{\href{https://github.com/teo2605/GD_tight_rates}{GitHub repository:}\label{footnote:github_repo_interp_triplets} \url{https://github.com/teo2605/GD_tight_rates}} containing \textsc{Matlab} scripts for numerical verification and symbolic validation of the primary algebraic manipulations.%
%
%\end{document}
\section{Performance bounds}\label{sec:Results_wc_analysis}
    %!TEX root = ../../ms.tex
%
In \cref{subsec:rates_cvx,subsec:rates_strg_cvx,subsec:rates_hypoconvex} we provide the exact performance bounds for constant stepsize schedules. 
\cref{table:summary_bounds_and_tightness} summarizes the different results involved in their proofs. 
In \cref{subsec:rates_hypoconvex_variable_step_sizes} we give the tight bounds for smooth not-necessarily convex functions when using variable stepsizes, upper bounded by some threshold. 
\cref{subsec:stepsize_thresholds} describes the challenging stepsize thresholds defining the regimes for constant stepsizes. %
\cref{subsec:opt_gamma} provides optimal stepsize rules with respect to the theoretical worst-case scenarios: 
best constant stepsize policies in \cref{sec:optimal_constant_stepsizes} 
and a dynamic stepsize schedule with better worst-case guarantees in \cref{sec:dynamic_stepsizes}. %
\begin{table}[!ht]
  \caption[width=\textwidth]{
  Overview of performance bounds proofs and their tightness for all regimes with constant stepsizes (the thresholds $\gbari[k][]$, formally defined in \cref{subsec:stepsize_thresholds}, increase with $k$). 
  The demonstrations use interpolation inequalities connecting consecutive iterations (distance-1) or, additionally, the ones situated at distance-2. 
  Tightness is showed by constructing a worst-case function or by proving its existence via interpolating triplets.}%
  \label{table:summary_bounds_and_tightness}
  \centering
  \setlength{\arrayrulewidth}{0.2mm}
  \setlength{\tabcolsep}{5pt}
  \def\arraystretch{1.6}
  % \hspace{-1.5cm}
  \resizebox{\linewidth}{!}{
  \begin{tabular}{@{}c@{}|ccccc@{}}
    \toprule
    \makecell{\textbf{Class}} & \makecell{\textbf{Stepsizes}} & \textbf{Key result}  & 
    \makecell{\textbf{Inequalities} }& \makecell{\textbf{Tightness} \\ \textbf{proof}}  & 
    \makecell{
    \textbf{Worst-case} \\ \textbf{example}
    }\\
    %%%%%%%%%%%%%%%%%%%%%%%%%%%%%%%%%%%%%%%%%%%%%%%%%%%%%%%%%%%%%%%%%%%%%%%%
    %%%%%%%%%%%%%%%%%%%%%%%%%%%%%%%%%%%%%%%%%%%%%%%%%%%%%%%%%%%%%%%%%%%%%%%%
    %%%%%%%%%%%%%%%%%%%%%%%%%%%%%%%%%%%%%%%%%%%%%%%%%%%%%%%%%%%%%%%%%%%%%%%%    
    \midrule
    \multirow{3}{*}{\makecell{convex \\ $\mu=0$ \\ \cref{corollary:cvx_rate}} } 
    & $(0, \frac{3}{2}]$ & \makecell{\cref{lemma:sufficient_decrease_h_leq_1} \& \\ \cref{lemma:sufficient_decrease_h_geq_1}} & distance-1 & 
    \multirow{2}{*}{\makecell{\cref{corr:convex_tightness}} } & 
    \multirow{3}{*}{\makecell{1D function}} \\
    & $(\frac{3}{2}, \gbari[N][])$ & \cref{lemma:G4SD_sufficient_decrease_cvx} & distance-1$\&$2 & & \\
     & $[\gbari[N][], 2)$ & \cref{lemma:general_lemma_ct_h_to_Linear} & distance-1$\&$2 & \multirow{1}{*}{\makecell{\cref{prop:tightness_linear_regime}} } & \\
    %%%%%%%%%%%%%%%%%%%%%%%%%%%%%%%%%%%%%%%%%%%%%%%%%%%%%%%%%%%%%%%%%%%%%%%%
    %%%%%%%%%%%%%%%%%%%%%%%%%%%%%%%%%%%%%%%%%%%%%%%%%%%%%%%%%%%%%%%%%%%%%%%%
    %%%%%%%%%%%%%%%%%%%%%%%%%%%%%%%%%%%%%%%%%%%%%%%%%%%%%%%%%%%%%%%%%%%%%%%%
    \midrule
    \multirow{5}{*}{\makecell{strongly \\ convex \\ $\mu>0$ \\ \cref{thm:strg_cvx_rate}} }
    & $(0, \gbari[1][]]$ & \makecell{\cref{lemma:dist_1_mu_1_step_decrease}} & distance-1 & \multirow{3}{*}{\makecell{\cref{prop:Tightness_strongly_convex}} } & \multirow{5}{*}{\makecell{1D function} } \\
    %%%%%%%%%%%%%%%%%%%%%%%%%%%%%%%%%%%%%%%%%%%%%%%%%%%%%%%%%%%%%%%%%%%%%%%%
    & $[\gbari[1][], \gbari[N-1][])$ &  \makecell{\cref{lemma:dist_1_mu_1_step_decrease} \& \\ 
    \cref{lemma:general_lemma_ct_h_to_Linear_mu} } & distance-1$\&$2 & &  \\
    %%%%%%%%%%%%%%%%%%%%%%%%%%%%%%%%%%%%%%%%%%%%%%%%%%%%%%%%%%%%%%%%%%%%%%%%
    & $[\gbari[N-1][], \gbari[N][])$ &  \makecell{\cref{lemma:general_lemma_ct_h_to_Linear_mu} } & distance-1$\&$2 & &  \\
    %%%%%%%%%%%%%%%%%%%%%%%%%%%%%%%%%%%%%%%%%%%%%%%%%%%%%%%%%%%%%%%%%%%%%%%%
    & $[\gbari[N][], \frac{2}{1+\kappa})$ &  \makecell{\cref{lemma:general_lemma_ct_h_to_Linear}  } & distance-1$\&$2 & \multirow{2}{*}{\makecell{\cref{prop:tightness_linear_regime}} } &  \\
    %%%%%%%%%%%%%%%%%%%%%%%%%%%%%%%%%%%%%%%%%%%%%%%%%%%%%%%%%%%%%%%%%%%%%%%%
    & $[\frac{2}{1+\kappa}, 2)$ &  \makecell{\cref{lemma:dist_1_L_1_step_decrease} } & distance-1 & & 
    \\
    %%%%%%%%%%%%%%%%%%%%%%%%%%%%%%%%%%%%%%%%%%%%%%%%%%%%%%%%%%%%%%%%%%%%%%%%
    %%%%%%%%%%%%%%%%%%%%%%%%%%%%%%%%%%%%%%%%%%%%%%%%%%%%%%%%%%%%%%%%%%%%%%%%
    %%%%%%%%%%%%%%%%%%%%%%%%%%%%%%%%%%%%%%%%%%%%%%%%%%%%%%%%%%%%%%%%%%%%%%%%
    \midrule
     \multirow{5}{*}{\makecell{weakly convex\\ $\mu<0$ \\ \cref{thm:wc_GM_hypo_ct_all}}
     } & $(0,1]$& \makecell{\cref{lemma:sufficient_decrease_h_leq_1}} & distance-1 & \makecell{\cref{prop:tightness_short_steps}}   & 1D function  \\
    %%%%%%%%%%%%%%%%%%%%%%%%%%%%%%%%%%%%%%%%%%%%%%%%%%%%%%%%%%%%%%%%%%%%%%%%
    & $[ 1, \gbari[1][]{]}$& \makecell{\cref{lemma:sufficient_decrease_h_geq_1}}  & distance-1  & \makecell{\cref{prop:tightness_long_ct_steps}}  & 2D triplets  \\
    %%%%%%%%%%%%%%%%%%%%%%%%%%%%%%%%%%%%%%%%%%%%%%%%%%%%%%%%%%%%%%%%%%%%%%%%
    &  $[ \gbari[1][], \gbari[N-1][])$  & 
    \makecell{\cref{lemma:sufficient_decrease_h_geq_1} \&\\ 
            \cref{lemma:hypo_GN4SD}} & distance-1$\&$2 & \makecell{\cref{conj:3D_triplets}} & 3D triplets \\
    %%%%%%%%%%%%%%%%%%%%%%%%%%%%%%%%%%%%%%%%%%%%%%%%%%%%%%%%%%%%%%%%%%%%%%%%
    &  $[ \gbari[N-1][], \gbari[N][])$& \makecell{\cref{lemma:hypo_GN4SD}} & distance-1$\&$2  & \makecell{\cref{prop:2D_pure_quad}} & 2D function  \\    
    %%%%%%%%%%%%%%%%%%%%%%%%%%%%%%%%%%%%%%%%%%%%%%%%%%%%%%%%%%%%%%%%%%%%%%%%
    &  $[ \gbari[N][], 2 )$& \makecell{\cref{lemma:general_lemma_ct_h_to_Linear}} & distance-1$\&$2 & \makecell{\cref{prop:tightness_linear_regime}} & 1D function  \\ 
    %%%%%%%%%%%%%%%% Variable step-sizes  %%%%%%%%%%%%%%%%%%%%%%%%%%%%%%%%%%%%%%%%%%%%%%%%%%%%%%%%
    \bottomrule
\end{tabular} %
  }%
\end{table}%

The convergence rates are influenced by: 
(i) curvature ratio $\kappa{=}\frac{\mu}{L}$, 
(ii) the (normalized) fixed stepsizes $\gamma_{i} L \in (0,2)$, and 
(iii) the number of iterations $N$. Given a \textit{fixed budget} (number of iterations $N$), the rates show stepsize dependent accuracies in finding a stationary point, of the type
\begin{align}\label{eq:general_performance_bounds_type}\tag{$\mathcal{P}_*$}
   \frac{1}{2L}\min_{0 {} \leq {} i {} \leq {} N} \big\{ \|\nabla f(x_i)\|^2\big\}  
        {} \leq {}  
    \frac{f(x_0)-f_*}
    { 1 + \tau_{N}(\kappa, (\gamma_i L)_{i=0}^{N-1})},
\end{align}
where $\tau_{N}(\cdot)$  is some function we determine through our analysis. 
\begin{remark}\label{remark:rates_gap_f*_fN}    
    To be consistent with standard performance bounds in the literature, our final theorems use the optimality gap $f(x_0)-f_*$, as in \eqref{eq:general_performance_bounds_type}. 
    We derive these results by providing a bound on the \shownewtextfinal{gap} $f(x_0) - f(x_N)$, which is a measure directly available from the first-order oracle:
    \begin{align}\label{eq:general_performance_bounds_type_fN}\tag{$\mathcal{P}_N$}
   \frac{1}{2L}\min_{0 {} \leq {} i {} \leq {} N} \big\{ \|\nabla f(x_i)\|^2\big\}  
        {} \leq {}  
    \frac{f(x_0)-f(x_N)}
    { \tau_{N}(\kappa, (\gamma_i L)_{i=0}^{N-1})}.%
\end{align}%
The two bounds differ only by a ``+1'' term in the denominator. 
In nonconvex settings, e.g., multimodal functions, possessing multiple local optima, \eqref{eq:general_performance_bounds_type_fN} can be more effective for characterizing local convergence as it avoids the distortion of the numerator from a distant global optimum $f_*$ with a much lower value of $f_*$.%
\end{remark}
%
%!TEX root = ../../ms.tex
\subsection{Performance bounds for convex functions}\label{subsec:rates_cvx}
In \cref{corollary:cvx_rate} we establish the exact performance bound for smooth convex functions. This results arises as a limiting case of \cref{thm:strg_cvx_rate} and \cref{thm:wc_GM_hypo_ct_all}. Since the smooth convex setting is the most commonly encountered regime in analyses of gradient descent, we present it as a theorem.%
\begin{theorem}%[Exact worst-case rate for \textbf{convex} functions] 
\label{corollary:cvx_rate}
Let $f \in \mathcal{F}_{0, L}$ and consider $N$ iterations of \eqref{eq:GM_it_fixed_steps} with constant stepsizes $\gamma L \in (0,2)$ starting from $x_0$. Then the following bound holds:%
\begin{align}\label{eq:cvx_full_ct}
\begin{aligned}
   \frac{1}{2L} \|\nabla f(x_N)\|^2
        {} \leq {}  
    \frac{f(x_0)-f_*}
    { 1 + \gamma L \min \big\{
    {2 N} 
    \,,\,
    \frac{-1+(1-\gamma L)^{-2N}}{\gamma L}
         \big\}}.%
\end{aligned}
\end{align}%
\end{theorem}%
We prove \cref{corollary:cvx_rate} in \cref{proof:corollary:cvx_rate}. In \cref{corr:convex_tightness} we provide a worst-case function attaining the performance bound \eqref{eq:cvx_full_ct}. To our knowledge, this result establishes the first \textit{proved} tight rate in the convex case for the full range of constant stepsizes. 
The denominator in \eqref{eq:cvx_full_ct} is similar to the numerically conjectured one from \cite[Conjecture 3.1]{drori_performance_2014}, even though for a different performance metric, with $x_* \in \argmin_x f(x)$:%
\begin{align}\label{eq:conjectured_rate_cvx_fvals_vs_dist_drori_teboulle_2014}
    f(x_N) - f_*
        {}\leq{} 
    \frac{L}{2}
    \frac{\|x_0-x_*\|^2}{1 + \gamma L \min \big \{ 2 N \,,\, \frac{-1+(1-\gamma L)^{-2N}}{\gamma L} \big \} }.%
\end{align}%
Since the denominator is the same, both metrics are optimized by the same optimal constant stepsize (see \cref{thm:optimal_ct_stepsize_and_rate_convex}).%
%
%
%!TEX root = ../../ms.tex
\subsection{Performance bounds for strongly convex functions}%
\label{subsec:rates_strg_cvx}%
In \cref{thm:strg_cvx_rate} we derive the exact performance bound for smooth strongly convex functions. %
\begin{theorem}\label{thm:strg_cvx_rate}%
Let $f \in \mathcal{F}_{\mu, L}$, with $\mu \in (0, L)$, and consider $N$ iterations of \eqref{eq:GM_it_fixed_steps} with stepsizes $\gamma L \in (0, 2)$ starting from $x_0$. Then the following bound holds:%
\begin{align}\label{eq:strg_cvx_full_ct}
\begin{aligned}
    \frac{1}{2L} \|\nabla f(x_N)\|^2
        {} \leq {}  
    \frac{f(x_0)-f_*}
    { 1 + \gamma L \min \big\{
    \frac{-1+(1-\gamma \mu)^{-2N}}{\gamma \mu}
    \,,\,
    \frac{-1+(1-\gamma L)^{-2N}}{\gamma L}
         \big\}}.
\end{aligned}
\end{align}%
\end{theorem}%
We prove \cref{thm:strg_cvx_rate} in \cref{proof:thm:strg_cvx_rate}. In \cref{prop:Tightness_strongly_convex}, we provide a worst-case function that reaches the performance bound \eqref{eq:strg_cvx_full_ct}. %
% We show the tightness of the bound \eqref{eq:strg_cvx_full_ct} with a constructive example in \cref{prop:Tightness_strongly_convex}.

\shownewtextfinal{%
    Tight performance bounds for stepsizes below $\frac{2}{L+\mu}$ are given in \cite[Table 2]{Taylor_2017_SIAM_Composite_convex}\footnote{The work in \cite{Taylor_2017_SIAM_Composite_convex} provides tight performance bounds for the proximal gradient method, which generalizes gradient descent.}. For the metric analysed in this paper, however, only the rate corresponding to the stepsize $\gamma = \frac{1}{L}$ is conjectured there. In \cite{Taylor_2017_SIAM_Composite_convex}, the tight contraction factor $\max\{\gamma L - 1, 1 - \gamma \mu\}$ is derived and is optimized by choosing $\gamma = \frac{2}{L+\mu}$. Notably, this stepsize choice is suboptimal for the metric considered in our work (see \cref{thm:optimal_ct_stepsize_and_rate_strongly_convex}) and is also outperformed by the dynamic stepsize procedure proposed in \cref{sec:dynamic_stepsizes}. 
    
    Subsequent }%
to the publication of our work on \textit{arXiv}, \cite{kim2025proofexactconvergencerate} provided an alternative proof of \cref{thm:strg_cvx_rate}. In the same paper, exploiting the \textit{H-duality} phenomenon \cite{kim2023timereversed}, the author provides the first proof for the numerically conjectured rate from \cite[Conjecture 2]{taylor_smooth_2017}, which in our notation writes as:%
\begin{align}\label{eq:conjectured_rate_strg_cvx_fvals_vs_dist_taylor_2017}
    f(x_N) - f_*
        {}\leq{} &
    \frac{L}{2}
    \frac{\|x_0-x_*\|^2}{1 + \gamma L \min \big\{
    \frac{-1+(1-\gamma \mu)^{-2N}}{\gamma \mu}
    \,,\,
    \frac{-1+(1-\gamma L)^{-2N}}{\gamma L}
         \big\}},
\end{align}%
where $x_* = \argmin_x f(x)$. Since the denominator is the same, both metrics are optimized by the same optimal constant stepsize (see \cref{thm:optimal_ct_stepsize_and_rate_strongly_convex}). 

Letting $\mu \searrow 0$, the rate for convex functions from \cref{corollary:cvx_rate} is recovered.%
%
\subsection{Performance bounds for weakly convex functions}\label{subsec:rates_hypoconvex}%
In this section we derive and discuss performance bounds for smooth weakly convex functions, when using constant stepsizes. We begin by stating the bound in \cref{thm:wc_GM_hypo_ct_all}.%
\begin{theorem}\label{thm:wc_GM_hypo_ct_all}
Let $f \in \mathcal{F}_{\mu, L}$, with $\mu<0$, and consider $N$ iterations of \eqref{eq:GM_it_fixed_steps} with constant stepsizes $\gamma L \in (0, 2)$ starting from $x_0$. Then the following bound holds:
\begin{align} \label{eq:GM_hypo_ct_rates_f_star}
\hspace{-.5em}
    %%%%%%%%%%%%%%%%%%%%%%%%%%%%%%%%%%%
    \frac{1}{2L}
    \min_{0 {} \leq {} i {} \leq {} N} \big\{ \|\nabla f(x_i)\|^2\big\}  {} \leq {} 
        \frac{f(x_0)-f_*}
        {1 + \gamma L \min \big\{ \PN \,,\, \frac{-1+(1-\gamma L)^{-2N}}{\gamma L}  \big\} },
\end{align}%
with $P_N(l,u)$ defined as:%
\begin{align} \label{eq:P_den_full_rate_h_ct}
% \hspace{-0.5em}
P_N(l,u) {}{}\coloneqq{}{}
     p(l, u) \bigg[N - 
    \frac{-l u}{l-u}
    \sum\limits_{k=1}^N
        \Big[ \frac{-1+(1-u)^{-2k}}{u} - \frac{-1+(1-l)^{-2k}}{l} \Big]_{+}
    \bigg]
\end{align}%
where \\%
\vspace{-3em}
\begin{align}\label{eq:pi_hi_rate_GM}
    p(l,u)
        {}{}\coloneqq{}{}
        2 - \frac{-l u}{1-u - |1-l|}
        {}={}
    \left\{
    \arraycolsep=5pt
    \def\arraystretch{1.75}
        \begin{array}{ll}
            2 - \frac{-l u}{l-u}, &  l \in (0,1]; \\
%%%%%%%%%%%%%%%%%%%%%%%%%%%%%%%%%%%%%%%%%%%%%%%%%%%%%%%%%%%%%%%%%%%%%%%%%%%%%%%%%%%%%%%%
            \frac{(2-l)(2-u)}{2-l - u}, & l \in [1, 2).
        \end{array}
    \right.
\end{align}%
\end{theorem}%
We prove \cref{thm:wc_GM_hypo_ct_all} in \cref{proof:thm:wc_GM_hypo_ct_all}. 
 We then establish the tightness of these rates over the entire stepsize interval, which we partition into several subintervals. Explicit worst-case functions are provided 
for $\gamma L \in (0{,1}]$ (one-dimensional, \cref{prop:tightness_short_steps}), 
$\gamma L \in [\gbari[N-1]{,\gbari[N]})$ (two-dimensional, \cref{prop:2D_pure_quad}) and 
$\gamma L \in [\gbari[N]{, 2})$ (one-dimensional, \cref{prop:tightness_linear_regime}). For the remaining intervals, the worst-case behaviour is captured by interpolating triplets: 
$\gamma L \in (1{, \gbari[1]}{]}$ (two-dimensional, \cref{prop:tightness_long_ct_steps}) and 
$\gamma L \in [\gbari[1]{, \gbari[N]})$ (three-dimensional, \cref{conj:3D_triplets}).%

The rate's denominator combines a set of $N$ sublinear regimes with a linear regime also observed for convex and strongly convex functions (\cref{corollary:cvx_rate} and \cref{thm:strg_cvx_rate}), holding for stepsizes close to $2$. For $\gamma L \in (0, 1]$, each term in the sum appearing in \eqref{eq:P_den_full_rate_h_ct} is equal to zero. %
With $\gamma L \in (1, 2)$, only the terms up to some index $k$ become positive, while the remaining $(N-k)$ terms are negative. The expression is continuous with respect to the stepsize. A new term in the sum becomes nonnegative when evaluated at specific stepsize thresholds $\gbari[k][]$ rigorously defined in \cref{subsec:stepsize_thresholds} and satisfying the condition:%
\begin{align}\label{eq:condition_gbar}
    \frac{-1+(1-\kappa \gbari[k][])^{-2k}}{\kappa \gbari[k][]} - \frac{-1+(1-\gbari[k][])^{-2k}}{\gbari[k][]} {}={} 0.%
\end{align}%
Using the identity
$$
    {\sum\limits_{j{=}1}^{k}}
        \frac{{-1}{+}(1{-}u)^{-2j}}{u} {-} \frac{{-1}+(1{-}l)^{-2j}}{l}
        {}{=}{}
    \frac{{-1} {+} (1-u)^{-2k}}{1{-}(1{-}u)^2} 
    {-} \frac{{-1} {+} (1{-}l)^{-2k}}{1{-}(1{-}l)^2} +
    \big(\frac{1}{l} {-} \frac{1}{u}\big) {k},%
$$
the expression of $\PN$ from \eqref{eq:P_den_full_rate_h_ct} can be equivalently written as %
$$ \PN 
        {}={}
    p(\gamma L, \gamma \mu) \min\limits_{0 \leq k \leq N}
    \bigg\{         
    \frac{
    \frac{-1+(1-\gamma L)^{-2k}}{\gamma L [1-(1-\gamma L)^2]} 
    - \frac{-1+(1-\gamma \mu)^{-2k}}{\gamma \mu [1-(1-\gamma \mu)^2]}}{ \frac{1}{\gamma L} - \frac{1}{\gamma \mu} } + N - k
    \bigg\}.%
$$
For stepsize thresholds satisfying \eqref{eq:condition_gbar}, it reduces to summing an exponential term with a linear one:%
$$
    P_N(\gbari[k][], \kappa \gbari[k][])
        {}={}
    \frac{(2-\gbari[k][])(2-\kappa \gbari[k][])}{2-\gbari[k][] - \kappa \gbari[k][]} (N - k) {}+{} 
    \frac{-1+(1-\gbari[k][])^{-2k}}{\gbari[k][]}.
$$%
When reaching the maximum threshold $\gbari[N][]$ corresponding to $k=N$, $\PN$ becomes $\frac{-1+(1 - \gbari[N][])^{-2N}}{\gbari[N][]}$ and the rate shifts to the regime linear in $(1-\gbari[N][])$. Summarizing, the expression of $\PN$ can be expanded as follows:%
$$
    \left\{%
    \arraycolsep=4pt
    \def\arraystretch{1.8}
    \begin{array}{ll}
        \bigl( 2 + \frac{\gamma \mu \gamma L}{\gamma L- \gamma \mu}\bigr) N, & \gamma L \in (0, 1]; \\
        %%%%%%%%%%%%%%%%%%%%%%%%%%%%%%%%%%%%%%%%%%%%%%%%%%%%%%%%%%%%%%        
        \frac{(2-\gamma L)(2-\gamma \mu)}{2-\gamma L - \gamma \mu} (N - k) {}+{} 
        \dfrac{ \frac{-1+(1-\gamma L)^{-2k}}{\gamma L [1-(1-\gamma L)^2]} - \frac{-1+(1-\gamma \mu)^{-2k}}{\gamma \mu [1-(1-\gamma \mu)^2]}}{ \frac{1}{1-(1-\gamma L)^2} - \frac{1}{1-(1-\gamma \mu)^2} } ,
        & \makecell{
        \gamma L \in \big[\gbari[k][], \gbari[k+1][] \big)\\
        k=0,{\dots},N-1;
        } \\
        %%%%%%%%%%%%%%%%%%%%%%%%%%%%%%%%%%%%%%%%%%%%%%%%%%%%%%%%
        \frac{-1 + (1-\gamma L)^{ -2 N }}{\gamma L},  &  \gamma L = \gbari[N][].
    \end{array}
    \right.
$$%
By fixing stepsize $\gamma L$ and varying $N$, another perspective on convergence rates emerges: determining the \shownewtextfinal{maximum number} of iterations to reach a desired accuracy with a \textit{given stepsize}. Let $\bar{N}$ be the largest index $k$ corresponding to a strictly positive term in expression of $\PN$ from \eqref{eq:P_den_full_rate_h_ct}. Then for $N \leq \bar{N}$ the rate is linear in $(1-\gamma L)$, whereas for $N \geq \bar{N}+1$ it shifts to a sublinear regime. %

%
\subsection{Performance bounds for non-strongly convex functions and variable stepsizes}\label{subsec:rates_hypoconvex_variable_step_sizes}
\begin{figure}
    \centering
        \centering
        \includegraphics[width=\textwidth]{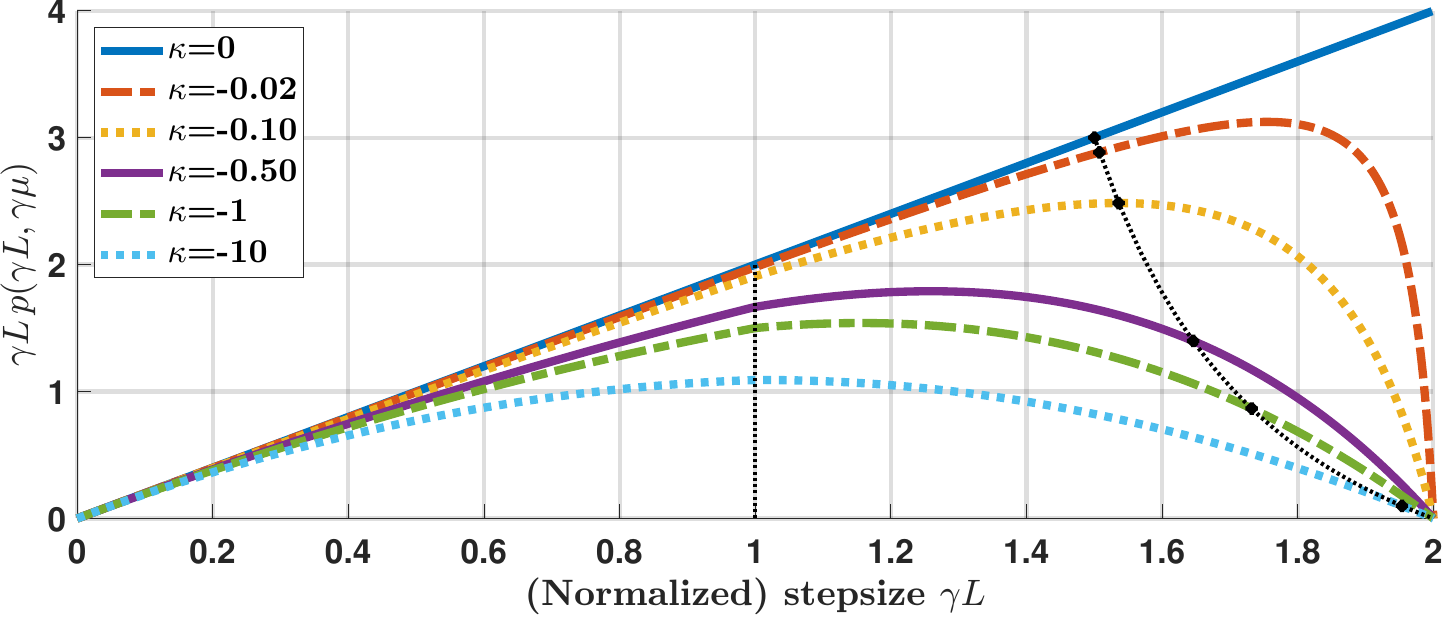}
        \caption{The dominant term $\gamma L \, \pigLgmu$ (see \eqref{eq:pi_hi_rate_GM}) in the bounds for weakly convex functions for several curvature ratios. $\gamma L = 1$ delimits the two main sublinear regimes, while the one-step analysis is tight up to the threshold $\gbari[1]$, marked in black dots, above which intervene transient contributions (see \eqref{eq:P_den_full_rate_h_ct}).}%
        \label{fig:constant_func_denominator}
    \hfill
\end{figure}

%%%%%%%%%%%%%%%%%%% One-step Theorem (up to threshold) %%%%%%%%%%%%%%%%%%%%%%%%%%%
%
In this section, we derive and discuss performance bounds for smooth weakly convex and convex functions under variable (but predefined) stepsizes, up to a well-defined threshold. In \cref{thm:wc_GM_hypo}, we establish the corresponding guarantee for $\mu \leq 0$, and \cref{corr:cvx_rates_variable_stepsizes} specializes this guarantee to the convex setting (i.e., $\mu=0$).%
\begin{theorem} \label{thm:wc_GM_hypo}
Let $f\in \mathcal{F}_{\mu, L}$, with $\mu \in (-\infty, 0]$, 
and consider $N$ iterations of \eqref{eq:GM_it_fixed_steps} starting from $x_0$ with stepsizes $\gamma_i L \in (0,\gbari[1]]$, $i=0,\dots,N-1$, where $\gbari[1] {}{}\coloneqq{}{}\tfrac{3}{1+\kappa+\sqrt{1-\kappa+\kappa^2}} \in [\tfrac{3}{2}, 2)$. Then the following bound holds, where $p(\gamma L,\gamma \mu)$ is defined in \eqref{eq:pi_hi_rate_GM}:%
\begin{align}\label{eq:GM_hypo_rate_both_p_and_q}
    \frac{1}{2L} \min_{0 {} \leq {} i {} \leq {} N} \big\{ \|\nabla f(x_i)\|^2\big\}
        {} \leq {} 
    \frac{f(x_0)-f_*}
        {1+\sum_{i=0}^{N-1} \gamma_i L \, \pigLgmu[\gamma_i]}.
\end{align}%
\end{theorem}%
\cref{thm:wc_GM_hypo} is demonstrated in \cref{proof:thm:wc_GM_hypo}. The bound in \eqref{eq:GM_hypo_rate_both_p_and_q} is proved to be tight for $\gamma_i L \in (0, 1]$ in \cref{prop:tightness_short_steps}. For the range $\gamma_i L \in (1, \gbari[1]]$, we conjecture its tightness based on a two-dimensional worst-case example (\cref{conjecture:variable_mid_stepsizes}).%
\begin{remark}
     Restricting it to constant stepsizes $\gamma L \leq \gbari[1]$, \cref{thm:wc_GM_hypo} recovers the result from \cref{thm:wc_GM_hypo_ct_all} corresponding to 
     $\frac{-1+(1-\gamma \mu)^{-2k}}{\gamma \mu} - \frac{-1+(1-\gamma L)^{-2k}}{\gamma L} {}\leq{} 0$ for all $k \geq 1$ in \eqref{eq:P_den_full_rate_h_ct}, 
     namely $\PN {}={} p(\gamma L,\gamma \mu)\, N$. The extension to variable stepsizes relies on the proof which, within the stepsize range from \cref{thm:wc_GM_hypo}, only requires inequalities relating consecutive iterations.
\end{remark}%
\cref{fig:constant_func_denominator} illustrates the behavior of the leading term $\gamma L \, \pigLgmu$, emphasizing the stepsize threshold $\gbari[1]$ and additionally showing that optimal rates, i.e., with larger denominators, are achieved in the regime of stepsizes greater than $1$.%
\begin{remark}
\cref{thm:wc_GM_hypo} extends to arbitrary $\mu < 0$ the state-of-the-art rate for smooth not-necessarily convex functions ($\mu = -L$) from \cite[Theorem 2]{abbaszadehpeivasti2021GM_smooth}:%
\begin{align}\label{eq:GM_nonconvex_smooth}
    \frac{1}{2L}
    \min_{0 {} \leq {} i {} \leq {} N} \big\{ \|\nabla f(x_i)\|^2\big\}
    {} \leq {} 
        \frac{f(x_0)-f_*}{1 + 
        \sum_{i=0}^{N-1} \gamma_i L \big(2 - \frac{\gamma_i L}{2} 
        \max\{1, \gamma_i L\}
        \big)}, \, \forall \gamma_i L \in (0, \sqrt{3}].
\end{align}%
The convergence rate \eqref{eq:Nesterov_rate}, stated, for instance, in \cite[Section 1.2.3]{Nesterov_cvx_2018}, is derived for smooth functions using only the upper curvature $L$:%
\begin{align} \label{eq:Nesterov_rate}
\frac{1}{2L}
    \min_{0 {} \leq {} i {} \leq {} N} \big\{ \|\nabla f(x_i)\|^2\big\}
{}\leq {} 
    \frac{f(x_0)-f_*}{1 + 
    \sum_{i=0}^{N-1} \gamma_i L (2 - \gamma_i L )}, {}\qquad{} \gamma_i L \in (0, 2).%
\end{align}%
\cref{thm:wc_GM_hypo_ct_all} recovers this result in the absence of lower curvature information, i.e., for $\mu \searrow -\infty$ and effectively becomes \cref{thm:wc_GM_hypo} covering the full stepsize domain $(0, 2)$ (as $\gbari[1][(-\infty)] = 2$). The proof of this celebrated result lacks the utilization of quadratic \textit{lower} bounds and thus remains generally non-tight. The convergence rate from \cref{thm:wc_GM_hypo_ct_all} exhibits a smooth interpolation between the result in \eqref{eq:Nesterov_rate} (when $\mu \searrow -\infty$), the rate for smooth not-necessarily convex functions ($\mu=-L$) from \eqref{eq:GM_nonconvex_smooth} and the rate on smooth convex functions ($\mu=0$).%
\end{remark}%
\begin{corollary}\label{corr:cvx_rates_variable_stepsizes}
    Let $f \in \mathcal{F}_{0, L}$ and consider $N$ iterations of \eqref{eq:GM_it_fixed_steps} starting from $x_0$ with fixed stepsizes $\gamma_i L \in (0, \tfrac{3}{2}]$, $i=0,\dots,N-1$. Then the following bound holds:%
\[
       \frac{1}{2L} \|\nabla f(x_N)\|^2
            {} \leq {} 
        \frac{f(x_0)-f_*}{1 + 
        2\sum_{i=0}^{N-1} \gamma_i L}.%
\]%
\end{corollary}%
\cref{corr:cvx_rates_variable_stepsizes} is obtained by taking $\mu=0$ in \cref{thm:wc_GM_hypo}. Although not explicitly presented as such, Teboulle and Vaisbourd propose in \cite{Teboulle2022} two sufficient decrease results (Lemma 1 and Lemma 5) yielding the same rate from \cref{corr:cvx_rates_variable_stepsizes} through telescoping summation (see also \cref{remark:descent_lemmas_Teboulle_22}). %
\subsection{Stepsize thresholds}\label{subsec:stepsize_thresholds}%
\begin{figure}[!htb]
    \centering
    \begin{subfigure}[t]{0.49\textwidth}
        \centering
        \includegraphics[width=\linewidth]{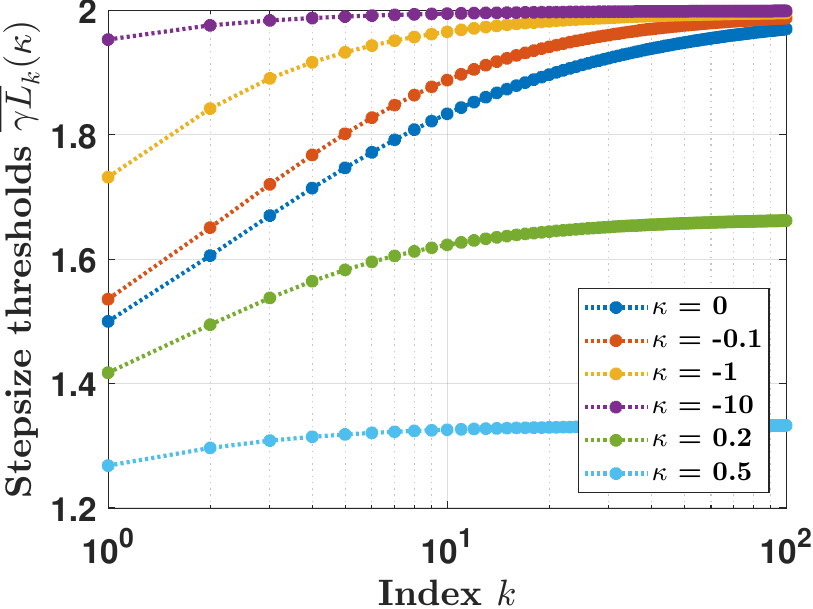}
        \caption{Stepsize thresholds $\gbari[k]$ converge to $\frac{2}{1+\kappa}$ for strongly convex functions and to $2$ for non-strongly convex functions.}
        \label{fig:stepsize_thresholds}
    \end{subfigure}\hfill
    \begin{subfigure}[t]{0.49\textwidth}
        \centering
        \includegraphics[width=\linewidth]{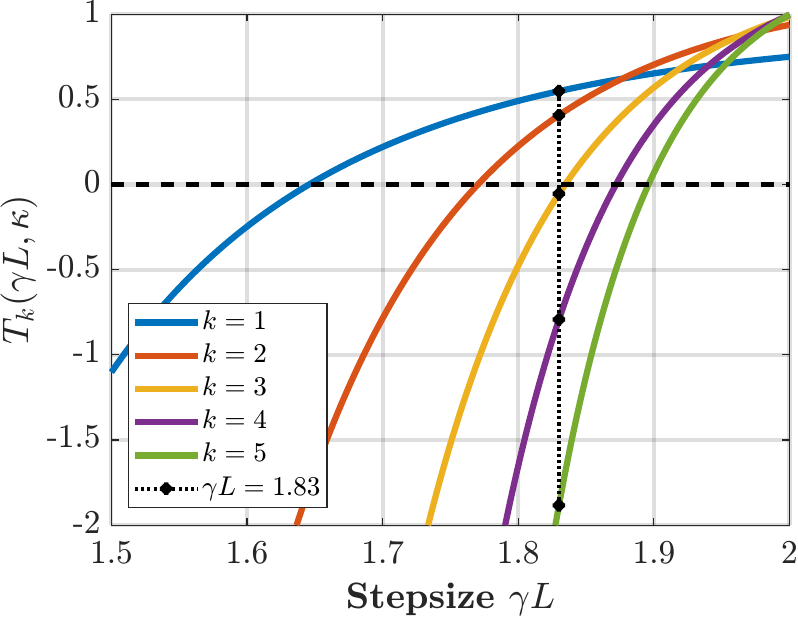}
        \caption{Dependence of $\TigLmu[k]$ on $\gamma L$ ($\kappa = -0.5$). For $\gamma L = 1.83$, $T_k \geq 0$ for $k\in\{1,2\}$; $\gbari[k]$ mark the x-intercept.}
        \label{fig:Plot_Tis_vs_h}
    \end{subfigure}
    \caption{\shownewtextfinal{Stepsize thresholds and the sign transition of $T_k$ as a function of the normalized stepsize $\gamma L$, for fixed curvature ratio $\kappa$.}}
    \label{fig:stepsize_and_Tk}
\end{figure}

Stepsize thresholds are given by the roots of \eqref{eq:condition_gbar}. For $N$ iterations, there are exactly $N$ such thresholds, delimiting the regimes for $\gamma L > 1$. These thresholds play a central role in the proofs for (strongly) convex and weakly convex functions in \cref{corollary:cvx_rate,thm:strg_cvx_rate,thm:wc_GM_hypo_ct_all}, as they partition the stepsize domain into $N+1$ distinct regimes, each corresponding to a different worst-case instance. We begin by introducing the following auxiliary expressions.%
\begin{definition}\label{def:def_TN}% 
Let $k$ be a positive integer and define function $E_k:\mathbb{R}\backslash\{0\} \rightarrow \mathbb{R}$:%
\begin{align}
    \Ei{}{}\coloneqq{}{}\sum\limits_{j=1}^{2k} x^{-j} 
        {}={} 
     \left\{
    \arraycolsep=6pt
    \def\arraystretch{1.25}
    \begin{array}{cc}
      \frac{-1+x^{-2k}}{1-x}, & x \neq 1; \\
      2k, & x = 1.
    \end{array}
    \right. \label{eq:def_Ei} \tag{$E_k$}
\end{align}%
Let $V \coloneqq \{(\ell, \kappa): \, \kappa \in (-\infty, 1), \,\ell \in (1, \frac{2}{1+[\kappa]_{+}})\}$ and define $T_k: V \rightarrow \mathbb{R}$ as %
\begin{align}
    T_k(\ell,\kappa)
        {}{}\coloneqq{}{}&
    \Ei[k][1-\kappa \ell] - \Ei[k][1-\ell] \label{eq:def_TN} \tag{$T_k$}.
\end{align}
\end{definition}%
$\TigLmu[k]$ denotes exactly the individual terms summed up in \eqref{eq:P_den_full_rate_h_ct} and canceled by the $k$-th stepsize threshold in \eqref{eq:condition_gbar}. Moreover, $\TigLmu[N]$ is the difference of the arguments in the denominators of the rates for convex functions \eqref{eq:cvx_full_ct} and strongly convex functions \eqref{eq:strg_cvx_full_ct}, respectively.
\begin{proposition}\label{prop:monotonicity_Tk}%
    The function $\TigLmu[k][\cdot][\kappa]$ is strictly increasing in the first argument.%
\end{proposition}%
\begin{proof}%
See \cref{appendix:proof_prop_monotonicity_Tk}.%
\qed\end{proof}%
The sign of $\TigLmu$ changes in the stepsize interval $\gamma L \in \big(1,\frac{2}{1+[\kappa]_{+}}\big)$, on which it therefore possesses a unique root.%
\begin{definition}[Stepsize thresholds]\label{def:stepsize_thresholds}%
Let $\gbari[\infty]{}\coloneqq{} \frac{2}{1+[\kappa]_{+}}$. The unique roots of $\TigLmu$ in the interval $(1,\gbari[\infty])$ are referred to as stepsize thresholds and denoted by $\gbari[k]$, where $\gbari[0] {}{}\coloneqq{}{} 1$ and:%
\begin{align}\label{def:gamma_bar} \tag{$\xoverline{\gamma L}_k$}
    \big\{\gbari[k]\big\} 
        {}{}\coloneqq{}{} 
    \big\{ \gamma L \in (1,\gbari[\infty]) \,|\, \TigLmu[k] = 0 \big\} {}\quad{} \forall k = 1,2,{\dots}\,\,.%
\end{align}%
\end{definition}%
The definition of stepsize thresholds \eqref{def:gamma_bar} resembles exactly the condition \eqref{eq:condition_gbar} marking the transitions between the different regimes indexed by $k$ in \cref{thm:wc_GM_hypo_ct_all}.%
\begin{proposition}\label{prop:properties_step_size_thr}
    The stepsize thresholds $\gbari$ satisfy the following properties:%
    \begin{enumerate}[align=right,itemindent=1em,labelsep=2pt,labelwidth=1em,leftmargin=0pt,nosep]
        \item $\gbari[k] < \gbari[k+1]$, for all integers $k \geq 0$;
        \item $\lim_{k\rightarrow \infty} \gbari[k] = \gbari[\infty]$;
        \item \shownewtextfinal{For any $\gamma L \in \big[\gbari[k],\gbari[k+1]\big)$: 
        $\TigLmu[i][\gamma L] \geq 0$ $\forall i \leq k$, 
        and $\TigLmu[i][\gamma L] < 0$ $\forall i \geq k+1$.}
     \end{enumerate}%
\end{proposition}%
\begin{proof}%
See \cref{appendix:proof_properties_step_size_thr}.%
\qed\end{proof}%
\noindent \cref{fig:Full_picture_all_regimes} illustrates the dependence on $\kappa$ of the first $N=4$ stepsize thresholds. Complementary, \cref{fig:stepsize_thresholds} details the rapid increase of the thresholds with the number of iterations, reaching: (i) $\frac{2}{1+\kappa}$ for strongly convex functions and (ii) $2$ in the non-strongly convex case. 
The distance between two consecutive thresholds decreases when increasing their index. 
\cref{fig:Plot_Tis_vs_h} depicts $\TigLmu$ defined in \eqref{eq:def_TN} for $\kappa=-0.5$ and several indices $k$, offering insight on how to determine the stepsize thresholds. For example, with $\gamma L =1.83$, $\TigLmu[k]$ is positive for $k \leq 2$ and negative for $k \geq 3$; thus, the interval corresponding to these conditions is $[\gbari[2][], \gbari[3][])$.%
\subsection{Improved stepsize schedules}\label{subsec:opt_gamma}%
Building on our worst-case performance analysis, we provide several stepsize policies that minimize these upper bounds. %
% Recommendations are provided for both constant (\cref{sec:optimal_constant_stepsizes}) and variable stepsizes (\cref{sec:dynamic_stepsizes}). %
Recommendations are provided for both constant stepsizes (\cref{sec:optimal_constant_stepsizes}), which depend on the iteration count, and variable stepsizes (\cref{sec:dynamic_stepsizes}); the latter improve upon the constant choices while also remaining independent of the horizon.%
\subsubsection{Optimal constant stepsizes}\label{sec:optimal_constant_stepsizes}%
In this section, we derive the optimal constant stepsize that maximizes the worst-case performance established in \cref{sec:Results_wc_analysis}. \cref{thm:optimal_ct_stepsize_and_rate_convex} and \cref{thm:optimal_ct_stepsize_and_rate_strongly_convex} provide this stepsize for convex and strongly convex functions, respectively. \cref{prop:gamma_star} derives the (asymptotically) optimal constant stepsize for weakly convex functions.%
\begin{proposition}[Constant stepsize for convex functions]\label{thm:optimal_ct_stepsize_and_rate_convex}%
    Let $f \in \mathcal{F}_{0, L}$ and consider $N$ iterations of $\eqref{eq:GM_it_fixed_steps}$ starting from $x_0$ with constant stepsizes. The worst-case rate for convex functions \eqref{eq:cvx_full_ct} is minimized by the choice $\gamma L = \gbari[N][(0)]$, implying the best worst-case \shownewtextfinal{bound}:%
    \begin{align}\label{eq:cvx_opt_ct}
\begin{aligned}
    \frac{1}{2L} \|\nabla f(x_N)\|^2
        {} \leq {}  
    \frac{f(x_{0})-f_*}{1 + 2N  \, \gbari[N][(0)]}.%
\end{aligned}
\end{align}%
\end{proposition}%
\begin{proposition}[Constant stepsize for strongly convex functions]\label{thm:optimal_ct_stepsize_and_rate_strongly_convex}%
    Let $f \in \mathcal{F}_{\mu, L}$, with $\mu \in (0, L)$, and $N$ iterations of $\eqref{eq:GM_it_fixed_steps}$ starting from $x_0$ with constant stepsizes. The worst-case rates for strongly convex functions \eqref{eq:strg_cvx_full_ct} is minimized by the choice $\gamma L = \gbari[N]$, implying the best worst-case \shownewtextfinal{bound}:%
    \begin{align}\label{eq:strg_cvx_opt_ct}
\begin{aligned}
    \frac{1}{2L} \|\nabla f(x_N)\|^2
        {} \leq {}  
    \frac{f(x_{0})-f_*}{(1-\gbari[N])^{-2N}}.%
\end{aligned}
\end{align}%
\end{proposition}%
To the limit $N \rightarrow \infty$, $\gbari[N]$ becomes $\frac{2}{1+\kappa}$, yielding a rate with $(\frac{1-\kappa}{1+\kappa})^{2N}$ behaviour. %
When evaluated at the $N$-th stepsize threshold $\gbari[N]$, \cref{def:stepsize_thresholds} implies that the two arguments in the minimization defining the performance bounds of \cref{corollary:cvx_rate} and \cref{thm:strg_cvx_rate} coincide.

The same optimal constant stepsize $\gbari[N]$ is minimizing the worst-case rate for the more standard performance criterion $\frac{f(x_N)-f_*}{\|x_0-x_*\|^2}$, since the denominators are the same (see \cref{corollary:cvx_rate} for convex functions 
and \cref{thm:strg_cvx_rate} for strongly convex functions). We refer the reader to the discussion following \cite[Conjecture 3.1]{drori_performance_2014} for convex functions and \cite[Sections 4.1.1-2]{taylor_smooth_2017} for strongly convex functions. These rates were recently proved in \cite{kim2025proofexactconvergencerate}.%
\begin{proposition}[Constant stepsize for weakly convex functions]\label{prop:gamma_star}
Let $\kappa < 0$, $\bar{\kappa} \coloneqq \tfrac{-9-5\sqrt{5} + \sqrt{190+90\sqrt{5}}}{4} \approx -0.1001$, and define $\gstar[(\kappa)]$ as the unique solution belonging to the interval $[1,2)$ of the equation:%
    \begin{align}
    \begin{aligned} \label{eq:equation_solve_for_h_star_one_step}
        -\kappa(1+\kappa) (\gamma L)^3 +
    \big[3 \kappa + (1+\kappa)^2\big] (\gamma L)^2 -
    4(1+\kappa) (\gamma L) +
    4 = 0.%
    \end{aligned}
    \end{align}%
With respect to minimizing the tight performance bound for weakly convex functions from \cref{thm:wc_GM_hypo_ct_all}, $\gstar[(\kappa)]$ is:%
\begin{enumerate}%
    \item the optimal constant stepsize for any $N$, when $\kappa \leq \bar{\kappa}$;%
    \item the \textit{asymptotically} optimal constant stepsize as $N \rightarrow \infty$, i.e., when the algorithm runs for a sufficiently large number of iterations, if $\kappa \in (\bar{\kappa},0)$.%
\end{enumerate}%   
\end{proposition}%
\begin{proof}
    See \cref{appendix:optimal_constant_stepsize_nonconvex}. \qed%
\end{proof}%
Unlike the optimal stepsize $\tfrac{2}{\sqrt{3} L}$ provided in \cite[Theorem 3]{abbaszadehpeivasti2021GM_smooth} for the usual smooth not-necessarily convex case $\mu = -L$, the recommendation in \cref{prop:gamma_star} leverages lower curvature information and maximizes the leading term $\gamma L \, \pigLgmu$ in the rate's denominator from \cref{thm:wc_GM_hypo_ct_all}. As $\kappa$ approaches $-\infty$, implying the absence of lower curvature information, the quantity $\gamma L (2-\gamma L)$ is maximized, resulting in the celebrated optimal stepsize $\tfrac{1}{L}$. Similarly, substituting $\kappa=-1$ in equation \eqref{eq:equation_solve_for_h_star_one_step} yields the optimal stepsize for smooth not-necessarily convex functions from \cite{abbaszadehpeivasti2021GM_smooth}. %
\begin{figure}[!htb]
    \centering
    \begin{subfigure}[t]{0.49\textwidth}
        \centering
        \includegraphics[width=\linewidth]{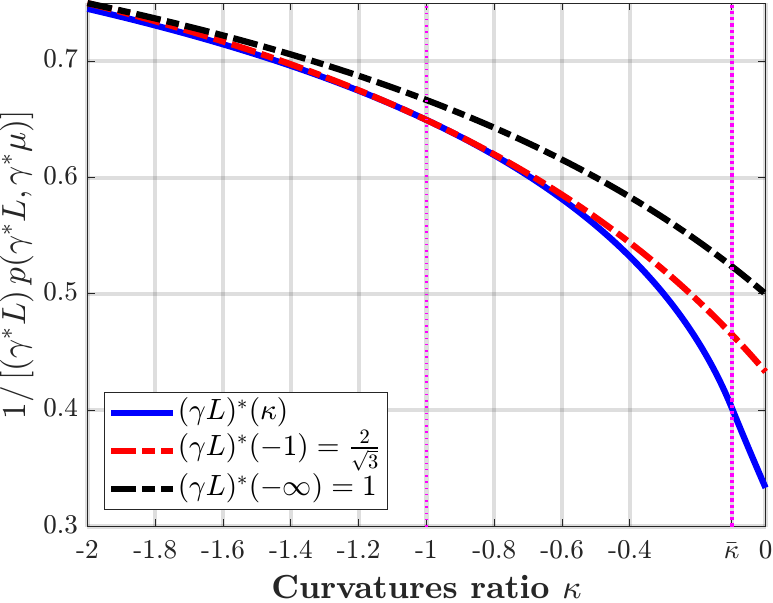}
        \caption{Comparison of (normalized) optimal stepsize recommendations: (i) $\gstar[(\kappa)]$ from \cref{prop:gamma_star}, (ii) $\gamma L=\tfrac{2}{\sqrt{3}}$ from \cite{abbaszadehpeivasti2021GM_smooth}, and (iii) $\gamma L=1$ (classical).}
        \label{fig:Compare_opt_ct}
    \end{subfigure}\hfill
    \begin{subfigure}[t]{0.49\textwidth}
        \centering
        \includegraphics[width=\linewidth]{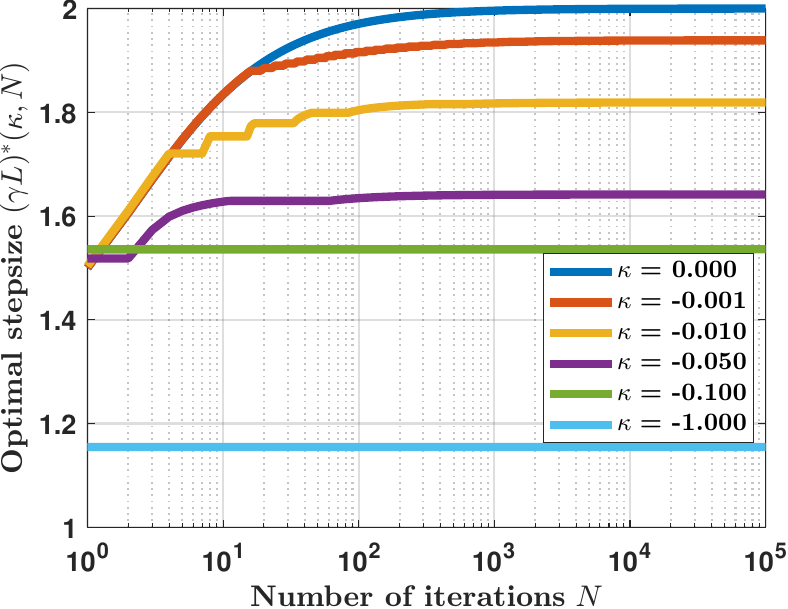}
        \caption{Optimal stepsize versus $N$ (log scale): asymptotic convergence to $\gstar[(\kappa)]$; near convexity, a transient staircase behavior appears.}
        \label{fig:Asymptotic_behaviour_vs_N}
    \end{subfigure}
    \caption{\shownewtextfinal{Optimal constant stepsize selection for weakly convex functions.}}
    \label{fig:optimal_stepsize_two_panels}
\end{figure}%

\cref{fig:Full_picture_all_regimes} depicts the continuous dependence of $\gstar[(\kappa)]$ on the curvature ratio; above $\bar{\kappa} \approx -0.1$, the dashed line marks its asymptotic optimality. %
\cref{fig:Compare_opt_ct} illustrates the maximum guaranteed improvement achievable by optimizing the leading term $p^{-1}(\gamma L, \gamma \mu)$ and compares: (i) the (asymptotically) optimal stepsize $\gstar[(\kappa)]$ from \cref{prop:gamma_star}, (ii) the optimal stepsize for smooth functions $\gstar[(-1)]=\tfrac{2}{\sqrt{3}}$, and (iii) the celebrated optimal stepsize $\gstar[(-\infty)] = 1$. As nonconvexity decreases ($\kappa \nearrow 0$), the optimal stepsize from \cref{prop:gamma_star} exhibits a greater improvement compared to stepsize $\gstar[(-1)]$, which does not employ lower curvature information.%

In contrast to the constant optimal stepsize recommendations for (strongly) convex functions, in the weakly convex case there exists an asymptotically optimal schedule independent of the number of iterations. A transient regime appears only for $-0.1 \lesssim \kappa < 0$, where the optimal constant stepsize also varies with $N$, i.e., $(\gamma L)^* = (\gamma L)^* (\kappa,N)$, increasing monotonically. \cref{fig:Asymptotic_behaviour_vs_N} shows this dependency for several curvature ratios, where the worst-case is numerically minimized for a given iteration budget. As $\kappa \nearrow 0$, the transient regime to reach the asymptotically optimal stepsize from \cref{prop:gamma_star} expands. During the initial iterations, a staircase-like behavior is observed and the optimal stepsize is confined within a specific subdomain determined by the stepsize thresholds, i.e., $(\gamma L)^*(\kappa,N) \in [1, \gbari[N])$. In the convex case, the staircase behavior disappears, as the optimal constant stepsize $(\gamma L)^*(0,N)$ becomes $\gbari[N][(0)]$, lying at the intersection of two linear regimes (see \cref{thm:optimal_ct_stepsize_and_rate_convex}). %
\subsubsection{Dynamic stepsizes}\label{sec:dynamic_stepsizes}%
Stepsizes with $\gamma L > 1$ admit faster rates, motivating analysis in this regime. For (strongly) convex functions, variable stepsize policies can attain superior worst-case performance to constant ones by occasionally exceeding $2$ \cite{altschuler2023acceleration_cvx,altschuler2023acceleration_strg_cvx,grimmer2024accelerated_compact,zhang2024acceleratedgradientdescentconcatenation,zhang2024anytimeaccelerationgradientdescent}, but they rely on intricate proofs, may depend on a fixed horizon, and their non-monotone nature makes them sensitive to curvature misestimation. %

We introduce a monotonically increasing stepsize sequence inspired by \cref{thm:wc_GM_hypo} and the work of Teboulle-Vaisbourd \cite{Teboulle2022}, that avoids these issues. It applies to smooth strongly convex, convex, and weakly convex functions, is horizon-free, achieves better worst-case guarantees than the optimal constant stepsize, and admits significantly simpler convergence proofs. %
The construction is derived by telescoping a two-point descent inequality (\cref{lemma:sufficient_decrease_h_geq_1}), yielding a recurrence that cancels a key term and produces the largest admissible stepsize at each iteration. %
\begin{definition}[Dynamic stepsizes] \label{def:variable_stepsizes}
Let $\kappa \in (-\infty, 1)$. We consider the sequence $\left\{s_k(\kappa)\right\}_{k=-1}^{\infty}$, with $s_{-1}=0$, recursively defined as:%
\begin{align}\label{eq:def_sequence_stepsizes}%\tag{$s_k$}
    s_{k+1}(\kappa)
        {}{}{\coloneqq}{}{}
    \big\{
        s_{+} \in \big(1 , \tfrac{2}{1+[\kappa]_{+}}\big){:}%
        \frac{s_{+} [ (2-s_{+})(2-\kappa s_{+}) - 1 ] }{2-s_{+}(1+\kappa)}
            {+}
        \frac{s_{k}(\kappa)}{2-s_{k}(\kappa)(1+\kappa)} = 0
    \big\}.
\end{align}
\end{definition}
\begin{proposition}\label{prop:properties_variable_stepsizes_sequence}
For any fixed $\kappa \in (-\infty, 1)$, the sequence $\left(s_k(\kappa)\right)_{k=-1}^{\infty}$ is: %
\begin{enumerate}[
  align=right,
  label=(\roman*),
  ref=(\roman*),
  itemindent=2em,
  labelsep=2pt,
  labelwidth=1em,
  leftmargin=0pt,
  nosep
]
\item uniquely defined; 
\item monotonically increasing; and 
\item $\lim_{k\rightarrow \infty} s_k(\kappa) = \frac{2}{1+[\kappa]_{+}}$.
\end{enumerate}
Moreover, for any $k \geq 1$ it holds that
\begin{enumerate}[
  resume,
  align=right,
  label=(\roman*),
  itemindent=2em,
  labelsep=2pt,
  labelwidth=1em,
  leftmargin=0pt,
  nosep
]
\item\label{eq:dynamic_prop_cvx} $s_{k-1}(0) {}\geq{} 2 - \frac{1}{\frac{3}{2}k}$; 
\item\label{eq:dynamic_prop_strg_cvx} $s_{k-1}(\kappa) {}\geq{} \frac{2}{1+\kappa} - (\frac{1-\kappa}{1+\kappa})^{2k}$, if $\kappa \in (0,1)$. 
\end{enumerate}
% (i) uniquely defined; 
% (ii) monotonically increasing; and 
% (iii) $\lim_{k\rightarrow \infty} s_k(\kappa) = \frac{2}{1+[\kappa]_{+}}$. Moreover, for any $k \geq 1$ it holds that
% (iv) $s_{k-1}(0) {}\geq{} 2 - \frac{1}{\frac{3}{2}k}$; 
% (v) $s_{k-1}(\kappa) {}\geq{} \frac{2}{1+\kappa} - (\frac{1-\kappa}{1+\kappa})^{2k}$, if $\kappa \in (0,1)$.
\end{proposition}%
\begin{proof}
    See \cref{proof:prop:properties_variable_stepsizes_sequence}. \qed
\end{proof}%
Note that $s_0=\gbari[1]$ is the largest value for which \cref{thm:wc_GM_hypo} holds. 
Obtaining a closed-form expression for $s_k$ is difficult, as each step requires solving a quadratic equation (in the convex case) or cubic one (when $\kappa \neq 0$).
%%%%%%%%%%%%%%%%%%%%%%%%%%%%%%%%%%%%%%%%%%%%%%%%%%%%%%%%%%%%%%%%%%%%%%%%%%%%%%%%%%%%%%%%%%%%%%%%%%%%%%%%%%%
%%%%%%%%%%%%%%%%%%%%%%%%%%%%%%%%%%%%%%%%%%%%%%%%%%%%%%%%%%%%%%%%%%%%%%%%%%%%%%%%%%%%%%%%%%%%%%%%%%%%%%%%%%%
\begin{theorem}\label{prop:tight_rate_all_kappas_increasing_sequence}
Let $f \in \mathcal{F}_{\mu, L}$, with $\mu \in (-\infty, L)$, and $\kappa = \frac{\mu}{L}$. Consider $N$ iterations of \eqref{eq:GM_it_fixed_steps} starting from $x_0$ with stepsizes $\gamma_i L$ defined as follows:%
    $$
        \gamma_i L {}\coloneqq{} \left\{ \begin{array}{cc}
           \min \{ s_i(\kappa) \,,\, \gstar[(\kappa)] \},  &  \kappa < 0 \\
           s_i(\kappa),  & \kappa \geq 0
        \end{array} \right., \quad \forall \, i = 0, {\dots}, N-1,
    $$
    where $\gstar[(\kappa)]$ is defined in \cref{prop:gamma_star} and $s_i(\kappa)$ is given in \cref{def:variable_stepsizes}. Then the following bound holds:%
\begin{align}\label{eq:tight_rate_all_kappas_increasing_sequence}
    \frac{1}{2L} \min_{0 {} \leq {} i {} \leq {} N} \big\{ \|\nabla f(x_i)\|^2\big\} {} \leq {} 
        \frac{f(x_0)-f_*}
        {1 + \sum_{i=0}^{N-1} \frac{\gamma_i L \left(2-\gamma_i L\right)\left(2-\kappa \gamma_i L\right)}{2-\gamma_i L (1+\kappa)}}.
\end{align}%
\end{theorem}%
\begin{proof}%
   See \cref{proof:prop:tight_rate_all_kappas_increasing_sequence}. \qed%
\end{proof}%
Based on numerical observations, we conjecture that the bound in \eqref{eq:tight_rate_all_kappas_increasing_sequence} is tight. %
\cref{prop:tight_rate_all_kappas_increasing_sequence} extends \cref{thm:wc_GM_hypo} to larger stepsizes (i.e., $>\gbari[1]$) and also applies to the strongly convex case ($\kappa > 0$). 
\cref{prop:tight_rate_all_kappas_increasing_sequence} is proved in \cref{proof:prop:tight_rate_all_kappas_increasing_sequence}, together with the resulting corollaries for convex, strongly convex and weakly convex functions (\cref{thm:convergence_rates_strg_convex_variable_stepsizes,corollary:convergence_rates_convex_variable_stepsizes,thm:GM_hypo_rate_sequence}).%
%%%%%%%%%%%%%%%%%%%%%%%%%%%%%%%%%%%%%%%%%%%%%%%%%%%%%%%%%%%%%%%%%%%%%%%%%%%%%%%%%%%%%%%%%%%%%%%%%%%%%%%%%%%
%%%%%%%%%%%%%%%%%%%%%%%%%%%%%%%%%%%%%%%%%%%%%%%%%%%%%%%%%%%%%%%%%%%%%%%%%%%%%%%%%%%%%%%%%%%%%%%%%%%%%%%%%%%
\begin{corollary}\label{corollary:convergence_rates_convex_variable_stepsizes}
    Let $f \in \mathcal{F}_{0, L}$  and consider $N$ iterations of \eqref{eq:GM_it_fixed_steps} starting from $x_0$ with stepsizes %
    $\gamma_i L = s_{i}(0) = 1 + \frac{2}{1+\sqrt{9 - 4 s_{i-1}(0)}}$
    for all $i=0,\dots,N-1$. Then the following bound holds:%
    \begin{align}\label{eq:wc_rate_var_stepsizes_convex}%\tag{wcvar}
\begin{aligned}
    \frac{1}{2L} \|\nabla f(x_{N})\|^{2}
        {}\leq{}
    \frac{f(x_{0})-f_*}{1 + \frac{s_{N-1}(0)}{2-s_{N-1}(0)}}.%
\end{aligned}%
\end{align}%
\end{corollary}%
According to \crefdefpart{prop:properties_variable_stepsizes_sequence}{eq:dynamic_prop_cvx}, the bound can be relaxed by $\frac{f(x_0)-f_*}{3N}$, which recovers the standard rate of $\mathcal{O}(\frac{1}{N})$.%
\begin{corollary}\label{thm:convergence_rates_strg_convex_variable_stepsizes}
    Let $f \in \mathcal{F}_{\mu, L}$, with $\mu \in (0, L)$, and consider $N$ iterations of \eqref{eq:GM_it_fixed_steps} starting from $x_0$ with stepsizes $\gamma_i L = s_i(\kappa)$ from \cref{def:variable_stepsizes}, with $i=0,\dots,N-1$. Then the following bound holds:%
    \begin{align}\label{eq:wc_rate_var_stepsizes_strg_convex}%\tag{wcvar}
\begin{aligned}
    \frac{1}{2L} \|\nabla f(x_{N})\|^{2}
        {}\leq{}
    \frac{f(x_{0})-f_*}{1 + \frac{s_{N-1}(\kappa)}{2-s_{N-1}(\kappa)(1+\kappa)}}.
\end{aligned}
\end{align}
\end{corollary}%
As a consequence of \crefdefpart{prop:properties_variable_stepsizes_sequence}{eq:dynamic_prop_strg_cvx}, the gradient norm converges to zero at a rate of  $\mathcal{O}( (\frac{1-\kappa}{1+\kappa})^{2N})$. %

Following the approach in \cite[Table 1]{Teboulle2022}, we numerically compare our dynamic stepsize policy with the optimal constant stepsize dependent on the number of iterations. 
This comparison is presented in \cref{table:compare_rates_ct_vs_dynamic_cvx} (convex case) and \cref{table:compare_rates_ct_vs_dynamic_kappa_0.001} (strongly convex case). 
Such a numerical comparison is necessary to evaluate the worst-case performance, as our dynamic policy lacks a closed-form expression.%
\begin{table}[!ht]
\centering
\caption{Comparison of complexity bounds (denominators) for convex functions between the policies: (i) celebrated stepsize $\gamma L = 1$; (ii) optimal constant stepsizes (\cref{thm:optimal_ct_stepsize_and_rate_convex}) and (iii) dynamic schedule (from \cref{corollary:convergence_rates_convex_variable_stepsizes}); higher is better. The last column displays the ratio of the these last two denominators.%
}\label{table:compare_rates_ct_vs_dynamic_cvx}%
\def\arraystretch{1.25}%
\resizebox{.85\textwidth}{!}{%
\begin{tabular}{@{}cccccccc@{}}
\toprule
\multirow{2}{*}{$N$} & 
\multicolumn{2}{c}{Standard choice} &
\multicolumn{2}{c}{\cref{thm:optimal_ct_stepsize_and_rate_convex}} & \multicolumn{2}{c}{\cref{corollary:convergence_rates_convex_variable_stepsizes}} &  \\ \cline{2-7}
    & $\gamma L = 1$ &   $1 + 2N$      & $\gbari[N][]$       & $1 + 2N \gbari[N][]$      & $s_{N-1}$  & $1+\frac{s_{N-1}}{2-s_{N-1}}$ & Ratio ($\%$)\TBstrut\\ \midrule
1     & 1 & 3             & 1.500               & 4.000           & 1.500      & 4.000                         & 100.000      \\
2     & 1 & 5             & 1.606               & 7.423           & 1.732      & 7.464                         & 100.549      \\
5     & 1 & 11             & 1.747               & 18.471          & 1.893      & 18.619                        & 100.806      \\
10    & 1 & 21             & 1.834               & 37.681          & 1.947      & 37.933                        & 100.670      \\
20    & 1 & 41             & 1.897               & 76.885          & 1.974      & 77.235                        & 100.456      \\
30    & 1 & 61             & 1.924               & 116.426         & 1.983      & 116.826                       & 100.343      \\
40    & 1 & 81             & 1.939               & 156.106         & 1.987      & 156.535                       & 100.275      \\
50    & 1 & 101            & 1.949               & 195.859         & 1.990      & 196.310                       & 100.230      \\
100   & 1 & 201            & 1.971               & 395.109         & 1.995      & 395.612                       & 100.127\Bstrut\\ \bottomrule 
\end{tabular}%
}%
\end{table}%
\begin{table}[!ht]
\centering
\caption{Comparison of complexity bounds (denominators) for strongly convex functions between the policies: (i) celebrated stepsize $\gamma L = \frac{2}{1+\kappa}$; (ii) optimal constant stepsizes (\cref{thm:optimal_ct_stepsize_and_rate_strongly_convex}) and (iii) dynamic schedule (\cref{thm:convergence_rates_strg_convex_variable_stepsizes}); higher is better. The last column displays the ratio of these last two denominators.}%
\label{table:compare_rates_ct_vs_dynamic_kappa_0.001}
\def\arraystretch{1.25}%
\resizebox{\textwidth}{!}{%
\begin{tabular}{c|cccccccc}
\toprule
 \multirow{2}{*}{$\kappa$} & \multirow{2}{*}{$N$} & \multicolumn{2}{c}{Standard choice} & \multicolumn{2}{c}{\cref{thm:optimal_ct_stepsize_and_rate_strongly_convex}} & \multicolumn{2}{c}{\cref{thm:convergence_rates_strg_convex_variable_stepsizes}}          &              \\ \cline{3-8}
              & & $\frac{2}{1+\kappa}$ &   $\big(\frac{1-\kappa}{1+\kappa}\big)^{-2N}$     & $\gbari[N][]$   & $(1-\gbari[N][])^{-2N}$  & $s_{N-1}$ & $1+\frac{s_{N-1}}{2-s_{N-1}(1+\kappa)}$ & Ratio ($\%$)\TBstrut\\ \midrule
%%%%%%%%%%%%%%%%%%%%%%%%%%%%%%%%%%%%%%%%%%%%%%%%%%%%%%%%%%%%%%%%%%%%%%%%%%%%%%%%%%%%%%%%%%%%
\multirow{9}{*}{$10^{-3}$} 
& 1 & 1.998 & 1.004 & 1.500 & 4.006   & 1.500 & 4.006   & 100.000 \\
& 2 & 1.998 & 1.008  & 1.605 & 7.447   & 1.731 & 7.488   & 100.551 \\
& 5 & 1.998 & 1.020 & 1.746 & 18.633  & 1.892 & 18.784  & 100.813 \\
& 10 & 1.998 & 1.041 & 1.833 & 38.381  & 1.946 & 38.643  & 100.682 \\
& 20 & 1.998 & 1.083  & 1.896 & 79.879  & 1.973 & 80.257  & 100.473 \\
& 30 & 1.998 & 1.127 & 1.923 & 123.416 & 1.982 & 123.865 & 100.363 \\
& 40 & 1.998 & 1.174 & 1.938 & 168.871 & 1.986 & 169.373 & 100.297 \\
& 50 & 1.998 & 1.221 & 1.948 & 216.257 & 1.989 & 216.805 & 100.253 \\
& 70 & 1.998 & 1.492 & 1.960 & 317.040 & 1.992 & 317.670 & 100.199\Bstrut\\   \midrule
%%%%%%%%%%%%%%%%%%%%%%%%%%%%%%%%%%%%%%%%%%%%%%%%%%%%%%%%%%%%%%%%%%%%%%%%%%%%%%%%%%%%%%%%%%%%
\multirow{9}{*}{$10^{-4}$} 
         & 1  & 1.9998 & 1.000   &     1.500     &       4.001     &      1.500      &       4.001    &   100.000\Tstrut \\
         & 2  & 1.9998 & 1.001    &     1.606     &       7.426     &      1.732      &       7.467    &   100.550 \\
         & 5  & 1.9998 & 1.002    &     1.747     &      18.487     &      1.893      &      18.636    &   100.807 \\
        & 10  & 1.9998 & 1.004    &     1.834     &      37.750     &      1.947      &      38.004    &   100.671 \\
        & 20  & 1.9998 & 1.008    &     1.897     &      77.178     &      1.974      &      77.531    &   100.457 \\
        & 30  & 1.9998 & 1.012    &     1.924     &     117.101     &      1.983      &     117.505    &   100.345 \\
        & 40  & 1.9998 & 1.016    &     1.939     &     157.323     &      1.987      &     157.759    &   100.277 \\
        & 50  & 1.9998 & 1.020    &     1.949     &     197.780     &      1.990      &     198.240    &   100.232 \\
        & 70   &       1.9998 &        1.028  &        1.961   &        279.309    &       1.993       &    279.801   &    100.176\Bstrut\\ \bottomrule
\end{tabular}%
}%
\end{table}%
\begin{corollary}\label{thm:GM_hypo_rate_sequence}
    Let $f\in \mathcal{F}_{\mu, L}$, with $\mu \in (-\infty, 0)$, and consider $N$ iterations of \eqref{eq:GM_it_fixed_steps} starting from $x_0$ with stepsizes $\gamma_i L {=} \min \{ s_i(\kappa) \,,\, \gstar[(\kappa)] \}$, with $i=0,{\dots},N-1$, $s_i(\kappa)$ defined in \cref{def:variable_stepsizes} and $\gstar[(\kappa)]$ given in \cref{prop:gamma_star}. Then the following bound holds:%
\begin{align}\label{eq:GM_hypo_rate_sequence}
    \frac{1}{2L} \min_{0 {} \leq {} i {} \leq {} N} \big\{ \|\nabla f(x_i)\|^2\big\}
        {} \leq {} 
    \frac{f(x_0)-f_*}
        {1 + p_* N + 
 \max \limits_{0 \leq k \leq N} \big\{ \frac{s_{k-1}(\kappa)}{2-s_{k-1}(\kappa)\,(1+\kappa)} - k p_* \big\}
         },%
\end{align}%
where $p_* {}\coloneqq{} p(\gstar[], \kappa \gstar[])$ is the maximized coefficient $\pigLgmu$ from \eqref{eq:pi_hi_rate_GM}.%
\end{corollary}%
\begin{table}[!ht]
\caption{Comparison of complexity bounds (denominators) for weakly convex functions (example $\kappa = -10^{-3}$) between the policies: (i) asymptotically constant optimal stepsize $\gstar[(\kappa)]$ from \cref{prop:gamma_star}, resulting in a denominator denoted by $P_N^A=P_N(\gstar[(\kappa)], \kappa \gstar[(\kappa)])$; (ii) constant optimal stepsize depending on the number of iterations $(\gamma L)^*(N,\kappa)$ determined by numerically maximizing the denominator $\PN$, denoted by $P_N^*=P_N\left((\gamma L)^*(N,\kappa), \kappa (\gamma L)^*(N,\kappa)\right)$; and (iii) dynamic sequence from \cref{thm:GM_hypo_rate_sequence}, leading to the denominator $D_N$. The transition between the increasing sequence and the constant schedule is marked by the horizontal line. The last column displays the ratio $D_N(\kappa) / P_N^* $.}
\label{table:nonconvex_policies_kappa_-0.001}
\centering
\def\arraystretch{1.25}%
\resizebox{\textwidth}{!}{%
\begin{tabular}{@{}cccccccc}
\hline
\multirow{2}{*}{N} & \multicolumn{2}{c}{\cref{prop:gamma_star}}               & \multicolumn{2}{c}{\cref{prop:gamma_star}}     & \multicolumn{2}{c}{\cref{thm:GM_hypo_rate_sequence}}  &              \\ \cline{2-7}
                   & $\gstar[(\kappa)]$ & $P_N^A$ & $(\gamma L)^*(N,\kappa)$ & $P_N^*$ & $\min\{s_{N-1},\gstar[(\kappa)]\}$ & $D_N(\kappa)$ & Ratio ($\%$)\Bstrut\\ \hline
1                  & 1.939              & 1.135                                            & 1.500                    & 3.994                                            & 1.500                              & 3.994         & 100.000      \\
2                  & 1.939              & 1.288                                            & 1.606                    & 7.400                                            & 1.733                              & 7.440         & 100.547      \\
5                  & 1.939              & 1.882                                            & 1.748                    & 18.310                                           & 1.893                              & 18.456        & 100.798      \\
8                  & 1.939              & 2.750                                            & 1.809                    & 29.509                                           & 1.935                              & 29.721        & 100.719      \\ \cline{6-8}
9                  & 1.939              & 3.121                                            & 1.823                    & 33.254                                           & 1.939                              & 33.483        & 100.689      \\
10                 & 1.939              & 3.542                                            & 1.835                    & 36.999                                           & 1.939                              & 37.246        & 100.667      \\
20                 & 1.939              & 12.546                                           & 1.885                    & 74.170                                           & 1.939                              & 74.867        & 100.940      \\
30                 & 1.939              & 44.438                                           & 1.894                    & 111.360                                          & 1.939                              & 112.489       & 101.014      \\
40                 & 1.939              & 144.899                                          & 1.898                    & 148.645                                          & 1.939                              & 150.111       & 100.986      \\
50                 & 1.939              & 182.520                                          & 1.905                    & 186.003                                          & 1.939                              & 187.733       & 100.930      \\
100                & 1.939              & 370.629                                          & 1.916                    & 373.255                                          & 1.939                              & 375.841       & 100.693      \\ \hline
\end{tabular}%
}%
\end{table}%
In the weakly convex case, the sequence $\left\{s_k(\kappa)\right\}_{k=0}^{\infty}$ is truncated when it exceeds the asymptotically optimal constant stepsize for nonconvex functions from \cref{prop:gamma_star}. %
The numerical evidence in \cref{table:nonconvex_policies_kappa_-0.001} supports that the rate from \cref{thm:GM_hypo_rate_sequence} outperforms both: (i) optimal constant stepsize schedule $(\gamma L)^*(N,\kappa)$ (obtained by numerically maximizing the denominator from \cref{thm:wc_GM_hypo_ct_all}) and (ii) asymptotically optimal constant stepsize $\gstar[(\kappa)]$. As $\kappa$ approaches $0$, $\gstar[(\kappa)]$ converges to $2$ and the result from \cref{corollary:convergence_rates_convex_variable_stepsizes} is retrieved. %
%
%
\section{Preliminary results for the proofs}\label{subsec:Preliminaries}
Interpolation conditions are a set of inequalities that tightly characterize a function belonging to a specific function class $\mathcal{F}$ by focusing on a set of triplets $\mathcal{T}=\{(x_i,g_i,f_i)\}_{i\in \mathcal{I}}$. The interpolation problem itself, described in \cite{rockafellar1970convex,lambert2004finite}, has its historical roots in Whitney's work on the extension of continuous functions \cite{whitney1934analytic,whitney1934differentiable}.
As demonstrated by Taylor et al. \cite{taylor_smooth_2017}, interpolation conditions are a key tool for precise convergence analysis required to \textit{exactly} solve the infinite-dimensional performance estimation problems (PEP) introduced in the seminal work of Drori and Teboulle \cite{drori_performance_2014}.
\begin{definition}[{\cite[Definition 2]{taylor_smooth_2017}}] The set $\mathcal{T}{}\coloneqq{}\{(x_i, g_i, f_i)\}_{i \in \mathcal{I}}$ is called $\mathcal{F}_{\mu, L}$-interpolable if and only if there exists a function $f \in \mathcal{F}_{\mu, L}$ such that $\nabla f(x_i)=g_i$ and $f(x_i)=f_i$ for all $i\in \mathcal{I}$.
\end{definition}%
\cref{thm:interp_hypo} generalizes the interpolation conditions for smooth (strongly) convex functions \cite[Theorem 4]{taylor_smooth_2017} and for smooth non-necessarily functions ($\mu=-L$) \cite[Theorem 3.10]{Taylor_2017_SIAM_Composite_convex} to any function where $\mu < 0$. The proof is based on standard arguments and is given in \cref{subsec:proof_interp_conds_hypo}. A graphical interpretation of interpolating functions for the case $\mu=-L$ is given in \cite[page 71]{PhD_AT_2017}.%
\begin{theorem}[Interpolation conditions for smooth functions]\label{thm:interp_hypo}
Set $\mathcal{T}=\big\{(x_i, g_i, f_i)\big\}_{i \in \mathcal{I}}$ is $\mathcal{F}_{\mu,L}$-interpolable, where $\mu \in (-\infty, L)$, if and only if for every pair of indices $(i,j)$, with $i,j \in \mathcal{I}$, it holds:%
\begin{align}\label{eq:Interp_hypoconvex}
f_{i}{-}f_{j}{-}\langle g_{j}, x_{i}-x_{j} \rangle 
        {}\geq{}
    \frac{1}{2L} \|g_{i}{-}g_{j}\|^{2} {}+{} \frac{\mu}{2L(L{-}\mu)} \|g_{i}{-}g_{j} {-} L (x_i {-}x_j)\|^2.
\end{align}%
\end{theorem}%
\cref{thm:interp_hypo} is tighter than \cite[Theorem 2.2]{Themelis_DR_2020}, which, to the best of our knowledge, is the first to explicitly use lower curvature information to improve stepsize ranges \shownewtextfinal{in weakly convex settings} (the authors use it in the context of deriving convergence rates for the Douglas-Rachford splitting method).%
\begin{corollary} \label{corr:direct_sum_interp_cond} 
Let $f \in \mathcal{F}_{\mu,L}$, with $\mu \in (-\infty, L)$. For all $x_i, x_j \in \mathbb{R}^d$ it holds:%
\begin{align}\label{eq:direct_sum_interp_cond}
    0 {}\geq{} \left\langle \nabla f(x_i) - \nabla f(x_j) - L(x_i-x_j) \,,\, 
                           \nabla f(x_i) - \nabla f(x_j) - \mu(x_i-x_j) \right\rangle.
\end{align}%
\end{corollary}%
\begin{proof}
    Inequality \eqref{eq:direct_sum_interp_cond} results by summing up inequalities \eqref{eq:Interp_hypoconvex} written for the pairs $(x_i,x_j)$ and $(x_j,x_i)$ and performing simplifications.
\qed\end{proof}
\cref{corr:direct_sum_interp_cond} is an extension to negative lower curvatures $\mu$ of the co-coercivity property from \cite[Theorem 2.1.12]{Nesterov_cvx_2018}. When $L + \mu > 0$, the class $\mathcal{F}_{\mu,L}$ fits in the framework of $(\frac{L \mu}{L+\mu}, \frac{1}{L+\mu})$-semimonotone operators (e.g., see \cite[Proposition 4.13]{evens2023_precond_PPA}), which satisfy inequality \eqref{eq:direct_sum_interp_cond} by definition.

Interpolation conditions are used in conjunction with the  \eqref{eq:GM_it_fixed_steps} iterations by substituting $x_{j}=x_{i} - \sum_{k=i}^{j-1} \gamma_k \nabla f(x_k)$, as detailed in \cref{def:dist-N_interp_cond}. Since the proofs are limited to constant stepsizes, we constrain $\gamma_k$ to $\gamma$.%
\begin{definition}\label{def:dist-N_interp_cond}
Inequalities \eqref{eq:distance-N_interpolation_conditions_I_ij} and \eqref{eq:distance-N_interpolation_conditions_I_ji} are necessary and sufficient conditions connecting gradient descent iterations separated by $(j-i)$ steps:
% %%%%%%%%%%%%%%%%%%%%%%% Improved decent version %%%%%%%%%%%%%%%%%%%%%%%%%%%%%%%%%%%%%%%%
\begin{align}\label{eq:distance-N_interpolation_conditions_I_ij}
    & \hspace{-.75em}
    \begin{aligned}
    f(x_{i})-f(x_{j})
        {}\geq {}  
    \gamma \big \langle \nabla f(x_{j}), \sSum_{k=i}^{j-1} \nabla f(x_{k}) \big \rangle +
    \frac{1}{2L} \left\|\nabla f(x_{i}) - \nabla f(x_{j})\right\|^2 + 
    \\% {}\qquad{}
    \frac{\mu}{2L(L{-}\mu)} \big\|\nabla f(x_{i}) - \nabla f(x_{j}) - \gamma L \sSum_{k=i}^{j-1} \nabla f(x_{k})\big\|^2; 
    \end{aligned} \tag{$\Iij$} \\
    %%%%%%%%%%%%%%%%%%%%%%%%%%%%%%%%%%%%%%%%%%%%%%%%%%%%%%%%%%%%
    & \hspace{-.75em}
    \begin{aligned} 
    f(x_{j})-f(x_{i}) 
        {}\geq {}
    {-}\gamma \big \langle \nabla f(x_{i}), \sSum_{k=i}^{j-1} \nabla f(x_{k}) \big \rangle +
    \frac{1}{2L} \left\|\nabla f(x_{i}) - \nabla f(x_{j})\right\|^2 + 
    \\% {}\qquad{}
    \frac{\mu}{2L(L{-}\mu)} \big\|\nabla f(x_{i}) - \nabla f(x_{j}) - \gamma L \sSum_{k=i}^{j-1} \nabla f(x_{k})\big\|^2.
    \end{aligned} \tag{$\Iji$}
    \label{eq:distance-N_interpolation_conditions_I_ji}
\end{align}%
We refer to inequalities \eqref{eq:distance-N_interpolation_conditions_I_ij} and \eqref{eq:distance-N_interpolation_conditions_I_ji} as \textit{distance-1} or \textit{distance-2} inequalities if the distances between their indices satisfy $|i-j|\leq 1$ or $|i-j|\leq 2$, respectively.%
\end{definition}%
In contrast to traditional convergence study, we separate the analysis concerning the optimal objective value $f_*$ and solely focus on the function value gap up to step $N$, namely $f(x_0) - f(x_N)$. By making a slight modification, as detailed in the descent result from \cref{lemma:sufficient_decrease_N_*}, we derive a convergence rate incorporating $f_*$ (for a proof of the lemma, see for example \cite[Section 1.2.3]{Nesterov_cvx_2018}).%
\begin{lemma}\label{lemma:sufficient_decrease_N_*}
Let $f \in \mathcal{F}_{\mu,L}$, with $\mu \in (-\infty, L)$. One iteration of \eqref{eq:GM_it_fixed_steps} with stepsize $\gamma L \in (0,2)$ decreases the function values and at each point $x_i$ it holds:
\begin{align}\label{eq:sufficient_decrease_N_*}
    f(x_i) - f_* {}\geq{} \frac{1}{2L} \|\nabla f(x_i)\|^2.
\end{align}%
\end{lemma}%
A convergence rate about the last gradient norm ${\|}\nabla f(x_N){\|}$ is obtained for (strongly) convex functions, since the gradient norm is non-increasing, as shown for example in \cite[Lemma 2]{Teboulle2022}. %
\begin{lemma}[Monotonicity of gradient norm for convex functions] \label{lemma:monotonicity_gradient_cvx_case}
Let $f \in \mathcal{F}_{0,L}$ and one iteration with stepsize \shownewtext{$\gamma >0$}, connecting $x_{i}$ and $x_{i+1}$. Then:%
\begin{align}\label{eq:gradient_norm_decrease_convex_case}
    \|\nabla f(x_{i})\|^2 - \|\nabla f(x_{i+1})\|^2 \geq \tfrac{2-\gamma L}{\gamma L} \|\nabla f(x_{i}) - \nabla f(x_{i+1})\|^2.%
\end{align}
Moreover, if $\gamma \leq \frac{2}{L}$, then the gradient norm is non-increasing.%
\end{lemma}%%
\section{Proofs of performance bounds}\label{sec:Proofs}
    \shownewtextfinal{In general, solving the associated PEP which provides (numerical) exact performance bounds in closed form is difficult. We propose a systematic procedure which eliminates the standard PEP numerical guessing by effectively canceling specific terms in a linear combination of a well defined set of inequalities from \cref{def:dist-N_interp_cond} to mimic the shape of the performance metric, as described in \cref{remark:main_idea_cancellations}. }%

\shownewtextfinal{The PEP framework facilitated the identification and validation of worst-case scenarios, streamlining the necessary inequalities in the proofs, the subsequent demonstrations being inspired by the numerical solutions (which are not unique), with an improved understanding.} %
Additionally, the intricate range of sublinear regimes outlined in \cref{thm:wc_GM_hypo_ct_all} could only be validated numerically, but not deducible from the PEP numerical analyses. For a comprehensive understanding of the PEP setup for gradient descent, we refer the interested reader to \cite{taylor_smooth_2017} or \cite{abbaszadehpeivasti2021GM_smooth}.

Similar proofs using distance-1 interpolation conditions to obtain tight convergence rates on different performance metrics are given in \cite[Theorem 2.1, Section 3]{Taylor_Jota_PGM_rates_proofs} where the proximal gradient descent is analysed, assuming strong convexity on the smooth function in the split.%

While using ``distance-2'' inequalities is a valid method for proving tight bounds with stepsizes greater than $\gbari[1]$, we mention that the only fundamental requirement is to use conditions connecting non-consecutive iterates. In fact, PEP experiments reveal several alternative strategies, employing: both distance-2 inequalities but with different multipliers, only distance-2 inequalities of the type $Q_{[i+2,i]}$, or inequalities relating the final iterate $Q_{[N,i]}$ \cite{kim2025proofexactconvergencerate}. Our specific approach is advantageous because of its unique cancellation technique, described in \cref{remark:main_idea_cancellations}, which offers a methodical derivation. However, beyond this procedural benefit, we have not yet developed a deeper intuition.%
\paragraph*{\textbf{Proofs organization.}} \, The lemmas from \cref{subsec:dist_1_lemmas_nonconvex,subsec:dist_1_lemmas_strongly_convex,subsec:dist_2_lemmas} are employed in the rest of the section to prove the convergence rates outlined in \cref{sec:Results_wc_analysis}. \cref{table:summary_bounds_and_tightness} provides a summary of the specific lemmas utilized in the demonstrations and the corresponding worst-case function types. The proofs are organized based on the stepsize thresholds delimiting the different regimes, even for (strongly) convex functions where the rates are expressed using minima in the denominators. %
\noindent\paragraph*{\textbf{Main idea: Target inequality. {}}}\,
Inequalities \eqref{eq:distance-N_interpolation_conditions_I_ij} and \eqref{eq:distance-N_interpolation_conditions_I_ji} contribute via nonnegative weights 
$\alphaij[i][j]$ and $\alphaij[j][i]$, respectively, to form a generalized inequality linking $N$ iterations\footnote{For linear combination of inequalities, we adopt the convention of writing the summed expression on the left-hand side, while stating the corresponding inequalities explicitly on the right-hand side.} such that there exists nonnegative weights $\sigma_i \geq 0$ satisfying:%
\begin{align}\label{eq:I^k_N_general_ineq}\tag{\ENQ}
\sum_{0\leq i < j \leq N} \big(\alphaij[i][j] \Iij[i][j] + \alphaij[j][i] \Iij[j][i]\big)\colon \,\,
    f(x_0) - f(x_N) {}\geq{}  \sum_{i=0}^N \sigma_i \frac{\|\nabla f(x_i)\|^2}{2L}.%
\end{align}%
\begin{proposition}[General procedure]\label{prop:general_procedure_min_grad_nrm}
    Inequality \eqref{eq:I^k_N_general_ineq} implies:
\[
    \frac{f(x_0) - f(x_N)}{\sum_{i=0}^N \sigma_i}
        {}\geq{}
    \frac{1}{2L} \,
    \min_{0\leq i \leq N} \big\{\|\nabla f(x_i)\|^2\big\},
\]
with equality only if all gradient norms $\|\nabla f(x_i)\|^2$ with $\sigma_i>0$ are equal. Moreover, 
\[
    \frac{f(x_0) - f_*}{1 + \sum_{i=0}^N \sigma_i}
        {}\geq{}
    \frac{1}{2L} \,
    \min_{0\leq i \leq N} \big\{\|\nabla f(x_i)\|^2\big\}.
\]
\end{proposition}%
\begin{proof}
   The first inequality results after taking minimum gradient norm in \eqref{eq:I^k_N_general_ineq}. %
The adjustment in the denominator of the second inequality results from combining \eqref{eq:sufficient_decrease_N_*} with \eqref{eq:I^k_N_general_ineq} and taking again the minimum gradient norm. \qed
\end{proof}%
\noindent Following \cref{prop:general_procedure_min_grad_nrm}, all convergence proofs reduce to demonstrating an inequality of the form \eqref{eq:I^k_N_general_ineq} with positive weights $\sigma_i$. We utilize interpolation inequalities \eqref{eq:distance-N_interpolation_conditions_I_ij} and \eqref{eq:distance-N_interpolation_conditions_I_ji} with (i) $|j-i|=1$ to prove \cref{thm:wc_GM_hypo} and (ii) $|j-i| \in \{1,2\}$ to prove \cref{thm:wc_GM_hypo_ct_all} and \cref{thm:strg_cvx_rate}. %

\smallskip
\noindent{\textbf{Notation.}}\, For readability reasons, we sometimes denote $f_i {=} f(x_i)$ and $g_i {=} \nabla f(x_i)$. Recall: $\kappa {=} \frac{\mu}{L}$ is the curvature ratio; $\gamma L \in (0,2)$ is the (normalized) stepsize and $\gbari[\infty] = \tfrac{2}{1+[\kappa]_{+}}$.%
\subsection{Distance-1 lemmas for weakly convex functions}\label{subsec:dist_1_lemmas_nonconvex}
For $|j-i|=1$, inequalities \eqref{eq:distance-N_interpolation_conditions_I_ij} and \eqref{eq:distance-N_interpolation_conditions_I_ji} simplify to:%
\begin{align} 
% \hspace{-.75em}
    f_{i}{-}f_{i+1}
        {}\geq {} &
    \frac{\gamma \mu \, \gamma L {-} 2 \gamma \mu {+} 1}{L{-}\mu} \frac{\|g_{i}\|^2}{2} +
    \frac{1}{L{-}\mu} \frac{\|g_{i+1}\|^{2}}{2} +     
    \frac{\gamma L{-}1}{L{-}\mu} \langle g_{i}, g_{i+1}\rangle; 
\label{eq:hypo_interp_(0,1)} \tag{\Iij[i][i+1]}\\
    %%%%%%%%%%%%%%%%%%%%%%%%%%%%%%%%%%%%%%%%%%%%%%%%%%%%%%%%%%%%%%%%%%%%%%%%%%
% \hspace{-1em}
    f_{i+1}{-}f_{i}
        {}\geq {} &
    \frac{\gamma \mu \gamma L {-} 2 \gamma L {+} 1}{L-\mu} \frac{\|g_{i}\|^2}{2} +
    \frac{1}{L{-}\mu} \frac{\|g_{i+1}\|^{2}}{2} +     
    \frac{\gamma \mu {-} 1}{L{-}\mu} \langle g_{i}, g_{i+1}\rangle.
\label{eq:hypo_interp_(1,0)} \tag{\Iji[i][i+1]}
\end{align}
Note that their \lhs only contains consecutive function values, whereas the \rhs is a combination of gradients, with symmetric coefficients in $\gamma L$ and $\gamma \mu$. %
\begin{lemma}\label{lemma:sufficient_decrease_h_leq_1}
    Let $f \in \mathcal{F}_{\mu,L}$, with $\mu \in (-\infty, L)$, and consider one iteration of \eqref{eq:GM_it_fixed_steps} with stepsize $\gamma L \in(0,1]$. Then the following inequality holds, with equality only if $\nabla f(x_{i}) = \nabla f(x_{i+1})$ or $\gamma L = 1$:
\begin{align}\label{eq:one_step_rate_h_leq_1} \tag{N2SD}
\hspace{-.5em}
    f(x_{i}) {-} f(x_{i+1})
        {}\geq {}
    \frac{\gamma L \, \gamma \mu {-} 2\gamma \mu {+} \gamma L}{L{-}\mu} \frac{\|\nabla f(x_{i})\|^2}{2} +
    \frac{\gamma L}{L{-}\mu} \frac{\|\nabla f(x_{i+1})\|^{2}}{2}.
\end{align}%
    Furthermore, if $\mu \leq 0$, then the gradient norms from \eqref{eq:one_step_rate_h_leq_1} are weighted by positive quantities.
\end{lemma}%
%%%%%%%%%%%%%%%%%%%%%%%%%%%%%%%%%%%%%%%%%%%%%%%%%%%%%%%%%%%%%%%%%%%%%%%%%%%%%%%%%
\begin{proof}\label{proof:lemma:sufficient_decrease_h_leq_1}
    Inequality \eqref{eq:hypo_interp_(0,1)} is reformulated by expressing the inner product in terms of gradient norms, i.e., $2 \big\langle g_{i}, g_{i+1}\big\rangle = \|g_{i}\|^2 + \|g_{i+1}\|^2 - \|g_{i} - g_{i+1}\|^2$, as:%
\begin{align} \label{eq:tight_short_step_sizes}
    f_{i}{-}f_{i+1}
        {}\geq {}
    \frac{\gamma L \, \gamma \mu {-} 2\gamma \mu {+} \gamma L}{L{-}\mu} \frac{\|g_{i}\|^2}{2} +
    \frac{\gamma L}{L{-}\mu} \frac{\|g_{i+1}\|^{2}}{2} +
    \frac{1 {-} \gamma L}{L{-}\mu} \frac{\|g_{i} {-} g_{i+1}\|^{2}}{2}.
\end{align}
Since $\gamma L \in (0,1]$, the mixed term can be neglected, resulting in inequality \eqref{eq:one_step_rate_h_leq_1}, with equality only if $g_{i} = g_{i+1}$ or $\gamma L = 1$. The numerator of first term rewrites as $-\gamma \mu(1-\gamma L) + (\gamma L - \gamma \mu)$, hence it is positive when $\mu \leq 0$.
\qed\end{proof}%
\noindent\cref{lemma:sufficient_decrease_h_leq_1} does not hold for stepsizes $\gamma L \in (1, 2)$.%
\begin{lemma}\label{lemma:sufficient_decrease_h_geq_1}
    Let $f \in \mathcal{F}_{\mu,L}$, with, $\mu \in (-\infty,L)$, and consider one iteration of \eqref{eq:GM_it_fixed_steps} with stepsize $\gamma L \in [1,\gbari[\infty])$. Then the following inequality holds:%
\begin{align} \label{eq:one_step_rate_h_geq_1} \tag{N4SD}
% \hspace{-.25em}
    f(x_{i}){-}f(x_{i+1})
        {}\geq {}
    \frac{\gamma [(2{-}\gamma L)(2{-}\gamma \mu){-}1]}{2{-}\gamma L {- }\gamma \mu}
    \frac{\|\nabla f(x_{i})\|^2}{2} +
    \frac{\gamma}{2{-}\gamma L {-} \gamma \mu} \frac{\|\nabla f(x_{i+1})\|^{2}}{2}.
\end{align}%
When $\gamma L > 1$, equality holds only if:%
\begin{align}\label{eq:equality_case_h_geq_1}
    {\langle \nabla f(x_{i})}, {\nabla f(x_{i+1})\rangle}
        {}{=}{} 
    \frac{(1{-}\gamma \mu)(1{-}\gamma L)}{2{-}\gamma L{-}\gamma \mu} \|\nabla f(x_{i})\|^2 {}{+}{} 
    \frac{1}{2{-}\gamma L{-}\gamma \mu} \|\nabla f(x_{i+1})\|^2.%
\end{align}%
Furthermore, if $\gamma L \in [1,\gbari[1]]$, then the gradient norms from \eqref{eq:one_step_rate_h_geq_1} are weighted by nonnegative quantities.%
\end{lemma}%
\begin{proof}\label{proof:lemma:sufficient_decrease_h_geq_1}
The mixed term in inequality \eqref{eq:tight_short_step_sizes} enters with a negative contribution for $\gamma L > 1$. To counterbalance it, the complementary interpolation inequality \eqref{eq:hypo_interp_(1,0)} is used, written with the equivalent expression:%
\begin{align*}
   f_{i+1}-f_{i}
        {}\geq {}
    \frac{\gamma L \gamma \mu {-} 2 \gamma L {+} \gamma \mu}{L{-}\mu} \frac{\|g_{i}\|^2}{2} +
    \frac{\gamma \mu}{L{-}\mu} \frac{\|g_{i+1}\|^{2}}{2} +
    \frac{1 {-} \gamma \mu}{L{-}\mu} \frac{\|g_{i} {-} g_{i+1}\|^{2}}{2}.%
\end{align*}%
Consider the sum of distance-1 inequalities $\bar{I}_{i}^{1} {}\coloneqq{} \Iij[i][i+1] + \Iji[i][i+1]$:%
\begin{align}\label{eq:dist_1_blended_mixed_grads_diff} \tag{$\bar{I}_{i}^{1}$}
\hspace{-0.5em}
0
    {}\geq{}
\frac{2 \gamma L \gamma \mu {-} \gamma L {-} \gamma \mu}{L{-}\mu} \frac{\|g_{i}\|^2}{2} {}+{}
\frac{\gamma L {+} \gamma \mu}{L{-}\mu} \frac{\|g_{i+1}\|^2}{2} {}+{} 
\frac{2{-}\gamma L {-} \gamma \mu}{L{-}\mu} \frac{\| g_{i} {-} g_{i+1}\|^2}{2}.
\end{align}%
By weighting \eqref{eq:dist_1_blended_mixed_grads_diff} with $\beta(\gamma L,\gamma \mu){}\coloneqq{} \tfrac{\gamma L{-}1}{2{-}\gamma L {-}\gamma \mu} \geq 0$ and summing it with \eqref{eq:tight_short_step_sizes}, the coefficient of $\| g_{i} {-} g_{i+1}\|^2$ gets canceled and we obtain \eqref{eq:one_step_rate_h_geq_1}.
The equality condition \eqref{eq:equality_case_h_geq_1} follows from the equality case of \eqref{eq:dist_1_blended_mixed_grads_diff}. Nevertheless, the coefficient of $\|g_{i}\|^2$ is nonnegative for $\gamma L \in [1,\gbari[1]]$, where $\gbari[1]$ is the largest stepsize $\gamma L \in [1, 2)$ for which $(2-\gamma L)(2-\gamma \mu) - 1 \geq {} 0$, implying nonnegativity of the weight of $\|\nabla f(x_{i})\|^2$ in \eqref{eq:one_step_rate_h_geq_1}.\qed
\end{proof}%
The threshold $\gbari[1]$ has the explicit formula $\gbari[1] {}={} \frac{3}{1+\kappa + \sqrt{1 - \kappa + \kappa^2}}$.%
\begin{remark}[Sufficient decrease results] \label{remark:descent_lemmas_Teboulle_22} 
Inequalities \eqref{eq:one_step_rate_h_leq_1} and \eqref{eq:one_step_rate_h_geq_1} are extensions to \textit{weakly convex} functions of the ``double'' sufficient decrease formula \eqref{eq:2SD_Teboulle} \cite[Lemma 1]{Teboulle2022} and the ``quadruple sufficient decrease'' formula \eqref{eq:4SD_Teboulle} \cite[Lemma 5]{Teboulle2022}, respectively, obtained for convex functions ($\mu=0$):%
    \begin{align}
        \label{eq:2SD_Teboulle} \tag{2SD}
            f(x_{i}){-}f(x_{i+1})
                {}\geq {} &
            \gamma \frac{\|\nabla f(x_{i})\|^2}{2} +
            \gamma \frac{\|\nabla f(x_{i+1})\|^2}{2} 
            % \qquad \qquad \qquad \quad 
            ,\,\, {}\forall{} \gamma L \in (0,1]; \\
        %%%%%%%%%%%%%%%%%%%%%%%%%%%%%%%%%%%%%%%%%%%%%%%%%%%%%%%%%%%%%%%%%%%%%%%%%%%%%%%%%%%
        \label{eq:4SD_Teboulle} \tag{4SD}
            f(x_{i}){-}f(x_{i+1})
                {}\geq {} &
            \frac{\gamma (3 {-} 2 \gamma L)}
            {2{-}\gamma L}
            \frac{\|\nabla f(x_{i})\|^2}{2} +
            \frac{\gamma }{2{-}\gamma L} \frac{\|\nabla f(x_{i+1})\|^{2}}{2}
            % \qquad 
            ,\,\,{}\forall{} \gamma L \in [1, 2).
    \end{align}
\end{remark}%%
\subsection{Distance-1 lemmas for strongly convex functions}\label{subsec:dist_1_lemmas_strongly_convex}
Further on, the analysis focuses on constant stepsizes $\gamma L \in (1, 2)$, interval only partially covered by \cref{thm:wc_GM_hypo}. %
% \medskip
\paragraph*{\textbf{Additional notation. }}\, For brevity, we sometimes use $\n{=}\n(\gamma L){}{\coloneqq}{}1{-}\gamma L \in (-1, 0)$, $\m{=}\m(\gamma \mu){}{\coloneqq}{}1{-}\gamma \mu$. Then $T_k(\gamma L, \kappa) = E_k(1{-}\kappa \gamma L) {-} E_k(1{-}\gamma L)$ (see \cref{def:def_TN}), where
$E_k(x) = \sum_{j=1}^{2k} x^{-j}$, can be equivalently written as $\Tinm{}\coloneqq{}\TigLmu$, with $\Tinm=\Ei[k][\m]-\Ei[k][\n]$. Recall $\gbari[\infty] = \frac{2}{1+[\kappa]_{+}}$.%
\begin{definition}\label{def:N_bar}
Let $\Nbar
        {}{}\coloneqq{}{}
    \max \{ k \in \mathbb{N} \,|\, \TigLmu[k] {}\geq{} 0\}$.
\end{definition}%
\begin{proposition}\label{prop:properties:m+n>0}%[Properties for distance-2]
    Let $\mu \in (-\infty, L)$ and $\gamma L \in \big(1, \gbari[\infty] \big)$. Then $\m {}>{}  0$, $\m + \n {}>0$ and $\n^{-2} {}>{} \m^{-2}$.
\end{proposition}%
\begin{proof}%
Condition $\gamma L < \gbari[\infty]$ implies $\gamma L < \frac{2}{1+[\kappa]_{+}}$, hence $0 < 2 {-} \gamma L {-} \gamma L [\kappa]_{+} {}\leq{} 2 {-} \gamma L {-} \gamma L \kappa$, implying $\m + \n > 0$. Since $\n = 1-\gamma L \in (-1,0)$, we have $\m=1-\gamma \mu > 0$. Condition $L>\mu$ implies $\m>\n$ and together with $\m+\n>0$ we conclude that $\m^2 > \n^2$, therefore $\n^{-2} > \m^{-2}$.\qed
\end{proof}
\begin{lemma}\label{lemma:dist_1_L_1_step_decrease}
Let $f \in \mathcal{F}_{\mu,L}$, with $\mu \in (-\infty,L)$, and consider one iteration of \eqref{eq:GM_it_fixed_steps} from $x_i$ with stepsize $\gamma L \in (1, 2)$. Then it holds:
\begin{align} \label{eq:dist_1_L_1_step_decrease} 
    \begin{aligned}    
    \frac{f(x_{i}) - f(x_{i+1})}{\gamma} 
                {}\geq{}
            {-}\Ei[i][\n] \frac{\|\nabla f(x_{i})\|^2}{2} {}+{}
            \Ei[i+1][\n] \frac{\|\nabla f(x_{i+1})\|^2}{2} {}+{}            
            \\ {}+{}
            \frac{ 1 - (\m+\n) \Ei[i+1][\n] }{\m-\n} \frac{\|\nabla f(x_{i+1})-\n \nabla f(x_{i})\|^2}{2}.
    \end{aligned}
\end{align}%
Furthermore, if $\mu \in (0,L)$ and $\gamma L \in [\frac{2}{1+\kappa}, 2)$, then $1-(\m+\n) \Ei[i][\n] > 0$.
\end{lemma}
\begin{proof}
Let $\beta_{i+1}(\n) {{}{}\coloneqq{}{}} (-\n) \Ei[i+1][\n] {}>{} 0$. Inequality \eqref{eq:dist_1_L_1_step_decrease} results after performing the simplifications in the linear combination of interpolation inequalities:
$$
\Iij[i][i+1] + \beta_{i+1}(\n) \, (\Iij[i][i+1] + \Iji[i][i+1]),$$
where \eqref{eq:distance-N_interpolation_conditions_I_ij} and \eqref{eq:distance-N_interpolation_conditions_I_ji} are equivalently rewritten for the pair $(i,i+1)$ as:%
$$
    \begin{aligned}
    (1+\beta_{i+1}) {\,} \Iij[i][i+1] {}:{}
        \frac{f_{i} {-} f_{i+1}}{\gamma}
            {}\geq{}&
        (1+\n) \frac{\|g_{i}\|^2}{2} {+} \frac{1}{\m{-}\n} \frac{\|g_{i+1}{-}\n g_{i}\|^2}{2}; \\
%%%%%%%%%%%%%%%%%%%%%%%%%%%%%%%%%%%%%%%%%%%%%%%%%%%%%%%%%%%%%%%%%%%%%%%%%%%
    \beta_{i+1} {\,} \Iji[i][i+1] {}:{}
        \frac{f_{i+1} {-} f_{i}}{\gamma}
            {}\geq{}&
        {-}\frac{\|g_{i}\|^2}{2} {-} \frac{1}{\n} \frac{\|g_{i+1}\|^2}{2} {+} \frac{\m}{\n} \frac{1}{\m{-}\n} \frac{\|g_{i+1}{-}\n g_{i}\|^2}{2}.
    \end{aligned}
$$%
\noindent{}With $\gamma L \in {[}\frac{2}{1+\kappa}, 2)$, we have $\m{+}\n \leq 0$ (see \cref{prop:properties:m+n>0}) and therefore $1-(\m+\n)\Ei[i+1][\n] > 0$.%
\qed \end{proof}%
\begin{lemma}\label{lemma:dist_1_mu_1_step_decrease}
Let $f \in \mathcal{F}_{\mu,L}$, with $\mu \in (-\infty,L)$, and consider one iteration of \eqref{eq:GM_it_fixed_steps} from $x_i$ with stepsize $\gamma L \in \big(0, \gbari[\infty]\big)$. Then it holds:%
\begin{align} \label{eq:dist_1_mu_1_step_decrease}
    \begin{aligned}    
    \frac{f(x_{i}) - f(x_{i+1})}{\gamma} 
                {}\geq{} 
            -\Ei[i][\m] \frac{\|\nabla f(x_{i})\|^2}{2} {}+{}
             \Ei[i+1][\m] \frac{\|\nabla f(x_{i+1})\|^2}{2} {}+{}            
            \\ {}+{}
            \frac{{-1} +(\m+\n) \Ei[i+1][\m]}{\m-\n} \frac{\|\nabla f(x_{i+1})-\m \nabla f(x_{i})\|^2}{2}.
    \end{aligned}
\end{align}%
Furthermore, if $\mu \in (0,L)$ and $\gamma L \in \big(0, \gbari[i+1]\big)$, then ${-1} +(\m+\n) \Ei[i+1][\m] > 0$.%
\end{lemma}%
\begin{proof}%
Firstly, note that $\m=1-\gamma \mu > 0$; the case $\gamma L \in \big(1, \frac{2}{1+[\kappa]_{+}}\big)$ is covered by \cref{prop:properties:m+n>0}, while for $\gamma L \in (0,1]$ we have $1-\gamma \mu > 1-\gamma L \geq 0$. Let:%
$$
    \beta_{i+1}(\m) {}{}\coloneqq{}{} {-1} + \m \Ei[i+1][\m] {}={} \frac{1 + \m + \dots + \m^{2i}}{\m^{2i+1}} > 0.
$$%
Inequality \eqref{eq:dist_1_mu_1_step_decrease} results after performing the simplifications in the linear combination of interpolation inequalities:%
$$
\Iij[i][i+1] + \beta_{i+1}(\m) \, (\Iij[i][i+1] + \Iji[i][i+1]),$$
where \eqref{eq:distance-N_interpolation_conditions_I_ij} and \eqref{eq:distance-N_interpolation_conditions_I_ji} equivalently rewrite for the pair $(i,i+1)$ as:
\begin{align*}
    (1+ \beta_{i+1}) {\,} \Iij[i][i+1] {}:{}\,{}{} 
        \frac{f_{i} {-} f_{i+1}}{\gamma}
            {}\geq{} &
        \frac{\|g_{i}\|^2}{2} {+} \frac{1}{\m} \frac{\|g_{i+1}\|^2}{2} {+} \frac{\frac{\n}{\m}}{\m{-}\n} \frac{\|g_{i+1}{-}\m g_{i}\|^2}{2}; \\    %%%%%%%%%%%%%%%%%%%%%%%%%%%%%%%%%%%%%%%%%%%%%%%%%%%%%%%%%%%%%%%%%%%%%%%%%%%
     \beta_{i+1} {\,} \Iji[i][i+1] {}:{}\,{}{}
        \frac{f_{i+1} {-} f_{i}}{\gamma}
            {}\geq{} &
        {-}(1{+}\m) \frac{\|g_{i}\|^2}{2} {+} \frac{1}{\m{-}\n} \frac{\|g_{i+1}{-}\m g_{i}\|^2}{2}.
\end{align*}%
Condition $\gamma L < \gbari[i+1]$ implies $\Tinm[i+1] < 0$, hence $i \geq \Nbar$ (recall \cref{def:N_bar}: $\Tinm[\bar{N}] \geq 0$ and $\Tinm[\bar{N}+1] < 0$, corresponding to $\gamma L < \gbari[\bar{N}+1]$). Then the coefficient of the mixed term in \eqref{eq:dist_1_mu_1_step_decrease} is positive:%
$$
    \frac{{-1} {+} (\m{+}\n) E_{i+1}(\m)}{\m{-}\n} 
        {}\geq{} 
    \frac{{-1} {+} (\m{+}\n) E_{\bar{N}+1}(\m)}{\m{-}\n}
        {}={}
    \frac{\Tinm[\bar{N}]{-}\n^2 \Tinm[\bar{N}+1]}{(\m{-}\n)^2} 
        {}>{} 0.
$$
The first inequality results from the monotonic increase with index $i$ of $E_{i+1}(\m)$, whereas the last one by the definition of $\Nbar$.%
\qed \end{proof}%
\subsection{Distance-2 lemmas}\label{subsec:dist_2_lemmas}
In this section we derive descent lemmas involving distance-2 interpolation inequalities, needed to prove rates on the stepsize range $\gamma L \in (1, \gbari[\infty])$.%
\begin{lemma}[Central distance-2 result]\label{lemma:telescoped_main_lemma}
Let $f \in \mathcal{F}_{\mu,L}$, with $\mu \in (-\infty, L)$. Consider $k+1$ iterations of \eqref{eq:GM_it_fixed_steps} with constant stepsize $\gamma L \in \big[\gbari[k], \gbari[\infty] \big)$ starting from $x_0$, where $k \geq 1$. Then it holds:%
\begin{align}
        & \frac{f(x_{0}) {-} f(x_{k})}{\gamma} +
        \frac{(\m{+}\n)\Tinm[k]}{2(\m{-}\n)} \frac{f(x_{k}) {-} f(x_{k+1})}{\gamma}
            \geq
        % \\ & {}\qquad{}
        \left(\Ei[k][\m] {+} \Ei[k][\n]\right) \frac{\|\nabla f(x_{k})\|^2}{4} +
        \nonumber \\ & {}\quad{}\,
        {}+{}\frac{\Tinm[k]}{2(\m{-}\n)} \Big[ \m \n \frac{\|\nabla f(x_{k})\|^2}{2} + 
        \langle \nabla f(x_{k}) , \nabla f(x_{k+1}) \rangle + 
        \frac{\|\nabla f(x_{k+1})\|^2}{2} \Big]. \label{eq:telescoped_main_lemma} \tag{D2}
\end{align}%
\end{lemma}%
\begin{proof}
Consider two iterations connecting $x_{i}$, $x_{i+1}$ and $x_{i+2}$, with $i \geq 0$. The interpolation inequalities \eqref{eq:distance-N_interpolation_conditions_I_ij} and \eqref{eq:distance-N_interpolation_conditions_I_ji} written for distance-1 (\Iij[i][i+1] and \Iji[i][i+1]) and distance-2 (\Iij[i][i+2] and \Iji[i][i+2]), respectively:%
\begin{align*}
    &\begin{aligned}    
    % & \alpha_{[i-1>i]} \times 
    \Iij[i][i+1]\colon 
        \frac{f_{i} - f_{i+1}}{\gamma}
            {}\geq{}
        (1+\n) \frac{\|g_{i}\|^2}{2} + \frac{1}{\m-\n} \frac{\|g_{i+1}-\n g_{i}\|^2}{2}; 
    \end{aligned}\\
    %%%%%%%%%%%%%%%%%%%%%%%%%%%%%%%%%%%%%%%%%%%%%%%%%%%%%%%%%%%%%%%%%%%%%%%%%%%
    &\begin{aligned}    
    % & \alpha_{[i-1<i]}  \times 
    \Iji[i][i+1]\colon 
        \frac{f_{i+1} - f_{i}}{\gamma}
            {}\geq{}
        -\frac{\|g_{i}\|^2}{2} - \frac{1}{\n} \frac{\|g_{i+1}\|^2}{2} + \frac{\m}{\n} \frac{1}{\m-\n} \frac{\|g_{i+1}-\n g_{i}\|^2}{2}; 
    \end{aligned}\\
    %%%%%%%%%%%%%%%%%%%%%%%%%%%%%%%%%%%%%%%%%%%%%%%%%%%%%%%%%%%%%%%%%%%%%%%%%%%
    &\begin{aligned}    
    &% \alpha_{[i-1>i+1]}  \times 
    \Iij[i][i+2]\colon 
        \frac{f_{i} - f_{i+2}}{\gamma}
            {}\geq{}
        (1+\n+\n^2) \frac{\|g_{i}\|^2}{2} + \n \frac{\|g_{i+1}\|^2}{2} + 
        \frac{1}{\n(\m-\n)} \cdot
        \\ & {}\qquad{}        
        \Big[
        \n(1-\m) \frac{\|g_{i+1}-\n g_{i}\|^2}{2}
        - (1-\n) \frac{\|g_{i+2}-\n g_{i+1}\|^2}{2}
        + \frac{\|g_{i+2}-\n^2 g_{i}\|^2}{2} 
        \Big];
    \end{aligned}\\
    %%%%%%%%%%%%%%%%%%%%%%%%%%%%%%%%%%%%%%%%%%%%%%%%%%%%%%%%%%%%%%%%%%%%%%%%%%%
    &\begin{aligned}    
    &% \alpha_{[i-1<i+1]}  \times 
    \Iji[i][i+2]\colon
        \frac{f_{i+2} - f_{i}}{\gamma}
            {}\geq{} 
        - \frac{\|g_{i}\|^2}{2} - \frac{1}{\n} \frac{\|g_{i+1}\|^2}{2}
        - \frac{1+\n}{\n^2} \frac{\|g_{i+2}\|^2}{2} + \frac{1}{\n(\m-\n)} \cdot
        \\ & {}\qquad{} 
        \Big[
        \m(1-\n) \frac{\|g_{i+1}-\n g_{i}\|^2}{2}
        - (1-\m) \frac{\|g_{i+2}-\n g_{i+1}\|^2}{2}
        + \frac{\m}{\n} \frac{\|g_{i+2}-\n^2 g_{i}\|^2}{2} 
        \Big],
    \end{aligned}
\end{align*}%
are summed up in a linear combination with \textit{nonnegative} corresponding weights:%
$$
        \eqref{eq:dist_2_atom_ineq} \equiv
\alphaij[i][i+1] \Iij[i][i+1] +
\alphaij[i+1][i] \Iji[i][i+1] +
\alphaij[i][i+2] \Iij[i][i+2] +
\alphaij[i+2][i] \Iji[i][i+2],
$$%
where:%
$$
\begin{aligned}
    \alphaij[i][i+1] {}{}\coloneqq{}{}& 
    1 - \n \Ei[i+1][\n] + \frac{(1-\m) \Tinm[i]}{2(\m-\n)} - \frac{\m(1+\n) \Tinm[i+1]}{2(\m-\n)};
    \\
   \alphaij[i+1][i] {}{}\coloneqq{}{}& 
    -\n \Ei[i+1][\n] + \frac{(1+\n) \Tinm[i]}{2(\m-\n)} + \frac{\n(1-\m) \Tinm[i+1]}{2(\m-\n)};
    \\
    \alphaij[i][i+2] {}{}\coloneqq{}{}&
     \frac{\m \Tinm[i+1]}{2(\m-\n)}; 
     \\
    \alphaij[i+2][i] {}{}\coloneqq{}{}&
    \frac{-\n \Tinm[i+1]}{2(\m-\n)},
\end{aligned}%
$$%
and \eqref{eq:dist_2_atom_ineq} is obtained after performing the simplifications:%
\begin{align}\label{eq:dist_2_atom_ineq}    
    %%%%%%%%%%%%%%%%%%%%%%%%%%%%%%%%%%%%%%%%%%%
    \begin{aligned}
        &
        \frac{f_{i} {-} f_{i+1}}{\gamma} + 
        (\m{+}\n) \Big[
        {-}\frac{\Tinm[i]}{2(\m{-}\n)} \frac{f_{i} {-} f_{i+1}}{\gamma} +
           \frac{\Tinm[i+1]}{2(\m{-}\n)} \frac{f_{i+1} {-} f_{i+2}}{\gamma}
        \Big]
            {}\geq{}
        \\ & {}\,{} 
        {}{-}{}\frac{\Tinm[i]}{2(\m{-}\n)} \Big[ \m \n \frac{\|g_{i}\|^2}{2} + 
        \langle g_{i}, g_{i+1} \rangle + \frac{\|g_{i+1}\|^2}{2}  \Big] 
        {-}\left(\Ei[i][\m] {+} \Ei[i][\n]\right) \frac{\|g_{i}\|^2}{4}
        \\ & {}\,{} 
        {}{+}{}\frac{\Tinm[i+1]}{2(\m{-}\n)} \Big[ \m \n \frac{\|g_{i+1}\|^2}{2} 
        {+} {\langle g_{i+1}, g_{i+2} \rangle} + \frac{\|g_{i+2}\|^2}{2}  \Big] {}+{} \left(\Ei[i+1][\m] {+} \Ei[i+1][\n]\right) \frac{\|g_{i+1}\|^2}{4}.
    \end{aligned}
\end{align}%
\noindent{}Telescoping it for indices $i=0,1,\dots,k-1$ yields inequality \eqref{eq:telescoped_main_lemma}. It remains to show the nonnegativity of multipliers. Since the stepsize thresholds increase with index $i$, condition $\gamma L \geq \gbari[k]$ ensures $\Tinm[i+1] \geq 0$, for all $i=0,{\dots},k-1$. This directly implies the positivity of $\alphaij[i][i+2]$ and $\alphaij[i+2][i]$. The positivity of distance-1 multipliers is evidenced by their equivalent representations as nonnegative sums (recall $\n \in (-1,0)$, $\m \in \left(0,\infty\right)$; $\m+\n>0$ and $\n^{-2} > \m^{-2}$ from \cref{prop:properties:m+n>0}):%
$$
\begin{aligned}
    &\begin{aligned}    
    2\alphaij[i][i+1] {}={}& 
    1 + \frac{\m+\m^{-2i}}{1+\m} + 
    \frac{(1+\n)(\m+\n-\n\m)}{\m \n^2} \frac{{-1} +\n^{-2(i+1)}}{{-1} +\n^{-2}} + 
    \\& {}\qquad\qquad{}
    \frac{(1+\n)(1+\m^2)}{\m(\m-\n)} 
    % \Big[ \frac{1-\n^{-2i}}{1-\n^{-2}} - \frac{1-\m^{-2i}}{1-\m^{-2}} \Big] 
    \sum_{j=0}^{i} \big(\n^{-2j} - \m^{-2j}\big)
    {}>{} 0;  
    \end{aligned} \\%
    %%%%%%%%%%%%%%%%%%%%%%%%%%%%%%%%%%%%%%%%%%
    &\begin{aligned}  
    2\alphaij[i+1][i] {}={}&
    \frac{{-1} +\n^{-2i}}{1-\n} + 
    \frac{\m(1+\n)}{1-\n} \frac{{-1} +\m^{-2(i+1)}}{1-\m}
    + \\& {}\qquad\qquad{}
     +\frac{-\n \big[1-\n+\m(1+\n)\big]}{1-\n} \frac{\n^{-2(i+1)} - \m^{-2(i+1)}}{\m-\n} > 0. 
    % {}\qquad\qquad\qquad\qquad\qquad{}
    {}\qquad{}
    {}\qed{}
    \end{aligned}%
\end{aligned}
$$%
\end{proof}%
%%%%%%%%%%%%%%%%%%%%%%%%%%%%%%%%%%%%%%%%%%%%%%%
\begin{remark}[\textbf{Derivation of \cref{lemma:telescoped_main_lemma}}]\label{remark:main_idea_cancellations}%
    Inequality \eqref{eq:telescoped_main_lemma} selectively involves gradient norms of indices $\{k,k+1\}$, function values of indices $\{0,k,k+1\}$ and the inner product $\langle \nabla f(x_{k}),\nabla f(x_{k+1})\rangle$. As a result, the linear combination of interpolation inequalities eliminates all other function values and gradients. Specifically, it removes: $k$ gradient norms, $k$ consecutive and $k$ distance-2 inner products of gradients, $k-1$ function values, while setting the coefficient of $f_0$ to $\frac{1}{\gamma}$, totaling $4k$ constraints. Remarkably, the linear combination employs exactly $4k$ inequalities as well, highlighting why we needed to involve distance-2 interpolation inequalities. Consequently, deriving \cref{lemma:telescoped_main_lemma} (and subsequent convergence rates) entails solving a linear system, contrasting with the PEP numerical analysis, which involves solving an SDP.%
\end{remark}%
%%%%%%%%%%%%%%%%%%%%%%%%%%%%%%%%%%%%%%%%%%%%%%%
Subsequent lemmas involving distance-2 interpolation inequalities are derived by augmenting \cref{lemma:telescoped_main_lemma} with distance-1 interpolation inequalities.%
\begin{lemma}[Multistep descent]\label{lemma:hypo_GN4SD}
Let $f \in \mathcal{F}_{\mu,L}$, with $\mu \in (-\infty, L)$, and consider $k+1 {}\geq{} 1$ iterations of \eqref{eq:GM_it_fixed_steps} with stepsize $\gamma L \in [\gbari[k], \gbari[\infty])$, {\shownewtext{satisfying $\gamma L > 1$,}} starting from $x_0$. Then the following inequality holds:%
\begin{align}%
% \hspace{-.25em}
    \frac{f(x_0) {-} f(x_{k+1})}{\gamma}
        {}\geq{}
    \frac{{-}\m^2 \n^2 \Tinm[k+1]}{\m^2{-}\n^2} \frac{\|\nabla f(x_{k})\|^2}{2} {}+{}
    \frac{\Tinm[k]{+}(\m{-}\n)}{\m^2{-}\n^2} \frac{\|\nabla f(x_{k+1})\|^2}{2}.
\label{eq:general_ineq_proof} \tag{GN4SD}
% \tag{$I_{N+1}^{N}$}
\end{align}%
Furthermore, if $\mu \in (-\infty, 0]$ and $\gamma L \in [\gbari[k], \gbari[k+1])$, then the gradient norms from \eqref{eq:general_ineq_proof} are weighted by nonnegative quantities.
\end{lemma}%
\begin{proof}%
For $k=0$, the result follows from \cref{lemma:sufficient_decrease_h_geq_1} and we obtain exactly \eqref{eq:one_step_rate_h_geq_1}. 
For $k\geq 1$, inequality \eqref{eq:general_ineq_proof} is derived by adjusting \eqref{eq:telescoped_main_lemma} through the linear combination of interpolation inequalities $\alphaij[k][k+1] \Iij[k][k+1] + \alphaij[k+1][k] \Iji[k][k+1]$:%
$$
    \begin{aligned}
        &  \Big[1 - \frac{(\m+\n)\,\Tinm[k]}{2(\m-\n)} \Big] \frac{f_{k} - f_{k+1}}{\gamma}
        {}\geq{}
        \\ & {}\qquad{}
        \Big(-\Ei[k][\m] - \Ei[k][\n] - \frac{\m \n \Tinm[k]}{\m-\n} {}-{} \frac{2\m^2\n^2 \Tinm[k+1]}{\m^2-\n^2} \Big) \frac{\|g_{k}\|^2}{4}
        \\ & {}\qquad{}        
        % \\ & {}\qquad{}
        {}+{} \Big(\frac{\m-\n + \Tinm[k]}{\m^2-\n^2} - \frac{\Tinm[k]}{2(\m-\n)} \Big) \frac{\|g_{k+1}\|^2}{2}
        {}-{}\frac{\Tinm[k]}{2(\m-\n)} \langle g_{k},g_{k+1} \rangle,
    \end{aligned}
$$%
where the corresponding weights are defined as follows:%
    $$
    \begin{aligned}
        \alphaij[k][k+1]
             {}{}\coloneqq{}{} &
        1 + \frac{-\n}{\m+\n} + \frac{(1-\m)(\m+\n)-2\n}{2(\m^2-\n^2)} \Tinm[k];  \\
        %%%%%%%%%%%%%%%%%%%%%%%%%%%%%%%%%%%%%%%%%%%%%%%%%%%%%%%%%%%%%%%%%%%%%%%%%%%%%%%%%%%%
        \alphaij[k+1][k]
            {}{}\coloneqq{}{} &
        \frac{-\n}{\m+\n} + \frac{(1+\n)(\m+\n)-2\n}{2(\m^2-\n^2)} \Tinm[k].
    \end{aligned}
    $$ %
    Multiplier $\alphaij[k+1][k]$ is represented as a sum of positive quantities ($\n < 0$).
    Multiplier $\alphaij[k][k+1]$ is directly written as sum of positive quantities if $(1-\m)(\m+\n)-2\n > 0$ (in particular, this covers the (strongly) convex case where $\m \leq 1$). Otherwise, $\alphaij[k][k+1]$ is equivalently expressed as a sum of positive terms as follows, with $\m>1$:%
$$
\begin{aligned}
    \alphaij[k][k+1]
        {}={}&
    \frac{\m}{2(\m+\n)} + \frac{1}{2(\m^2-\n^2)} 
        \bigg[ \frac{-\n(\m-\n)(1+\m^2)}{(\m-1)(1-\n)} {}+{}
        \\ & {}\quad{} 
        -\big[(1-\m)(\m+\n)-2\n\big] \Big( \frac{\m^{-2k}}{\m-1} + \frac{\n^{-2k}}{1-\n} \Big)  \bigg]
        {}>{} 0. {}\quad{} \qed
\end{aligned}
$$
\end{proof}%
\begin{lemma}\label{lemma:general_lemma_ct_h_to_Linear} 
Let $f \in \mathcal{F}_{\mu,L}$, with $\mu \in (-\infty, L)$, and consider $k+1 \geq 1$ iterations of \eqref{eq:GM_it_fixed_steps} with stepsize $\gamma L \in [\gbari[k], \gbari[\infty] \big)$, satisfying $\gamma L > 1$, starting from $x_0$. Then the following inequality holds:%
\begin{align}\label{eq:general_ineq_proof_to_linear_compact}
   \frac{f(x_0) {-} f(x_{k+1})}{\gamma}
        {}\geq{} 
    \Ei[k+1][\n] \frac{\|\nabla f(x_{k+1})\|^2}{2} {+} 
    \frac{\m^2 \Tinm[k+1]}{(\m{-}\n)^2}
    \frac{\| \nabla f(x_{k+1}) {-} \n \nabla f(x_{k}) \|^2}{2}.
\end{align}%
Furthermore, if $\gamma L \in [\gbari[k+1], \gbari[\infty] \big)$, then $\Tinm[k+1]\geq 0 $.%
\end{lemma}%
\begin{proof}%
For $k=0$, the result follows from \cref{lemma:dist_1_L_1_step_decrease}. For $k \geq 1$, inequality \eqref{eq:general_ineq_proof_to_linear_compact} is derived by adjusting \eqref{eq:telescoped_main_lemma} with the linear combination of interpolation inequalities $\alphaij[k][k+1] \Iij[k][k+1] + \alphaij[k][k+1] \Iji[k][k+1]$:%
$$
    % \hspace{-1em}
    \begin{aligned}
        &\Big[1 - \frac{(\m+\n)\,\Tinm[k]}{2(\m-\n)} \Big] \frac{f_{k} - f_{k+1}}{\gamma}
            {}\geq{}
        \Big(-\frac{\Tinm[k]}{2(\m-\n)} - \frac{\m^2 \n \Tinm[k+1]}{(\m-\n)^2}\Big) \langle g_{k},g_{k+1} \rangle {}+{}
        \\ & {}\quad\qquad{}
        {}+{} \Big[ \frac{2\m^2 \n^2 \Tinm[k+1]}{(\m-\n)^2} -\Ei[k][\m] - \Ei[k][\n] - \frac{\m \n \Tinm[k]}{\m-\n} \Big] \frac{\|g_{k}\|^2}{4} {}+{}
        \\ & {}\quad\qquad{}
        {}+{} \Big[ \frac{\m^2 \, \Tinm[k+1]}{(\m-\n)^2}  + \Ei[k+1][\n] -\frac{\Tinm[k]}{2(\m-\n)} \Big] \frac{\|g_{k+1}\|^2}{2},
    \end{aligned}
$$%
where the corresponding weights are defined as:
$$
\begin{aligned}
    \alphaij[k][k+1] 
        {}{}\coloneqq{}{} &
    1+(-\n)\Ei[k+1][\n] + \frac{1-\m}{2(\m-\n)} \Tinm[k];
    \\
    %%%%%%%%%%%%%%%%%%%%%%%%%%%%%%%%%%%%%%%%%%%%%%%%%%%%%%%%%%%%%%%%%%%%%%%%%%%%%%%%%%%%
    \alphaij[k+1][k]
        {}{}\coloneqq{}{} &
    (-\n) \Ei[k+1][\n] + 
    \frac{1+\n}{2(\m-\n)} \Tinm[k].
\end{aligned}%
$$%
Multiplier $\alphaij[k+1][k]$ is represented as a sum of positive quantities, while $\alphaij[k][k+1]$ is directly written as sum of positive terms if $\m \leq 1$ (the (strongly) convex case). Otherwise, $\m > 1$ and $\alphaij[k][k+1]$ is equivalently expressed as sum of positive terms as:%
$$
        \alphaij[k][k+1]
    {}={}
        \frac{1+(-\n) \n^{-2(k+1)}}{2(1-\n)} + \m \frac{\m \m^{-2(k+1)} + (-\n) \n^{-2(k+1)}}{2(\m-\n)}. {}\quad{} \qed
$$
\end{proof}%
\begin{lemma}\label{lemma:general_lemma_ct_h_to_Linear_mu} 
Let $f \in \mathcal{F}_{\mu,L}$, with $\mu \in (-\infty,L)$, and consider $k+1 \geq 1$ iterations of \eqref{eq:GM_it_fixed_steps} with stepsize $\gamma L \in [\gbari[k], \gbari[\infty])$, satisfying $\gamma L > 1$, starting from $x_0$. Then the following inequality holds:%
\begin{align}\label{eq:general_ineq_proof_to_linear_mu_compact}
% \hspace{-.75em}
    \frac{f(x_0) {-} f(x_{k+1})}{\gamma}
        {\geq}
    \Ei[k+1][\m] \frac{\|\nabla f(x_{k+1})\|^2}{2} {+} 
    \frac{-\n^2 \Tinm[k+1]}{(\m{-}\n)^2}
    \frac{\| \nabla f(x_{k+1}) {-} \m \nabla f(x_{k}) \|^2}{2}.
\end{align}%
Furthermore, if $\gamma L \in [\gbari[k], \gbari[k+1] \big)$, then $\Tinm[k+1]<0$.
\end{lemma}%
\begin{proof}
Case $k=0$ results from \cref{lemma:dist_1_mu_1_step_decrease}; further on, $k\geq 1$. Inequality \eqref{eq:direct_sum_interp_cond} written for the pair $(i,j)=(k,k+1)$ reads:%
   \begin{align}\label{eq:dist_1_correction_to_mu}
       0
        {}\geq{}
       \frac{\m \n}{\m-\n} \|g_{k}\|^2 + \frac{1}{\m-\n} \|g_{k+1}\|^2 - \frac{\m+\n}{\m-\n} \langle g_{k} \,,\, g_{k+1} \rangle.
   \end{align}%
   \begin{enumerate}[align=right,itemindent=1em,labelsep=2pt,labelwidth=1em,leftmargin=0pt,nosep]
       \item \textbf{Case $\gamma L \in [\gbari[k], \gbari[k+1])$.} By using $\m>0$, \eqref{eq:dist_1_correction_to_mu} equivalently rewrites as:%
   $$
       0
        {}\geq{}
       -\m^2 \frac{\|g_{k}\|^2}{2} + \frac{\|g_{k+1}\|^2}{2} + \frac{\m+\n}{\m-\n} \frac{\|g_{k+1}-\m g_{k}\|^2}{2}.
   $$%
   By multiplying it with $\frac{-\n^2}{\m^2-\n^2}\Tinm[k+1]$ (positive because $\Tinm[k+1] < 0$) and adding it to \eqref{eq:general_ineq_proof}, the coefficient of $\|g_{k}\|^2$ is canceled and we get:%
    $$
    \frac{f_0{-}f_{k+1}}{\gamma}
        {}{\geq}{}
    \frac{\Tinm[k]{+}(\m{-}\n){-}\n^2 \Tinm[k+1]}{\m^2{-}\n^2} \frac{\|g_{k+1}\|^2}{2} {}+{}
       \frac{{-}\n^2 \Tinm[k+1]}{(\m{-}\n)^2}  \frac{\|g_{k+1}{-}\m g_{k}\|^2}{2},
    $$
    which after simplifications is exactly \eqref{eq:general_ineq_proof_to_linear_mu_compact}.%
    \item \textbf{Case $\gamma L \in \big[\gbari[k+1],\gbari[\infty] \big)$.} After multiplying \eqref{eq:dist_1_correction_to_mu} with $\frac{-\n \m}{\m-\n} > 0$:
    $$
        0
            {}\geq{}
        \frac{\|g_{k+1}\|^2}{2} + 
        \frac{-\m^2}{(\m-\n)^2} \frac{\|g_{k+1}-\n g_{k}\|^2}{2} +
        \frac{-\n^2}{(\m-\n)^2} \frac{\|g_{k+1}-\m g_{k}\|^2}{2}.
    $$
    By scaling it with $\Tinm[k+1]>0$ and summing it with \eqref{eq:general_ineq_proof_to_linear_compact}, the coefficient of $\|g_{k+1} - \n g_{k}\|^2$ is canceled and after simplifications we get \eqref{eq:general_ineq_proof_to_linear_mu_compact}.%
    {}\qed{}%
    \end{enumerate}%
\end{proof}%
\subsection{Proof of performance bounds for weakly convex functions}%
\label{subsec:proofs:nonconvex_rates}%
\begin{proof}[\textbf{\cref{thm:wc_GM_hypo}.}]\label{proof:thm:wc_GM_hypo}
By telescoping the sufficient decrease inequalities \eqref{eq:one_step_rate_h_leq_1} and \eqref{eq:one_step_rate_h_geq_1} for indices $i = 1,\dots,N$ and stepsizes $\gamma_0,\dots,\gamma_{N-1}$, respectively:%
\begin{align}\label{eq:sum_lemmas_N_steps}
\begin{aligned}
% \hspace{-0.1em}
    f_0 {-} f_N
        &{} \geq {}
    \cfrac{\|g_{0}\|^2}{2} \Big[ 
        \frac{\gamma_0 L \, \gamma_0 \mu {-} 2\gamma_0 \mu {+} \gamma_{0} L}{L{-}\mu} \delta_0 {}+{} 
        \frac{\gamma_{0}[(2{-}\gamma_{0} L)(2{-}\gamma_{0} \mu){-}1]}{2{-}\gamma_{0} L {-}\gamma_{0} \mu} (1-\delta_0)
    \Big] %{}{+}{}
    \\ & {}+{} 
         \sum_{i=1}^{N-1} \cfrac{\|g_{i}\|^2}{2}
     \Big[ 
        \frac{\gamma_{i-1}}{L-\mu} \delta_{i-1} {}+{} 
        \frac{\gamma_{i-1}}{2-\gamma_{i-1} L -\gamma_{i-1} \mu} (1-\delta_{i-1}) {}+{}
        \\  & \qquad {}+{}
        \frac{\gamma_i L \, \gamma_i \mu - 2\gamma_i \mu + \gamma_{i} L}{L-\mu} \delta_i {}+{} 
        \frac{\gamma_{i}[(2-\gamma_i L)(2-\gamma_i \mu)-1]}{2-\gamma_{i} L -\gamma_{i} \mu} (1-\delta_{i})
    \Big] 
    \\& {}+{}
    \cfrac{\|g_{N}\|^2}{2} \Big[ 
        \frac{\gamma_{N-1}}{L-\mu} \delta_{N-1} {}+{} 
        \frac{\gamma_{N-1}}{2-\gamma_{N-1} L -\gamma_{N-1} \mu} (1-\delta_{N-1})
    \Big],%
\end{aligned}
\end{align}%
where $\delta_i = 1$ if $\gamma_i \leq \frac{1}{L}$, and $\delta_i = 0$ otherwise. The telescoped inequality resembles a decrease inequality of type \eqref{eq:I^k_N_general_ineq} if all gradient norms have positive weights; a sufficient condition is $\frac{\gamma_{i}[(2{-}\gamma_i L)(2{-}\gamma_i \mu){-}1]}{2{-}\gamma_{i} L {-}\gamma_{i} \mu} {\geq} 0$ for $i=0,{\dots},N-1$, which is satisfied by any choice of $\gamma_i L \in (0, \gbari[1]]$ (see \cref{lemma:sufficient_decrease_h_geq_1}). By taking the minimum gradient norm and reordering:%
\[%
\begin{aligned}
    f_{0}-f_{N}
        {}\geq {}&
    \min_{0\leq i \leq N} \Big\{\cfrac{\|g_i\|^2}{2L}\Big\}
    \Big[
        \sum_{i=0}^{N-1} 
        L 
        \Big(
            \frac{\gamma_i L}{L-\mu} + 
            \frac{\gamma_i L \, \gamma_i \mu - 2\gamma_i \mu + \gamma_{i} L}{L-\mu} 
        \Big) \delta_i {}+{}         
        \\  & \, \quad {}+{}
        \Big(%
            \frac{\gamma_{i} L}{2-\gamma_{i} L -\gamma_{i} \mu} + 
            \frac{\gamma_{i} L[(2-\gamma_i L)(2-\gamma_i \mu)-1]}{2-\gamma_{i} L -\gamma_{i} \mu}
        \Big)
        (1-\delta_{i})
    \Big]
    %%%%%%%%%%%%%%%%%
        \\ {}={} &
    \min_{0\leq i \leq N} \Big\{\cfrac{\|g_i\|^2}{2L}\Big\} \,
    \sum_{i=0}^{N-1} \gamma_i L \,\pigLgmu[\gamma_i],
\end{aligned}%
\]%
with $\pigLgmu[\gamma_i]$ defined in \eqref{eq:pi_hi_rate_GM}. Furthermore, the following particular rates hold:%
\begin{enumerate}[align=right,itemindent=1em,labelsep=2pt,labelwidth=1em,leftmargin=0pt,nosep]
 \item Case $\gamma_i L \in (0,1]$:
$$
\begin{aligned}
    \tfrac{1}{2L} \min_{0 \leq i \leq N} \big\{\|g_i\|^2\big\}
        {}{\leq}{}&
    \min \bigg\{
    \frac{f_0 - f_N}{\sum_{i{=}0}^{N{-}1} \gamma_i L \bigl( 2 {+} \frac{\gamma_i L \gamma_i \mu}{\gamma_i L {-} \gamma_i \mu}\bigr)}
    ,
    %%%%%%%%%%%%%%%%%%%%%%%%%%%%%%%%%%
    % \min_{0 \leq i \leq N} \Big\{\frac{\|g_i\|^2}{2L}\Big\}
    %     {}\leq{}
    \frac{f_0 - f_*}{1 {+} \sum_{i{=}0}^{N{-}1} \gamma_i L \bigl( 2 {+} \frac{\gamma_i L \gamma_i \mu}{\gamma_i L {-} \gamma_i \mu}\bigr)}
    \bigg\};
\end{aligned}
$$
\item Case $\gamma_i L \in [1,\gbari[1]{]}$:
$$
\begin{aligned}
    \tfrac{1}{2L} \min_{0 \leq i \leq N} \big\{\|g_i\|^2\big\}
        {}{\leq}{}&
    \min\bigg\{
    \frac{f_0 - f_N}{\sum_{i{=}0}^{N{-}1} \frac{\gamma_i L(2{-}\gamma_i L)(2{-}\gamma_i \mu)}{2{-}\gamma_i L {-} \gamma_i \mu}}
    ,
    %%%%%%%%%%%%%%%%%%%%%%%%%%%%%%%%%%
    % \min_{0 \leq i \leq N} \Big\{\frac{\|g_i\|^2}{2L}\Big\}
    %     {}\leq{}&
    \frac{f_0 - f_*}{1{+}\sum_{i{=}0}^{N{-}1} \frac{\gamma_i L(2{-}\gamma_i L)(2{-}\gamma_i \mu)}{2{-}\gamma_i L {-} \gamma_i \mu}}
    \bigg\}. \,{}\qed{}
\end{aligned}
$$%
\end{enumerate}%
\end{proof}%
\begin{proof}[\textbf{\cref{thm:wc_GM_hypo_ct_all}}]\label{proof:thm:wc_GM_hypo_ct_all}%
The interval $(0,\gbari[1]{]}$ is addressed by the proof of \cref{thm:wc_GM_hypo}, so it remains to prove the rates for constant stepsizes in the range %
\[
    [\gbari[1],2) {}={} \mathop{\cup}\limits_{k=1}^{N-1} [\gbari[k],\gbari[k+1]) \cup [\gbari[N],2).
\]%
\begin{enumerate}[align=right,itemindent=1em,labelsep=2pt,labelwidth=1em,leftmargin=0pt,nosep]
\item \textbf{Case $\gamma L \in [\gbari[k],\gbari[k+1])$,} with some $k \in \{1,\dots,N-1\}$. To inequality \eqref{eq:general_ineq_proof} we append $(N{-}k{-}1)$ inequalities \eqref{eq:one_step_rate_h_geq_1}, written for indices $i=k+1,\dots,N-1$. Then we get (using the usual notation with $\gamma$, $L$, $\mu$):
$$
    \begin{aligned}    
        \frac{f_{0}{-}f_{N}}{\gamma}
            {}{\geq}{}
        \frac{{-}(1{-}\gamma L)^{2}(1{-}\gamma \mu)^{2}\TigLmu[k+1]}{(1{-}\gamma \mu)^{2} {-} (1{-}\gamma L)^{2}}
        \frac{\|g_{k}\|^2}{2} {+}
        \frac{\TigLmu[k] {+} \gamma L {-} \gamma \mu}{(1{-}\gamma \mu)^{2} {-} (1{-}\gamma L)^{2}}
        \frac{\|g_{k+1}\|^2}{2} {}{+}{}
            \\ 
        \frac{(2{-}\gamma L)(2{-}\gamma \mu){-}1}{2{-}\gamma L {-} \gamma \mu} \frac{\|g_{k+1}\|^2}{2}
        {+} \frac{(2{-}\gamma L)(2{-}\gamma \mu)}{2{-}\gamma L {-} \gamma \mu} \sum_{i={k{+}2}}^{N{-}1} \frac{\|g_{i}\|^2}{2} +
        \frac{1}{2-\gamma L - \gamma \mu}\frac{\|g_{N}\|^2}{2},
    \end{aligned}
$$%
simplifying to
$$
    \begin{aligned}    
        \frac{f_{0}{-}f_{N}}{\gamma}
            {}{\geq}{}
        \frac{{-}(1{-}\gamma L)^{2}(1{-}\gamma \mu)^{2}\TigLmu[k+1]}{(1{-}\gamma \mu)^{2} {-} (1{-}\gamma L)^{2}}
        \frac{\|g_{k}\|^2}{2} +
        \frac{\TigLmu[k]}{(1{-}\gamma \mu)^{2} {-} (1{-}\gamma L)^{2}}
        \frac{\|g_{k+1}\|^2}{2} {}{+}{}
            \\  
        \frac{(2{-}\gamma L)(2{-}\gamma \mu)}{2{-}\gamma L {-} \gamma \mu} \frac{\|g_{k+1}\|^2}{2}
        {+} \frac{(2{-}\gamma L)(2{-}\gamma \mu)}{2{-}\gamma L {-} \gamma \mu} \sum_{i={k{+}2}}^{N{-}1} \frac{\|g_{i}\|^2}{2} +
        \frac{1}{2{-}\gamma L {-} \gamma \mu}\frac{\|g_{N}\|^2}{2},
    \end{aligned}
$$%
which is of type \eqref{eq:I^k_N_general_ineq} since all gradient norms have positive coefficients because $\TigLmu[k]\geq 0$ and $\TigLmu[k+1]<0$. By taking the minimum gradient norm:%
$$
     \min\limits_{k \leq i \leq N} \Big\{\frac{\|g_i\|^2}{2L}\Big\}
            {}\leq{}
        \frac{f_{0}-f_{N}}{ \gamma L \, \PN },
$$%
where
$$
\begin{aligned}
     \mathsmaller{\PN}
        {}{=}{}
     \tfrac{-(1{-}\gamma L)^{2}(1{-}\gamma \mu)^{2}\TigLmu[k+1] + \TigLmu[k] + \gamma L {-} \gamma \mu}{(1{-}\gamma \mu)^{2} {-} (1{-}\gamma L)^{2}}  {}{+}{}
        \tfrac{(2{-}\gamma L)(2{-}\gamma \mu)}{2{-}\gamma L {-} \gamma \mu} \mathsmaller{(N{-}k{-}1)}.
    \end{aligned}
$$
Using algebraic manipulations one can check the following equivalent expressions for $\PN$ discussed in \cref{subsec:rates_hypoconvex}:%
$$
\begin{aligned}
P_N(\gamma L, \gamma \mu)
{}={}&
        \frac{ \frac{{-1} +(1-\gamma L)^{-2k}}{\gamma L [1-(1-\gamma L)^2]} - \frac{{-1} +(1-\gamma \mu)^{-2k}}{\gamma \mu [1-(1-\gamma \mu)^2]}}{ \frac{1}{1-(1-\gamma L)^2} - \frac{1}{1-(1-\gamma \mu)^2} } 
{}+{}
\tfrac{(2-\gamma L)(2-\gamma \mu)}{2-\gamma L - \gamma \mu} (N - k); \\
%%%%%%%%%%%%%%%%%%%%%%%%%%%%%%%%%%%%%%%%%
P_N(\gamma L, \gamma \mu)
{}={}&
\tfrac{(2-\gamma L)(2-\gamma \mu)}{2-\gamma L - \gamma \mu} 
\Bigg[
N {}-{} \tfrac{-\gamma \mu \gamma L}{\gamma L - \gamma \mu}
        \bigg[ 
 \tfrac{\frac{{-1} +(1-\gamma \mu)^{-2k}}{1-(1-\gamma \mu)^2 } - k}{\gamma \mu} - 
 \tfrac{\frac{{-1} +(1-\gamma L)^{-2k}}{1-(1-\gamma L)^2} - k}{\gamma L}
 \bigg]
\Bigg].
\end{aligned}
$$%
The expression from \eqref{eq:P_den_full_rate_h_ct} is obtained from the identity:
$$
\tfrac{\frac{{-1} +(1-\gamma \mu)^{-2k}}{1-(1-\gamma \mu)^2 } - k}{\gamma \mu} 
{}-{}
\tfrac{\frac{{-1} +(1-\gamma L)^{-2k}}{1-(1-\gamma L)^2 } - k}{\gamma L} 
{}={}
\sSum_{i=0}^{k} \left[\tfrac{{-1} +(1-\gamma \mu)^{-2i}}{\gamma \mu} - \tfrac{{-1} +(1-\gamma \mu)^{-2i}}{\gamma \mu} \right],
$$%
where in the \rhs are summed up quantities $\TigLmu[i]$ with $i=\{0,1,\dots,k\}$, i.e., for all indices for which they are positive. Therefore, the \rhs can be rewritten using the positive quantities over $N$ terms such as in formula \eqref{eq:P_den_full_rate_h_ct}.%
\item \textbf{Case $\gamma L \in [\gbari[N], 2)$.} From \cref{lemma:general_lemma_ct_h_to_Linear} with $k=N-1$ it holds:%
\begin{align*}
   \frac{f_0 - f_N}{\gamma}
        {}\geq{} & 
    \Ei[N][1-\gamma L] \frac{\|g_N\|^2}{2} + 
    \frac{(1-\gamma \mu)^2 \, \TigLmu[N]}{(\gamma L - \gamma \mu)^2} \, 
    \frac{\| g_N - (1-\gamma L) g_{N-1} \|^2}{2},
\end{align*}%
where the mixed term can be neglected since $\TigLmu[N] \geq 0$ and obtain:%
$$
    f_0 - f_N
        {}\geq{}
    \big[ {-1} {}+{} (1-\gamma L)^{-2N}\big] \frac{\|g_{N}\|^2}{2L},
$$%
which is of type \eqref{eq:I^k_N_general_ineq} and implies the linear regime from \cref{thm:wc_GM_hypo_ct_all}. \qed
\end{enumerate}
\end{proof}%
\subsection{Proofs of performance bounds for strongly convex functions}\label{subsec:proof_thm_ct_strg_cvx}
If $\mu > 0$, then $\gbari[\infty] = \frac{2}{1+\kappa} < 2$.
\begin{proof}[\textbf{\cref{thm:strg_cvx_rate}}]\label{proof:thm:strg_cvx_rate}%
We split the proof over stepsize intervals:%
\[
    (0,2) 
        {}={}
    \big(0, \gbari[1]\big] \cup
    \Big( \mathop{\cup}\limits_{k=1}^{N-1} \big[\gbari[k],\gbari[k+1]\big)\Big)
    {}\cup{} \big[\gbari[N],\tfrac{2}{1+\kappa}\big) 
    {}\cup{} \big[\tfrac{2}{1+\kappa}, 2\big).
\]%
\begin{enumerate}[align=right,itemindent=1em,labelsep=2pt,labelwidth=1em,leftmargin=0pt,nosep]%
\item \textbf{Case $\gamma L \in \big(0, \gbari[1]\big]$}. By telescoping \eqref{eq:dist_1_mu_1_step_decrease} for indices $i=0,\dots,N-1$:%
$$
        \frac{f_{0} - f_{N}}{\gamma} 
            {}\geq{}
        \Ei[N][\m] \frac{\|g_{N}\|^2}{2} {}+{}        
            \sum_{i=0}^{N-1}
            \frac{ {-1} + (\m+\n) \Ei[i+1][\m] }{\m-\n} \frac{\|g_{i+1}-\m g_{i}\|^2}{2},
$$%
where all mixed terms have nonnegative coefficients because $\gamma L \in \big(0,  \gbari[1]\big]$ (see \cref{lemma:dist_1_mu_1_step_decrease}). Hence they can be neglected to obtain the lower bound of type \eqref{eq:I^k_N_general_ineq}:
\begin{align}\label{eq:sufficient_decrease_N_its_mu_strongly_convex_towards_rate_no_dist_2}
    \frac{f_0 - f_N}{\gamma}
        {}\geq{}
    \Ei[N][\m] \frac{\|g_{N}\|^2}{2}.
\end{align}%:
\item \textbf{Case $\gamma L \in \big[\gbari[k],\gbari[k+1]\big)$}, with some $k \in \{1,\dots,N-1\}$. By telescoping \eqref{eq:dist_1_mu_1_step_decrease} for indices $i=k+1,\dots,N-1$:%
$$
\begin{aligned}
        \frac{f_{k+1} {-} f_{N}}{\gamma} 
            {}\geq{} &
        {-}\Ei[k+1][\m] \frac{\|g_{k+1}\|^2}{2} {+}
        \Ei[N][\m] \frac{\|g_{N}\|^2}{2} {+} \\
        & \quad     \sum_{i{=}k{+}1}^{N{-}1}
            \frac{ {-1} {+} (\m{+}\n) \Ei[i{+}1][\m] }{\m{-}\n} \frac{\|g_{i+1}{-}\m g_{i}\|^2}{2},
\end{aligned}
$$%
with nonnegative weights of the mixed terms since $\gamma L \in \big(0,  \gbari[k+1]\big)$ (see \cref{lemma:dist_1_mu_1_step_decrease}). After appending it to \eqref{eq:general_ineq_proof_to_linear_mu_compact} from \cref{lemma:general_lemma_ct_h_to_Linear_mu}, we obtain:%
$$
\begin{aligned}
    \frac{f_0 {-} f_N}{\gamma}
        {}\geq{}
    \Ei[N][\m] \frac{\|g_{N}\|^2}{2} +
    \frac{-\n^2 \, \Tinm[k+1]}{(\m-\n)^2} \, 
    \frac{\| g_{k+1} - \m g_{k} \|^2}{2} +     \\
    \sum_{i=k+1}^{N-1}
            \frac{ {-1} + (\m+\n) \Ei[i+1][\m] }{\m-\n} \frac{\|g_{i+1}-\m g_{i}\|^2}{2}.
\end{aligned}
$$
Since $\Tinm[k+1] < 0$, all mixed terms can be neglected to get a \eqref{eq:I^k_N_general_ineq}-type inequality:%
\begin{align}\label{eq:sufficient_decrease_N_its_mu_strongly_convex_towards_rate}
    \frac{f_0 - f_N}{\gamma}
        {}\geq{}
    \Ei[N][\m] \frac{\|g_{N}\|^2}{2}.
\end{align}%
 \item \textbf{Case $\gamma L \in \big[\gbari[N],\frac{2}{1+\kappa}\big)$.} From \cref{lemma:general_lemma_ct_h_to_Linear} with $k=N-1$ it holds:%
$$
   \frac{f_0 - f_N}{\gamma}
        {}\geq{}
    \Ei[N][\n] \frac{\|g_N\|^2}{2} + 
    \frac{\m^2 \, \Tinm[N]}{(\m-\n)^2} \, 
    \frac{\| g_N - \n g_{N-1} \|^2}{2}.
$$%
Since $\TigLmu[N] \geq 0$, by neglecting the mixed term we obtain a \eqref{eq:I^k_N_general_ineq}-type inequality:
\begin{align}\label{eq:wc_linear_L}
    \frac{f_0 - f_N}{\gamma}
        {}\geq{}
    \Ei[N][\n] \frac{\|g_{N}\|^2}{2}.
\end{align}%
\item \textbf{Case $\gamma L \in \big[\frac{2}{1+\kappa},2\big)$.} By telescoping \cref{lemma:dist_1_L_1_step_decrease} for $i=0,1,\dots,N-1$ it holds:%
$$
    \frac{f_{0} - f_{N}}{\gamma} 
                {}\geq{}
            \Ei[N][\n] \frac{\|g_{i}\|^2}{2} +
            \sum_{i=0}^{N-1} \frac{ 1 - (\m+\n) \Ei[i+1][\n] }{\m-\n} \frac{\|g_{i+1}-\n g_{i}\|^2}{2},
$$%
where the mixed terms can be neglected and then the inequality reduces to \eqref{eq:wc_linear_L}.%
\end{enumerate}%
Since the transition between the two linear regimes is given by the sign of $\Tinm[N]=\Ei[N][\m]{-}\Ei[N][\n]$, we can glue together expressions \eqref{eq:sufficient_decrease_N_its_mu_strongly_convex_towards_rate_no_dist_2}, \eqref{eq:sufficient_decrease_N_its_mu_strongly_convex_towards_rate} and \eqref{eq:wc_linear_L}:%
$$
    \frac{f_0 - f_N}{\gamma}
        {}\geq{}
    \min\left\{ \Ei[N][\m] \,,\, \Ei[N][\n]\right\} \frac{\|g_{N}\|^2}{2},
$$%
equivalent to:
$$
    f_0 - f_N
        {}\geq{}
    \min
    \Big\{ \frac{{-1} +(1-\gamma \mu)^{-2N}}{\mu} \,,\, 
    \frac{{-1} +(1-\gamma L)^{-2N}}{L}
    \Big\} \frac{\|g_{N}\|^2}{2},
$$%
which implies the bound \eqref{eq:strg_cvx_full_ct} and the complementary one:
$$
\begin{aligned}
    \frac{1}{2L} \|\nabla f(x_N)\|^2
        {} \leq {}  
    \frac{f(x_0)-f(x_N)}
    { \min\big\{ \frac{L}{\mu} \left[{-1}+(1-\gamma \mu)^{-2N}\right] \,,\, {-1}+(1-\gamma L)^{-2N} \big\} }. {}\qquad{}\qed 
\end{aligned}%
$$%
\end{proof}%
\subsection{Results for convex functions}\label{subsec:cvx_case_comments}
\begin{corollary}[Multistep descent for convex functions] \label{lemma:G4SD_sufficient_decrease_cvx}
Let $f \in \mathcal{F}_{0,L}$. Consider $k+1 \geq 1$ iterations of \eqref{eq:GM_it_fixed_steps} with stepsize $\gamma L \in [\gbari[k], \gbari[k+1])$, satisfying $\gamma L > 1$, starting from $x_0$. Then the following inequality holds, with positive coefficients of gradient norms:%
    \begin{align}
        f(x_0) {-} f(x_{k+1}) 
            {}{\geq}{}
        \Big[2(k+1)\frac{{-}(1-\gamma L)^2}{2{-}\gamma L} + \sum_{i{=}0}^{k} (1{-}\gamma L)^{{-}2i} \Big] \frac{\|\nabla f(x_{k})\|^2}{2L} {}+{}
        \nonumber \\ {}+{}
        \Big[2(k+1)\frac{1}{2{-}\gamma L} - \sum_{i=0}^{k} (1{-}\gamma L)^{-2i} 
        \Big] \frac{\|\nabla f(x_{k+1})\|^2}{2L}. \label{eq:GN4SD_convex}\tag{G4SD}
    \end{align}%
Moreover, \eqref{eq:GN4SD_convex} is valid on the extended stepsize range $\gamma L \in (1, 2)$.
\end{corollary}
\begin{proof}
    It results by taking $\m = 1$ in \eqref{eq:general_ineq_proof}. \qed
\end{proof}
\eqref{eq:GN4SD_convex} is a multistep generalization of {\shownewtext{the descent lemma}} \eqref{eq:4SD_Teboulle} from \cite{Teboulle2022}.%
\begin{proof}[\textbf{\cref{corollary:cvx_rate}}]\label{proof:corollary:cvx_rate} 
Letting $\mu \nearrow 0$, all terms summed in \eqref{eq:P_den_full_rate_h_ct} become zero, simplifying to $P_N(\gamma L,0) {}={} 2N$, so that all sublinear regimes corresponding to stepsizes up to $\gbari[N][(0)]$ share the same expression. For stepsizes $\gamma L \geq \gbari[N][(0)]$, it holds the linear regime following the expression in the denominator $\Ei[N][1-\gamma L] = \frac{{-1} +(1-\gamma L)^{-2N}}{\gamma L}$. Additionally, the thresholds are determined by the roots of:%
$$
    \TigLmu[k][\gamma L][0]
        {}={}
    2k - \frac{{-1} +(1-\gamma L)^{-2k}}{\gamma L}, \,\, \forall k=1,\dots,N,
$$
hence the rate can be expressed as a minimum in the denominator of \eqref{eq:cvx_full_ct}. Same bound is recovered to the limit $\mu \searrow 0$ in strongly convex functions' rate \eqref{eq:strg_cvx_full_ct}. 
\qed \end{proof}%%
%
\section{Tightness of performance bounds}\label{sec:tightness}
    % \section{Tightness of the results} \label{subsec:tightness}
The tightness of convergence rates is typically assessed by identifying a worst-case function matching the upper bound when iterating the algorithm. This function can be determined either analytically or by providing an interpolable set of triplets. \cref{table:summary_bounds_and_tightness} summarizes the various demonstrations of tightness based on different stepsizes ranges.%

\begin{remark}
    The exactness of all proved convergence rates can be confirmed numerically by solving the associated performance estimation problems (PEPs), for example by using the specialized software packages PESTO \cite{PESTO} (for Matlab) or PEPit \cite{PESTO_python} (for Python). In this section we propose analytical proofs.
\end{remark}
For stepsizes in the range $\gamma L \in [\gbari[N], 2)$, a single performance bound applies across all function classes (weakly convex, convex and strongly convex).
As shown in \cref{prop:tightness_linear_regime}, this unified bound is proven to be tight using the same worst-case function $L \frac{\|x\|^2}{2}$, which, for simplicity, can be reduced to one dimension.%

\begin{proposition}[Tightness for $\gamma L \in {[} {\gbari[N]}, 2 )$]\label{prop:tightness_linear_regime}
Let $L > 0$, $N$ be a positive integer, $\gamma L \in [\gbari[N],2)$ and $f_L(x) {}\coloneqq{} L \frac{\|x\|^2}{2}$. Then after applying $N$ iterations of \eqref{eq:GM_it_fixed_steps} on function $f=f_L$ with constant stepsize $\gamma L$, starting from $x_0$, it holds:
$$
    \tfrac{1}{2L} \min_{0 {} \leq {} i {} \leq {} N} \{ \|\nabla f(x_i)\|^2\}  
        {}={}
    \tfrac{1}{2L} \|\nabla f(x_N)\|^2
        {}={}
    \frac{f(x_0) - f_*}{(1-\gamma L)^{-2N}}.
$$
\end{proposition}%
\begin{proof}%
    One can directly check that $f_L$ satisfies the necessary conditions for equality in \cref{lemma:general_lemma_ct_h_to_Linear}: 
    $g_N=(1-\gamma L)g_{N-1}$ and $f_0 - f_N = [-1 + (1-\gamma L)^{-2N}] \frac{\|g_N\|^2}{2L}$.
\qed\end{proof}%

For the remainder of the stepsize interval, $\gamma L \in (0, \gbari[N])$, the proofs are organized as follows: 
\cref{subsec:tightness_proofs_convex_strg_convex} covers constant stepsizes for convex and strongly convex functions;
\cref{subsec:tightness_hypoconvex_up_to_thr} addresses variable stepsizes below $\gbari[1]$ for convex and weakly convex functions;
\cref{subsec:tightness_hypo_above_thr} focuses on the specific range of constant stepsizes $(1, \gbari[N])$ for weakly convex functions.% 
\subsection{Tightness for convex and strongly convex functions}\label{subsec:tightness_proofs_convex_strg_convex}
\begin{proposition}[Tightness of \cref{thm:strg_cvx_rate}]\label{prop:Tightness_strongly_convex}   
Let $L > 0$, $\mu \in (0,L)$, $\gamma L \in (0,2)$ and $N$ be a positive integer. 
Then there exists $x_0$ and $f \in \mathcal{F}_{\mu,L}$ such that after applying $N$ iterations of \eqref{eq:GM_it_fixed_steps} on function $f$ with constant stepsize $\gamma L \in (0,2)$, starting from $x_0$, it holds:%
$$
    \tfrac{1}{2L} \|\nabla f(x_N)\|^2
        {} = {}  
    \frac{f(x_0)-f_*}
    { 1 + \gamma L \min \left\{
    \frac{-1+(1-\gamma \mu)^{-2N}}{\gamma \mu}
    \,,\,
    \frac{-1+(1-\gamma L)^{-2N}}{\gamma L}
         \right\}}.%
$$%
\end{proposition}
%%%%%%%%%%%%%%%%%%%%%%%%%%%%%%%%%%%%%%%%%%% OLD PROOF TO CANCEL %%%%%%%%%%%%%%%%%%%%%%%%%%%%%%%%%%%
\begin{proof}
We define the following worst-case function depending on the stepsizes $\gamma L$: 
if $\gamma L \in [\gbari[N], 2)$, then it is given by the function $f_L$ (defined in \cref{prop:tightness_linear_regime}); 
if $\gamma L \in (0, \gbari[N])$, then it is the function $\varphi$ given below. 
Let $\Delta > 0$ and define the parameter 
$\tau {}{}\coloneqq{}{} \sqrt{ \frac{2}{L} \, \frac{\Delta}{1{+}\frac{L}{\mu} [{-}1{+}(1{-}\gamma \mu)^{-2N}]}}$. The function $\varphi\colon \mathbb{R} \rightarrow \mathbb{R}$ is given by:%
    \begin{align*}
        \varphi(x) 
            {}{}\coloneqq{}{}
        \left\{
        \arraycolsep=6pt
        \def\arraystretch{1.25}
            \begin{array}{ll}
                \mu \frac{x^2}{2} + (L-\mu) \tau |x| -(L-\mu)\frac{\tau^2}{2}, & \text{ if } |x| \geq \tau; \\
                L \frac{x^2}{2}, &  \text{ if } |x| < \tau.
            \end{array}
        \right.
    \end{align*}%
Note that $\varphi \in \mathcal{F}_{\mu,L}$, $x_{*} =\argmin \varphi(x) = 0$, $\varphi(x_*) = 0$ and
$$
    \nabla \varphi(x) 
        {}{}\coloneqq{}{}
    \left\{
    \arraycolsep=6pt
    \def\arraystretch{1.25}
        \begin{array}{ll}
            \mu x + (L-\mu) \tau \sgn(x), & \text{ if } |x| \geq \tau; \\
            L x, &  \text{ if } |x| < \tau.
        \end{array}
    \right.
$$% 
One can check that initializing the iterates at $x_0 {\coloneqq} \tau [ 1 {+} \frac{L}{\mu} \left({-1} {+} (1{-}\gamma \mu)^{-N} \right) ]$, with $\varphi(x_0) = \Delta$, and running $N$ iterations of \eqref{eq:GM_it_fixed_steps}, the following holds: 
all iterations belong to the branch $x \geq \tau$,
the final iterate is exactly $x_N = \tau$, having the gradient $\nabla \varphi(x_N) = L \tau$ and thus reaching the bound claimed in the first term.%
    {}\qed{}%
\end{proof}%
This one-dimensional worst-case example is inspired from \cite[Section 4.1.2]{taylor_smooth_2017}.%
\begin{proposition}[Tightness of \cref{corollary:cvx_rate}]\label{corr:convex_tightness}%
Let $L > 0$, $\gamma L \in (0,2)$ and $N$ be a positive integer. 
Then there exists $x_0$ and $f \in \mathcal{F}_{0,L}$ such that after applying $N$ iterations of \eqref{eq:GM_it_fixed_steps} on function $f$ with constant stepsize $\gamma L \in (0,2)$, starting from $x_0$, it holds:%
$$ 
    \tfrac{1}{2L} \|\nabla f(x_N)\|^2
        {} = {}  
    \frac{f(x_0)-f_*}
    { \min \left\{
    1 + 2N \gamma L
    \,,\,
    (1-\gamma L)^{-2N}
         \right\}}.%
$$%
\end{proposition}%
\begin{proof}
With $\Delta > 0$, let $\tau {}{}{\coloneqq}{}{} \sqrt{\frac{2}{L} \, \frac{\Delta}{1{+}2N \gamma L}}$ and $x_0 = \tau(1{+}N \gamma L)$. We define $\varphi{\colon}\mathbb{R} \rightarrow \mathbb{R}$, such that $\varphi(x_0) = \Delta$, $\varphi(x_*) = 0$ and %
$$
\begin{aligned}
        \varphi(x) 
            {}{}\coloneqq{}{}&
        \left\{
        \arraycolsep=4pt
        \def\arraystretch{1.25}
            \begin{array}{ll}
                L \tau |x| - L \frac{\tau^2}{2}, & \text{ if } |x| \geq \tau; \\
                L \frac{x^2}{2}, &  \text{ if } |x| < \tau,
            \end{array}
        \right.
            {}\,\,{}
         \nabla \varphi(x) 
        {}{}\coloneqq{}{}
    \left\{
    \arraycolsep=4pt
    \def\arraystretch{1.25}
        \begin{array}{ll}
            L \tau \sgn(x), & \text{ if } |x| \geq \tau; \\
            L x, &  \text{ if } |x| < \tau.
        \end{array}
    \right.
\end{aligned}%
$$%
Then we select $f = \varphi$ or $f = f_L$ (see \cref{prop:tightness_linear_regime}) whether $\gamma L < \gbari[N][(0)]$ or $\gamma L \geq \gbari[N][(0)]$, respectively. \qed
\end{proof}%
\subsection{Removing the optimal point in tightness analysis for smooth weakly convex functions}
\cref{lemma:characterization_min} shows that by incorporating inequalities \eqref{eq:sufficient_decrease_N_*} for all iterates $i \in \mathcal{I}$ the upper bounds remain tight and provides an existence guarantee of an interpolating function with global minimum $f_*$. It generalizes Theorem 7 from \cite{drori2021complexity}, obtained for the particular choices $\mu=-L$ and $\mu=0$, to any $\mu \leq 0$.%
\begin{lemma}\label{lemma:characterization_min}%
Let $\mathcal{T}=\{(x_i, g_i, f_i)\}_{i \in \mathcal{I}}$ be $\mathcal{F}_{\mu,L}$-interpolable, with $\mu \in (-\infty, 0]$. There exists at least one interpolating function with a finite global minimum
\begin{align}\label{eq:Charact_f*}
  f_* = \min_{x\in\mathbb{R}^d} f(x) = 
  \min_{i\in \mathcal{I}} \big\{ f_i - \tfrac{1}{2L} \|g_i\|^2 \big\}.
\end{align}
A global minimizer is given by 
$x_*=x_{i_*} - \tfrac{1}{L}g_{i_*}$, where $
i_* \in \arg \min\limits_{i\in\mathcal{I}} \big\{f_i - \frac{1}{2L} \|g_i\|^2 \big\}.
$
\end{lemma}%
\begin{proof}
    The proof consists on algebraic manipulations, see \cref{appendix:proof_charact_min}.
\qed\end{proof}%
\begin{lemma}\label{lemma:decouple_f*}
Assume an $\mathcal{F}_{\mu,L}$-interpolable set $\mathcal{T}=\{(x_i, g_i, f_i)\}_{i \in \mathcal{I}}$ satisfying the condition%
\begin{align}\label{eq:tightness_condition_N=i_*}
    N \in \argmin_{i\in \mathcal{I}} \big\{ f_i - \tfrac{1}{2L} \|g_i\|^2  \big \}.
\end{align}
Then proving tightness of the bound
$$
    \tfrac{1}{2L} \min_{0 {} \leq {} i {} \leq {} N} \{ \|\nabla f(x_i)\|^2\} 
        {} \leq {}
    \frac{f(x_0)-f_*}{1 + \PN}
$$%
reduces to showing tightness of the bound%
\begin{align}\label{eq:rate_fN_on_tightness}
    \tfrac{1}{2L} \min_{0 {} \leq {} i {} \leq {} N} \{ \|\nabla f(x_i)\|^2\} 
        {} \leq {}
    \frac{f(x_0)-f(x_N)}{\PN}.
\end{align}%
\end{lemma}%
\begin{proof}
    Condition \eqref{eq:tightness_condition_N=i_*} implies that \eqref{eq:Charact_f*} from \cref{lemma:characterization_min} holds with $N=i_*$, hence there exists a function for which $f(x_N) - f_* {}\geq{} \frac{1}{2L} \|\nabla f(x_N)\|^2$. This inequality produces the shift in the denominator when appended to rate \eqref{eq:rate_fN_on_tightness}.%
\qed\end{proof}%
Further on, in all tightness proofs we show that condition \eqref{eq:tightness_condition_N=i_*} holds and therefore we reduce the analysis to only prove the equality in \eqref{eq:rate_fN_on_tightness}.%
\subsection{Tightness analysis for weakly convex functions and stepsizes below \texorpdfstring{$\gbari[1][]$}{first threshold}}\label{subsec:tightness_hypoconvex_up_to_thr}%
\begin{figure}%
    \centering    
    \begin{subfigure}[t]{0.49\textwidth}
        \centering
        \includegraphics[width=\textwidth]{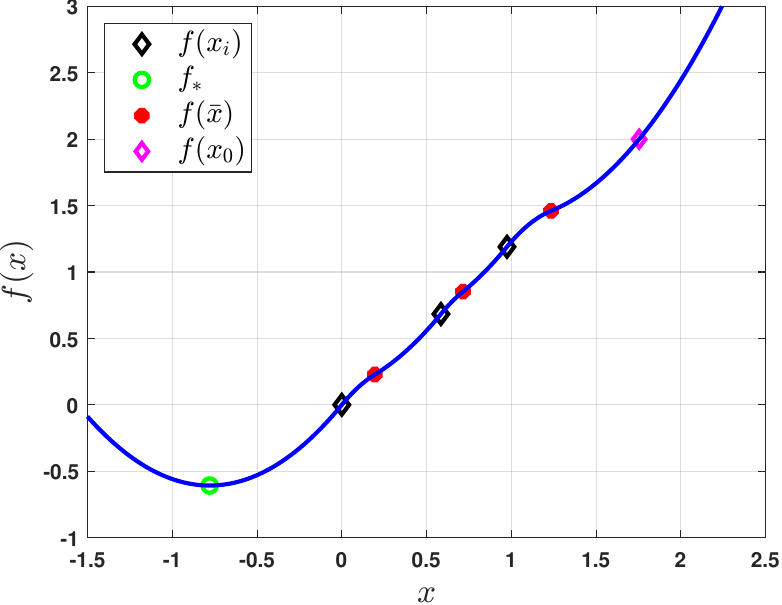}
        \caption{One-dimensional worst-case function example corresponding to short stepsizes $\gamma_i L \in (0, 1]$ (piecewise quadratic; see \cref{prop:tightness_short_steps}). Setup: $N=3$, $f_0-f_*=2$, $L=2$, $\mu=-4$, $\gamma_0 L=1$, $\gamma_1 L=0.5$, $\gamma_2 L=0.75$.}%
        \label{fig:func_example_short_steps_1D}
    \end{subfigure}
    \hfill
    \begin{subfigure}[t]{0.49\textwidth}
        \centering
        \includegraphics[width=1.1\textwidth]{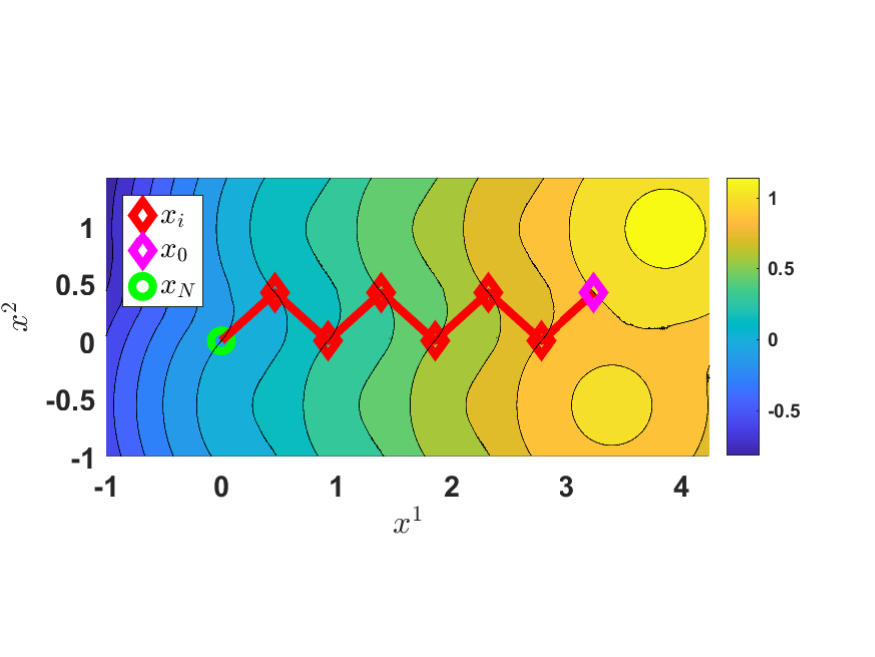}
        \caption{Function example for constant stepsize $\gamma L\in(1, \gbari[1]]$. Setup: $L=1$, $\mu=-0.5$, $\gamma L=1.5$, $N=7$. Component $x^1$ is successively reduced, while $x^2$ alternates between two values, due to alternating gradients $g_{i}$ and $g_{i+1}$.%
        }
        \label{fig:func_example_long_steps_2D}
    \end{subfigure}    
    \caption{%
    Worst-case examples for \cref{thm:wc_GM_hypo_ct_all,thm:wc_GM_hypo} with $\gamma L \leq \gbari[1]$.
    }\label{fig:wc_examples_1D_2D_hypo}%
\end{figure}

In this section, we construct worst-case examples for $\gamma L \leq \gbari[1]$. For a given choice of parameters, we illustrate these worst-case functions in \cref{fig:wc_examples_1D_2D_hypo}.%
\paragraph*{\textbf{Range $\gamma_i L \in (0,1]$.}} We show a one dimensional worst-case function example which extends to weakly convex functions the one from \cite[Proposition 4]{abbaszadehpeivasti2021GM_smooth}.%
\begin{proposition}[Tightness for $\gamma L \in (0, 1{]}$]\label{prop:tightness_short_steps}
Let $L > 0$, $\mu \in (-\infty,0)$, $N$ be a positive integer and $\gamma_i L \in (0,1]$, with $i=0,\dots,N-1$. Then there exists $x_0$ and $f \in \mathcal{F}_{\mu,L}$ such that after applying $N$ iterations of \eqref{eq:GM_it_fixed_steps} on function $f$ with stepsizes $\gamma_i L$, starting from $x_0$, it holds:%
$$
    \tfrac{1}{2L} \min_{0 {} \leq {} i {} \leq {} N} \{ \|\nabla f(x_i)\|^2\}  {} = {}
        \frac{f(x_0)-f_*}{%
        1 + \sum_{i=0}^{N-1} \gamma_i L \big(2 + \frac{\gamma_i L \gamma_i \mu}{\gamma_i L-\gamma_i \mu}\big)
        }.
$$
\end{proposition}%%
%%%
\begin{proof}
Motivated by the necessary conditions $g_{i}=g_{i+1}$ in \cref{lemma:sufficient_decrease_h_leq_1}, a one dimensional function with constant gradients $U$ at the iterates is built, where $U$ is the square root of the minimum gradient norm:
$
    U^2 {}{}{\coloneqq}{}{} { \frac{\Delta}{ 
        \frac{1}{2L}\sum_{i{=}0}^{N{-}1} \gamma_i L ( 2 {+} \frac{\gamma_i L \gamma_i \mu}{\gamma_i L- \gamma_i \mu} ) } },
$
with $\Delta>0$. 
The equality case in \cref{lemma:sufficient_decrease_h_leq_1} implies (together with $g_i=U,\, {}\forall{} i=0,1,\dots,N$):
$$
    \frac{f_{i}-f_{i+1}}{\Delta}
        {} = {}
    \frac{\gamma_i L \big( 2  + \frac{\gamma_i L \gamma_i \mu}{\gamma_i L - \gamma_i \mu} \big)}
        { \sum_{j=0}^{N-1} \gamma_j L \big( 2  + \frac{\gamma_j L \gamma_j \mu}{\gamma_j L- \gamma_j \mu} \big) }
    {} \Leftrightarrow {}
    \frac{f_{i}}{\Delta}
        {} = {}
    \frac{\sum_{j=i}^{N-1} \gamma_j L \big( 2  + \frac{\gamma_j L \gamma_j \mu}{\gamma_j L - \gamma_j \mu} \big)}
          {\sum_{j=0}^{N-1} \gamma_j L \big( 2  + \frac{\gamma_j L \gamma_j \mu}{\gamma_j L - \gamma_j \mu} \big)}.
$$
Without loss of generality, we set $f_N = 0$, $f_0=\Delta$ and $x_N=0$ and define the triplets $\big\{(x_i, g_i, f_i)\big\}_{i \in \{0,\dots,N\}}$, with $g_i=U$ and $x_i {}={} U \sum_{j=i}^{N-1} \gamma_j$. 
Consider the following points lying between consecutive iterates:%
$$
    \bar{x}_i {}\coloneqq{} x_i - \tfrac{-\mu}{L-\mu} \gamma_i U \in [x_{i+1},x_{i}], {}\quad{} i = 0,1,\dots,N-1.%
$$%
A worst-case function interpolating the triplets is the piecewise quadratic $f\colon\mathbb{R}\rightarrow \mathbb{R}$, with alternating curvature between $\mu$ and $L$ with inflection points $x_i$ and $\bar{x}_i$, constant gradients at the iterations and decreasing function values, where $i \in \{0,\dots,N-1\}$:%
\begin{align}\label{eq:wc_func_iterates_no_opt}
% \hspace{-.5em}
    f(x)
        {}={} 
    \left\{
    \arraycolsep=6pt
    \def\arraystretch{1.25}
    \begin{array}{ll}
        \frac{L}{2} (x-x_N)^2 + U (x - x_N) + f_N,
        & 
        \text{ if }x \in \left(-\infty, x_N\right]; \\
        \frac{\mu}{2} \left(x - x_{i+1} \right)^2 + U \left(x - x_{i+1} \right) + f_{i+1},
        &
        \text{ if }x \in [x_{i+1}, \bar{x}_i];
        \\
        \frac{L}{2} \left(x - x_{i} \right)^2 + U \left(x - x_{i}   \right) + f_{i}, & 
        \text{ if }x \in [\bar{x}_i, x_i]; \\
        \frac{L}{2} \left(x-x_0\right)^2 + U \left(x - x_{0} \right) + f_{0},
        &
        \text{ if }x \in \left[x_{0}, \infty\right). {}\quad{}\qed
    \end{array}
    \right.%
\end{align}%
\end{proof}%
\shownewtextfinal{An illustration of the one-dimensional worst-case construction from \cref{prop:tightness_short_steps} is given in \cref{fig:func_example_short_steps_1D}.}%
\begin{corollary}[Tightness of \cref{corr:cvx_rates_variable_stepsizes} for convex functions] \label{corr:convex_tightness_variable_stepsizes}%
Let $L > 0$, $N$ be a positive integer and $\gamma_i L \in (0,\frac{3}{2}]$, with $i=0,\dots,N-1$. Then there exists $x_0$ and $f \in \mathcal{F}_{0,L}$ such that after applying $N$ iterations of \eqref{eq:GM_it_fixed_steps} on function $f$ with stepsizes $\gamma_i L$, starting from $x_0$, it holds:%
$$
    \tfrac{1}{2L} \|\nabla f(x_N)\|^2
        {} = {} 
    \frac{f(x_0)-f_*}{1 + 
    2\sum_{i=0}^{N-1} \gamma_i L}.%
$$%
\end{corollary}%
\begin{proof}%
The proof results by setting $\mu=0$ in \cref{prop:tightness_short_steps} to obtain a linear function between iterations, extended outside by a quadratic of curvature $L$. \qed%
\end{proof}%
\paragraph*{\textbf{Range $\gamma_i L \in (1, \gbari[1]]$. }} \,We build a set of triplets matching the corresponding bounds from \cref{thm:wc_GM_hypo_ct_all,thm:wc_GM_hypo} and interpolating an $\mathcal{F}_{\mu,L}$-function, implying the validity of all interpolation inequalities relating the $N+1$ iterations. The upper bounds in \cref{thm:wc_GM_hypo} are obtained by analyzing only consecutive iterations, from which we track equality conditions collected in \cref{prop:nec_eq_conds_variable_gamma}.%
\begin{proposition}\label{prop:nec_eq_conds_variable_gamma}
Let $L > 0$, $\mu \in (-\infty,0)$, $N$ be a positive integer and $\gamma_i L \in (1,\gbari[1])$, with $i=0,\dots,N-1$. Then any worst-case function $f \in \mathcal{F}_{\mu,L}$ matching the bound from \cref{thm:wc_GM_hypo} satisfies the following identities on the triplets $\mathcal{T} {}={} \{(x_i, g_i, f_i)\}_{i\in\{0,\dots,N\}}$, with $x_{i+1}
        {}={}
    x_i - \gamma_i g_i$:%
\begin{align} 
&  \|g_i\|^2 
        {}={} 
    U^2 
        {}{}\coloneqq{}{}
    \frac{f_0 - f_N}{\frac{1}{2L} \sum_{j=0}^{N-1} \frac{\gamma_j L (2-\gamma_j L)(2-\gamma_j \mu)}{2-\gamma_j L - \gamma_j \mu}}, {}\quad{} \forall i = 0,\dots,N; \label{eq:defs_U_tightness_long_h_ct} \\
%%%%%%%%%%%%%%%%%%%%%%%%%%%%%%%%%%%%%%%%%%%%%%%%%%%%%%%%%%%%%%%%%%%
& \cgLgmu[\gamma_i]
        {}{}\coloneqq{}{}
    \frac{\left\langle g_{i} , g_{i+1} \right\rangle}{U^2}
        {}={} 
    \frac{1 + (1-\gamma_i \mu)(1-\gamma_{i} L)}{2-\gamma_i L - \gamma_i \mu}, {}\quad{} \forall i = 0,\dots,N-1;
 \label{eq:def_ci} \\
%%%%%%%%%%%%%%%%%%%%%%%%%%%%%%%%%%%%%%%%%%%%%%%%%%%%%%%%%%%%%%%%%%%
&f_i
        {}={}
    f_N + 
    \frac{U^2}{2L} \sum_{j=i}^{N-1} \frac{\gamma_j L (2-\gamma_j \mu)(2-\gamma_{j} L)}{2-\gamma_j L - \gamma_j \mu}, {}\quad{} \forall i = 0,\dots,N.
    \label{eq:cond_fi_medium_steps}%
\end{align}%
In particular, a two-dimensional function $f\colon \mathbb{R}^2\rightarrow \mathbb{R}$ satisfies these necessary conditions if, for any suitable $g_0$, the following holds:%
\begin{align}\label{eq:def_reccurence_gi_medium_steps}
    g_{i+1} 
        {}={} 
    R\big( (-1)^{i+1} \thetagLgmu[\gamma_i] \big) g_{i}, {}\quad{} \forall i = 0,1,\dots,N-1,
\end{align}%
where $R$ is the two-dimensional rotation matrix and%
\begin{align}\label{eq:def_theta_i}
    \thetagLgmu[\gamma_i] 
        {}{}\coloneqq{}{} 
    \arccos \big( \cgLgmu[\gamma_i] \big). 
\end{align}%
Moreover, %
\begin{align}\label{eq:cond_*_mid_steps}
    N{}\in{}\argmin_{0 \leq i \leq N} \big\{f_i - \tfrac{1}{2L}\|g_i\|^2\big\}.
\end{align}%
\end{proposition}%
\begin{proof}
    See \cref{proof:nec_eq_conds_variable_gamma_mid_steps}. \qed%
\end{proof}%
\begin{proposition}[Tightness for constant stepsize $\gamma L \in (1, {\gbari[1][]} {]} $]\label{prop:tightness_long_ct_steps}
Let $L > 0$, $\mu \in (-\infty,0)$, $N$ be a positive integer and $\gamma L \in (1,\gbari[1]{]}$. 
Then there exists $x_0$ and $f \in \mathcal{F}_{\mu,L}$ such that after applying $N$ iterations of \eqref{eq:GM_it_fixed_steps} on function $f$ with constant stepsize $\gamma L$, starting from $x_0$, it holds:%
$$
    %%%%%%%%%%%%%%%%%%%%%%%%%%%% 
    \tfrac{1}{2L} \min_{0 {} \leq {} i {} \leq {} N} \{ \|\nabla f(x_i)\|^2\}  {}={}
        \frac{f(x_0)-f_*}{
        1 + 
        \gamma L \frac{(2-\gamma L)(2-\gamma \mu)}{2-\gamma L - \gamma \mu} N
        }. %\label{eq:rate_long_ct_steps_f_star}
$$   
\end{proposition}%
\begin{proof}
We construct a two-dimensional worst-case function $f: \mathbb{R}^2 \rightarrow \mathbb{R}$. 
From \cref{prop:nec_eq_conds_variable_gamma}, the rotation angle between the gradients $\thetagLgmu$ is constant, hence they alternate between odd and even iterations such that:
$$
    \|g_i\|^2 {}={} \langle g_{i}, g_{i+(2\mathbb{N})} \rangle {}={} U^2 
        {}\,{} \text{ and } {}\,{}
    \langle g_{i}, g_{i+(2\mathbb{N}+1)} \rangle {}={} c U^2
    \qquad \forall i = 0,\dots,N-1.
$$
Consider the set of triplets $\mathcal{T}_{\mathcal{I}}{}\coloneqq{}\{\left(x_i, g_i, f_i\right)\}_{i\in\{0,\dots,N\}}$ defined as follows: (i) since $g_0$ can be set (almost) arbitrary, we take $g_{i} {}={} U \left[\sqrt{\tfrac{1{+}c}{2}};\, ({-1})^{i} \sqrt{\tfrac{1{-}c}{2}}\right]^\top$; %
(ii) by setting $x_N = [0,0]^\top$, the iterations are given by $x_i = \gamma \sum_{j=i}^{N{-}1} g_j$; (iii) by setting $f_N=0$ the function values simplify to $f_i {}={} f_0 (1 - \tfrac{i}{N})$.
In \cref{appendix:proof:tightness_long_steps_ct} we demonstrate that $\mathcal{T}_{\mathcal{I}}$ is $\mathcal{F}_{\mu,L}$-interpolable.%
\qed\end{proof}%
\cref{fig:func_example_long_steps_2D} shows a numerical example of the triplets from \cref{prop:tightness_long_ct_steps}'s proof. Excepting initial point $x_0$, all subsequent iterations are inflection points where the curvature transitions from $\mu$ to $L$. This pattern is reversed between odd and even iterations and the decrease mainly occurs to component $x^1$. %
\begin{remark}
The numerical example from \cref{fig:func_example_long_steps_2D} is obtained from the lower interpolation function derived in \cite[Theorem 3.14]{PhD_AT_2017}, originally defined for smooth strongly convex functions. This interpolation function can be extended to arbitrary lower curvatures $\mu$ using the minimal curvature subtraction argument employed in the proof of \cref{thm:interp_hypo} (refer to Remark 3.13 in \cite{PhD_AT_2017}). Given triplets $\mathcal{T}_{\mathcal{I}}{}\coloneqq{}\{(x_i, g_i, f_i)\}_{i\in \mathcal{I}}$, this function is the solution of the quadratically constrained optimization problem (QCQP):%
 $$
    f(x) {}{}\coloneqq{}{} \min_{g \in \mathbb{R}^2} \max_{0 \leq i \leq N}     
    \big\{
    f_i + \langle g_i \,,\, x-x_i \rangle + \tfrac{\mu}{2} \|x-x_i\|^2 +
    \tfrac{1}{2(L-\mu)} \|g-g_i-\mu(x-x_i)\|^2
    \big\}.
$$
\end{remark}%
When using non-constant stepsizes $\gamma_i L \in (1, \gbari[1]]$, based on the triplets construction from \cref{prop:nec_eq_conds_variable_gamma} we state in \cref{conjecture:variable_mid_stepsizes} the tightness of the convergence rate from \cref{thm:wc_GM_hypo}. The conjecture relies on numerically verification that all interpolation conditions hold (see \cref{footnote:github_repo_interp_triplets}).%
\begin{conjecture}[Lower bound for $\gamma_i L \in ({1,} {\gbari[1][]}{]}$] \label{conjecture:variable_mid_stepsizes}
There exists a two-dimensional function $f \in \mathcal{F}_{\mu,L}$ interpolating the triplets defined in \cref{prop:nec_eq_conds_variable_gamma}, satisfying%
\[
    \tfrac{1}{2L} \min_{0 \leq i \leq N} \{\|\nabla f(x_i)\|^2\}
        {}={} 
    \frac{f(x_0) - f(x_N)}{\sum_{i=0}^{N-1} \frac{\gamma_i L (2-\gamma_i L)(2-\gamma_i \mu)}{2-\gamma_i L - \gamma_i \mu}}.%
\]%
\end{conjecture}%
\subsection{Tightness analysis for weakly convex functions and stepsizes above \texorpdfstring{$\gbari[1][]$}{first threshold}%
}\label{subsec:tightness_hypo_above_thr}%
We firstly show the tightness for the stepsize interval $[{\gbari[N-1]}, {\gbari[N]})$ in \cref{prop:2D_pure_quad}, by constructing a two-dimensional quadratic function. Then \cref{conj:3D_triplets} demonstrates the exactness of our bounds for any stepsize in the range $(1, \gbari[N])$, by providing in its proof a three-dimensional example. We also provide alternative, equivalent expressions of denominator $\PN$ from \eqref{eq:P_den_full_rate_h_ct}. %
Using $\Nbar$ (recall \cref{def:N_bar}) emphasizes that the analysis of worst-case functions relies on fixing the stepsizes and varying the number of iterations; $\Nbar$ plays the role of index $k$ delimiting intervals in the proofs of \cref{thm:wc_GM_hypo_ct_all} and \cref{thm:strg_cvx_rate} from \cref{subsec:proofs:nonconvex_rates} and \cref{subsec:proof_thm_ct_strg_cvx}, respectively. %
%%%%%%%
\begin{proposition}[Tightness for $\gamma L \in {[}{\gbari[N-1][]},{\gbari[N][]})$]\label{prop:2D_pure_quad}
\oldtext{\color{gray}{\sout{Let $f \in \mathcal{F}_{\mu,L}$, with $\mu \in (-\infty, 0)$, and consider $N$ iterations of \eqref{eq:GM_it_fixed_steps} with constant stepsize $\gamma L \in [\gbari[N-1], \gbari[N])$, starting from $x_0$. Then the corresponding convergence rate from \mbox{\cref{thm:wc_GM_hypo_ct_all}} is tight:}}}%
\shownewtext{%
Let $L {>} 0$, $\mu \in (-\infty,0)$, $N$ be a positive integer and $\gamma L \in [\gbari[N-1], \gbari[N])$, such that $\Nbar {=} N{-1}$. Then there exists $x_0$ and $f \in \mathcal{F}_{\mu,L}$ such that after applying $N$ iterations of \eqref{eq:GM_it_fixed_steps} on function $f$ with constant stepsize $\gamma L$, starting from $x_0$, the corresponding performance bound from \cref{thm:wc_GM_hypo_ct_all} holds exactly:
}
\begin{equation*}
    \tfrac{1}{2L}\min_{0 {} \leq {} i {} \leq {} N} \{ \|\nabla f(x_i)\|^2 \}
        {}={}        
    \frac{f(x_0) - f_*} { 1 {}+{} \gamma L        
        \frac{ \frac{ (1-\gamma L)^2 \, \Ei[N][1-\gamma L] } {1-(1-\gamma L)^2}  - 
               \frac{ (1-\gamma \mu)^2 \, \Ei[N][1-\gamma \mu] } {1-(1-\gamma \mu)^2} 
        }{ \frac{1}{1-(1-\gamma L)^2} - \frac{1}{1-(1-\gamma \mu)^2}}
        },
\end{equation*}
with $\Ei$ given in \cref{def:def_TN}.
\end{proposition}%
\begin{proof}%
Let $\Delta > 0$ and $U = \sqrt{\frac{2L \, \Delta}{ \gamma L \PN }}$, with $\PN$ defined in \eqref{eq:P_den_full_rate_h_ct}. 
The following performance bound (without $f_*$) 
\begin{equation*}%\label{eq:rate_fN_quad}
    \tfrac{1}{2L}\min_{0 {} \leq {} i {} \leq {} N} \{ \|\nabla f(x_i)\|^2 \}  
        {}={}        
    \tfrac{f(x_0) - f(x_N)} { \gamma L } 
    \tfrac{ \frac{1}{1-(1-\gamma L)^2} - \frac{1}{1-(1-\gamma \mu)^2} } {
        \frac{ (1-\gamma L)^2 \, \Ei[N][1-\gamma L] } {1-(1-\gamma L)^2} - 
        \frac{ (1-\gamma \mu)^2 \, \Ei[N][1-\gamma \mu] } {1-(1-\gamma \mu)^2}
        }
\end{equation*}
is achieved by the $\mathcal{F}_{\mu,L}$ quadratic function%
\[%
    \begin{aligned}
        f(x) 
            {}={}& 
        \frac{1}{2} 
            (x-x_N)^\top 
            \begin{bmatrix}
                L & 0 \\
                0 & \mu
            \end{bmatrix}
            (x-x_N) + 
        g_{N}^{\top} (x-x_N) + f_N,
    \end{aligned}%
\]%
with arbitrarily fixed $x_{N}=[0,0]^\top$, $f_{N}=0$, $f(x_0) = \Delta$ and%
$$
g_{N} 
    {}{}\coloneqq{}{} 
% \min_{0 \leq i \leq N}  \left\{\|g_i\|\right\}
U \sqrt{\frac{[(1-\gamma \mu)^2-1][1-(1-\gamma L)^2]}{(1-\gamma \mu)^2-(1-\gamma L)^2}}
    \left[
    \begin{array}{c}
       \frac{1-\gamma L}{\sqrt{1-(1-\gamma L)^2}}  \\ 
       \frac{1-\gamma \mu}{\sqrt{(1-\gamma \mu)^2-1}}
    \end{array}
    \right]. {}\qquad\qquad{} \qed{}
$$%
\end{proof}%
\begin{remark}
The gradients of the worst-case function from \cref{prop:2D_pure_quad}, evaluated at iterations $x_i$, decrease in each component by the contraction factors $(1-\gamma L)$ and $(1-\gamma \mu)$, respectively:%
    \[
    \nabla f(x_{i+1}) 
        {}={}
    \begin{bmatrix}
        1-\gamma L & & 0 \\
        0 & & 1 - \gamma \mu
    \end{bmatrix}
    \nabla f(x_{i}).
    \]%
\end{remark}%

Proving the tightness of \cref{thm:wc_GM_hypo_ct_all} for stepsizes $\gamma L \in \big[\gbari[1], \gbari[N-1]\big)$ presents a significant challenge due to the nonlinear term in \eqref{eq:P_den_full_rate_h_ct}, which connects the quadratic function characterizing the first $\Nbar+1$ iterations (see \cref{prop:2D_pure_quad}) and the 2D function with alternating gradients directions from  \cref{prop:tightness_long_ct_steps}. This later function dominates the sublinear rate through the leading coefficient $\gamma L \, \pigLgmu$. 
The three-dimensional worst-case construction from \cref{conj:3D_triplets}'s proof also provides an alternative worst-case construction for this regime with stepsizes $\gamma L \in (1, \gbari[1][]{]}$. 
\begin{proposition}[Tightness for $\gamma L \in (1, {\gbari[N][]} {)}$]\label{conj:3D_triplets}
Let $L > 0$, $\mu \in (-\infty,0)$, $N$ be a positive integer and $\gamma L \in (1,\gbari[N]{)}$. 
There exists $x_0$ and $f \in \mathcal{F}_{\mu,L}$ such that after applying $N$ iterations of \eqref{eq:GM_it_fixed_steps} on function $f$ with constant stepsize $\gamma L$, starting from $x_0$, the performance bound from \cref{thm:wc_GM_hypo_ct_all} is attained:%
$$ 
    \tfrac{1}{2L} \min_{0 {} \leq {} i {} \leq {} N} \{ \|\nabla f(x_i)\|^2\}  
        {}={}        
    \frac{f(x_0) - f_*}{
    1 + 
    \frac{(2-\gamma L)(2-\gamma \mu)}{2-\gamma L - \gamma \mu} (N - \bar{N}) {}+{} 
        \frac{ \frac{E_{\bar{N}}(1-\gamma L)}{1-(1-\gamma L)^2} - \frac{E_{\bar{N}}(1-\gamma \mu)}{1-(1-\gamma \mu)^2} }{ \frac{1}{1-(1-\gamma L)^2} - \frac{1}{1-(1-\gamma \mu)^2} }
    },
$$%
where $\Ei$ given in \cref{def:def_TN}, together with the complementary rate:
$$
    \tfrac{1}{2L} \min_{0 {} \leq {} i {} \leq {} N} \{ \|\nabla f(x_i)\|^2\}  
        {}={}        
    \frac{f(x_0) - f(x_N)}{
    \frac{(2-\gamma L)(2-\gamma \mu)}{2-\gamma L - \gamma \mu} (N - \bar{N}) {}+{} 
        \frac{ \frac{E_{\bar{N}}(1-\gamma L)}{1-(1-\gamma L)^2} - \frac{E_{\bar{N}}(1-\gamma \mu)}{1-(1-\gamma \mu)^2} }{ \frac{1}{1-(1-\gamma L)^2} - \frac{1}{1-(1-\gamma \mu)^2} }
    }.%
$$%
\end{proposition}%
\begin{proof}
Let $\Delta > 0$ and $U = \sqrt{\frac{2L \, \Delta}{ \gamma L \PN }}$, with $\PN$ defined in \eqref{eq:P_den_full_rate_h_ct}. The value of $U$ denotes the minimum gradient norm over $N$ iterations, i.e., the rate, where $\Delta$ accounts for the initial function value gap $f(x_0) - f(x_N)$. 
    The definition of $\bar{N}=\Nbar$ implies $\gamma L \in [\gbari[\bar{N}-1], \gbari[\bar{N}])$, with $1 \leq \bar{N} \leq N-1$. 
    We define the set of triplets $\mathcal{T}_{\mathcal{I}}{}\coloneqq{}\{\left(x_i, g_i, f_i\right)\}_{i=\{0,1,\dots,N\}}$ as follows:
    \begin{align*}
    &\begin{aligned}
    g_i 
        {}={}& 
    U \sqrt{\tfrac{[(1-\gamma \mu)^2-1][1-(1-\gamma L)^2]}{(1-\gamma \mu)^2-(1-\gamma L)^2}}
    \left[
    \begin{array}{c}
       \tfrac{(1-\gamma L)^{i-\bar{N}}}{\sqrt{1-(1-\gamma L)^2}}  \\ 
       \tfrac{(1-\gamma \mu)^{i-\bar{N}}}{\sqrt{(1-\gamma \mu)^2-1}}  \\
       0
    \end{array}
    \right]
    {}\quad{} \forall i=0,1,\dots,\bar{N}+1; \\
    %%%%%%%%%%%%%%%%%%%%%%%%%%%%%%%%%%%%%%%
    g_{i+1}
        {}={}&
    (1-\gamma L) g_{i} + \gamma(L-\mu) \cos(\theta_{i+1})
    \Big[ \begin{array}{cc}  0  &  \\ &  R(\theta_{i+1}) \end{array} \Big] g_{i} {}\quad{} \forall i=\bar{N}, \dots, N; \\
    %%%%%%%%%%%%%%%%%%%%%%%%%%%%%%%%%%%%%%%
    q_{i} {}={}& (-1)^{i-(\bar{N}+1)} \sqrt{[(1-\gamma \mu)^2-1] \tfrac{1-(1-\gamma L)^{2(i-(\bar{N}+1))}}{1-(1-\gamma L)^2}} {}\quad{} \forall i=\bar{N}+1, \dots, N;  \\
    %%%%%%%%%%%%%%%%%%%%%%%%%%%%%%%%%%%%%%%
    \theta_{i} {}={}& \arctan (q_{i}) {}\quad{} \forall i=\bar{N}+1, \dots, N;
    \end{aligned} \\
    %%%%%%%%%%%%%%%%%%%%%%%%%%%%%%    
    & \begin{aligned}
        f_i {}={}
        \Delta
    \tfrac{ N-i + \frac{\gamma L \gamma \mu}{\gamma L - \gamma \mu} \sum_{j=1}^{\bar{N}-i}\TigLmu[j] }
    {N + \frac{\gamma L \gamma \mu}{\gamma L - \gamma \mu} \sum_{j=1}^{\bar{N}}\TigLmu[j]}  {}\quad{} i=0,1,\dots,N,
    \end{aligned}
\end{align*} 
$f_N=0$, $x_N=0$ and $x_{i+1} = x_i - \gamma g_i$, $i=0,\dots,N-1$ and $R(\theta) \coloneqq \left[\begin{array}{cr}
   \cos(\theta)  & {-}\sin(\theta) \\
   \sin(\theta)  &  \cos(\theta)
\end{array}\right]$ being the two-dimensional rotation matrix of angle $\theta$. 
For $\mathcal{T}_{\mathcal{I}}$, we show in \cref{sec:proof:3D_wc_nonconvex} that: 
(i) it attains the performance bound from \cref{thm:wc_GM_hypo_ct_all} and 
(ii) it is interpolable by a $\mathcal{F}_{\mu,L}$-function with three variables.\qed
\end{proof}
%%%%%%
\begin{remark}%
If the stepsize is not the threshold $\gbari[\bar{N}]$, then the first $\bar{N}+1$ gradients have uniquely defined norm and inner products between them. This follows the equality conditions from \cref{lemma:hypo_GN4SD}. Restricting them to be two-dimensional, they are unique up to applying orthonormal transformations. The first $\bar{N}+1$ iterations are exactly characterized by the worst-case function from \cref{prop:2D_pure_quad}, where a third component is set to zero. %
For the subsequent iterations, the first component remains linearly decreasing in $({1-}\gamma L)$, while the others are obtained by applying a certain rotation which, given the decrease in the first component, preserves: (i) the gradient norm to $U^2$ and (ii) the inner product between consecutive iterations to $c U^2$ such that:%
\[
\begin{aligned}
    \text{(i):} {}\,\,{}& 
        \|g_{i+1} - (1-\gamma L) g_i\|^2 
            {}={} 
        U^2 \big[1 + (1-\gamma L)^2 - 2(1-\gamma L) c \big]; \\
    \text{(ii):} {}\,\,{}& 
        \langle g_{i+1} - (1-\gamma L) g_i \,,\, g_i \rangle 
            {}={} 
        U^2 \big[c - (1-\gamma L) \big],
\end{aligned}
\]%
with $c = \tfrac{1+(1-\gamma L)(1-\gamma \mu)}{2- \gamma L - \gamma \mu}$ from \cref{prop:tightness_long_ct_steps}. The matrix $R(\theta_{i+1})$ rotates the second and third components of $g_i$ by angle $\theta_i$, whose absolute value monotonically increasing towards 
$\lim_{i \rightarrow \infty} |\theta_i| = \arctan (\sqrt{ \frac{{-}\mu}{L} \frac{2{-}\gamma \mu}{2{-}\gamma L}}).$
These components are then properly scaled to preserve $\|g_{i+1}\|=U$. Moreover, as $i$ approaches $\infty$, the first component vanishes and the gradients start alternating, resembling the same pattern observed for the two-dimensional worst-case described in \cref{prop:tightness_long_ct_steps}'s proof.%
\end{remark}%
\section{Conclusion}\label{sec:conclusion}%
    %!TEX root = ../../ms.tex
%
%
We have presented a comprehensive worst-case analysis of gradient descent applied to smooth weakly convex, convex and strongly convex functions across the entire range of constant stepsizes $\gamma L \in (0,2)$. %

For convex and strongly convex functions we establish and prove tight performance bounds over ranges previously only covered by conjectures.

In the weakly convex scenario, determining convergence rates is more challenging outside extreme cases already addressed in the literature: $\mu \in \{ -\infty, 0\}$. The state-of-the-art analysis of smooth not-necessarily convex functions is extended from the particular case $\mu=-L$ and stepsizes $\gamma L \in (0, \sqrt{3}]$ to encompass full domains. 

The proofs are based on incorporating not only consecutive interpolation inequalities but also those connecting iterations at a distance of two, the key idea being detailed in \cref{remark:main_idea_cancellations}. A key advantage of these proofs type is the obtaining of new descent lemma type inequalities, independent of the number of iterations $N$, as they would be in tight convergence analyses using inequalities involving the last iterate. %
\shownewtextfinal{When deriving closed-form expressions for tight performance analysis with constant stepsizes, it is uncommon and challenging to involve nonconsecutive iterations.}%

While all rates align with the numerical solutions obtained from solving performance estimation problems (PEPs), we demonstrate analytically their tightness over the complete stepsize range $\gamma L \in (0, 2)$.%

This analysis yields insights into the optimal constant stepsizes with respect to the demonstrated worst-case scenarios and we further improve upon them by proposing dynamic stepsize schedules which have the key benefit of being independent of the number of iterations.%
%
%
%
% \begin{acknowledgements}
%     % This research was funded within the framework of the Global PhD Partnership KU Leuven - UCLouvain.
%     The authors thank the anonymous referees and the handling editor for their detailed and constructive comments, which substantially improved the presentation and strengthened the results compared to the initial version of the paper.
% \end{acknowledgements}

% \medskip
% \noindent\textbf{Funding.} This work was supported by the framework of the Global PhD Partnership KU Leuven -- UCLouvain. \\
%
% \noindent\textbf{Code availability.} The symbolic script to check the derivations and the numerical examples are available on {\href{https://github.com/teo2605/GD_tight_rates}{GitHub repository:} \url{https://github.com/teo2605/GD_tight_rates}}.%

%apparently not working  
%\begin{appendices}
\appendix{%
% \clearpage
\section{Properties of stepsize sequences (proofs)}%
\label{appendix:properties_stepsize_sequences}%
    %!TEX root = ../../main.tex
\subsection{Stepsize thresholds} \label{appendix:properties_stepsize_thresholds}
\begin{definition}
    We use the writing $\TigLmu[k][\gamma L] {}={} \sum_{i=1}^k \delta_i(\gamma L, \kappa)$, where
    $$
        \delta_i(\gamma L,\kappa) 
            {}{}\coloneqq{}{} 
        \frac{2-\kappa \gamma L}{(1-\kappa \gamma L)^{2i}} - 
        \frac{2-\gamma L}{(1-\gamma L)^{2i}}.
    $$ 
\end{definition}%
\begin{proof}[\textbf{\cref{prop:monotonicity_Tk}}] \label{appendix:proof_prop_monotonicity_Tk}
We show $\frac{d\TigLmu}{d (\gamma L)} \geq 0$ for $\gamma L \in (1, \gbari[\infty])$, which resumes to prove the positivity of
  \begin{align*}
        \frac{d \delta_i(\gamma L, \kappa)}{d(\gamma L)}
            {}{}\coloneqq{}{}
      \frac{-\kappa \left( 1 - 2i \frac{2-\kappa \gamma L}{1-\kappa \gamma L} \right) }{(1-\kappa \gamma L)^{2i}} + \frac{ 1 + 2i \frac{2-\gamma L}{\gamma L-1}}{(1-\gamma L)^{2i}}.
  \end{align*}
  \begin{enumerate}[align=right,itemindent=1em,labelsep=2pt,labelwidth=1em,leftmargin=0pt,nosep]
      \item \textbf{Case $\kappa \geq 0$:} we lower bound $\frac{d \delta_i(\gamma L, \kappa)}{d(\gamma L)}$ by neglecting the denominators:
      \begin{align*}
          \delta_i(\gamma L, \kappa) 
            {}\geq{} 
          {-}\kappa \Big( 1 - 2i \frac{2-\kappa \gamma L}{1-\kappa \gamma L} \Big) + 1 + 2i \frac{2-\gamma L}{\gamma L-1}
            {}={}
          (1-\kappa) + 2i \Big[ \frac{\kappa (2-\kappa \gamma L)}{1-\kappa \gamma L}  + \frac{2-\gamma L}{\gamma L-1} \Big] {}\geq{} 0.
      \end{align*}%
      \item \textbf{Case $\kappa < 0$:} it is sufficient to show 
    \begin{align*}
        \frac{1+2i \frac{2-\gamma L}{\gamma L -1}}{(1-\gamma L)^{2i}}
            {}\geq{} 1 {}\geq{}
        -\kappa \frac{-1+2i (1 + \frac{1}{1-\kappa \gamma L})}{(1-\kappa \gamma L)^{2i}}.
    \end{align*}
      The \lhs holds because $\gamma L \in (1,2)$ and $(1-\gamma L)^{2i}<1$. The \rhs is equivalent to:%
      \begin{align*}
          (1-\kappa \gamma L)^{2i} {}\geq{} -\kappa \Big[-1+2i \Big(1 + \frac{1}{1-\kappa \gamma L}\Big)\Big].
      \end{align*}
Using $\gamma L \in (1,2)$ and the second order Taylor approximation for $(1-\kappa)^{2i}$, we get the bounds:%
      \begin{align*}
          (1-\kappa \gamma L)^{2i} {}\geq{}& (1-\kappa)^{2i} 
             {}\geq{} 
         1 + 2i (-\kappa) + i(2i-1) \kappa^2 \\
         %%%%%%%%%%%%%%%%%%%%%%%%%%%%%%%%%%%%%%%%%
          -\kappa \Big[-1+2i \Big(1 + \frac{1}{1-\kappa \gamma L}\Big)\Big]
            {}\leq{}&
          -\kappa \Big[-1+2i \left(1 + \frac{1}{1-\kappa}\right)\Big].
      \end{align*}
Therefore, it is sufficient to show:%
$$
    1 + 2i (-\kappa) + i(2i-1) \kappa^2 {}\geq{} 
    -\kappa \Big[-1+2i \Big(1 + \frac{1}{1-\kappa}\Big)\Big],%
$$%
which equivalently rewrites as the following inequality valid for all $i \geq 1$:
$$
    \kappa^2 i(i-1) + \Big( \kappa i + \frac{1}{1-\kappa} \Big)^2 + \frac{(1-\kappa)^3-1}{(1-\kappa)^2} \geq 0. {}\qquad \qed{}%
$$%
\end{enumerate}
\end{proof}

\begin{lemma} \label{lemma:unique_delta_k_neg_thr}
Let $\kappa \in (-\infty, 1)$ and $\gamma L \in (1,\gbari[\infty])$. Then there exists an index $k$ such that $\delta_k(\gamma L, \kappa) < 0$ and $\delta_j(\gamma L, \kappa) < 0$ for all $j > k$.    
\end{lemma}
\begin{proof}
    The sequence $\delta_i$ is a difference of two monotone sequences: $\frac{2-\kappa \gamma L}{(1-\kappa \gamma L)^{2i}}$ and $\frac{2-\kappa \gamma L}{(1-\kappa \gamma L)^{2i}}$, respectively. The former is decreasing for $\kappa < 0$ (weakly convex functions) and increasing for $\kappa \geq 0$ (convex functions), whereas the later is always decreasing. In particular, $\delta_0(\gamma L, \kappa) > 0$. 
    Because $(1-\gamma \mu)^2 {}\leq{} (1-\gamma L)^2$, to the limit it holds $\lim_{k \rightarrow \infty} \delta_k(\gamma L, \kappa) = -\infty$ and there exists some index $k$ with $\delta_k(\gamma L, \kappa) < 0$. We prove by contradiction that $\delta_j(\gamma L, \kappa) < 0$ for all $j > k$; assume there exists some index $l > k$ such that $\delta_l(\gamma L, \kappa) \geq 0$. The following inequalities hold:%
\begin{align*}
    \delta_{k} < 0 \,:\,& 
        \frac{2-\gamma \mu}{(1-\gamma \mu)^{2k}} 
            {}<{} 
        \frac{2-\gamma L}{(1-\gamma L)^{2k}} \\
    %%%%%%%%%%%%%%%%%%%%%%%%%%%
    \delta_{l} > 0 \,:\,& 
        \frac{2-\gamma \mu}{(1-\gamma \mu)^{2l}} 
            {}>{} 
        \frac{2-\gamma L}{(1-\gamma L)^{2l}}.
\end{align*}
By multiplying the inequalities and performing the simplifications we get:%
\begin{align*}
    (1-\gamma \mu)^{2(k-l)} > (1-\gamma L)^{2(k-l)},
\end{align*}
which implies, due to $(1-\gamma \mu)^2 {}>{} (1-\gamma L)^2$, that $k > l$, contradiction! \qed
\end{proof}%
\begin{proof}[\textbf{\cref{prop:properties_step_size_thr}}] \label{appendix:proof_properties_step_size_thr} 
\begin{enumerate}[align=right,itemindent=1em,labelsep=2pt,labelwidth=1em,leftmargin=0pt,nosep]
\item Case $k=0$ holds by definition since $\gbari[0]=1$. For any $k \geq 1$, by using \cref{lemma:unique_delta_k_neg_thr}, condition $\TigLmu[k][\gbari[k][]] = 0$ implies $\delta_k(\gamma L, \kappa) \leq 0$ and hence $\delta_{k+1}(\gamma L, \kappa) < 0$. Furthermore, $\TigLmu[k+1][\gbari[k][]] {} \leq {} \TigLmu[k][\gbari[k][]] = 0$. Then it holds:
\begin{align*}
    \TigLmu[k+1][\gbari[k][]] {} \leq {} 
    \TigLmu[k][\gbari[k][]] = 0 
    {}={} \TigLmu[k+1][\gbari[k+1][]].
\end{align*}
The monotonicity of $\TigLmu$ in stepsize $\gamma L$ (see \cref{prop:monotonicity_Tk}) implies $\gbari[k] < \gbari[k+1]$.

\item Because stepsize thresholds sequence is monotone, it has a limit.
\begin{enumerate}[align=right,itemindent=1.5em,labelsep=2pt,labelwidth=1em,leftmargin=0pt,nosep]
    %%%%%%%%%%%%%%%%%%%%%%%%%% 2. %%%%%%%%%%%%%%%%%%%%%%%%%%%%%%%%%%%%%%%%%%%
    \item 
    \textbf{Case $\kappa > 0$.} By definition, $\gbari[k] < \frac{2}{1+\kappa} < 2$ for all integers $k\geq 1$ and condition $\TigLmu[k][\gbari[k]]=0$ implies:%
    \begin{align*}
        \left(\frac{1-\gbari}{1-\kappa \gbari}\right)^{2k}
            {}={}& 
        \kappa + (1-\kappa) (1-\gbari)^{2k}
            {}\Leftrightarrow{} \\
            %%%%%%%%%%%%%%%%%%%%%%%%%%%%
        \left|\frac{1-\gbari}{1-\kappa \gbari}\right|
            {}={}&
        \exp \left[ \frac{\ln \left(\kappa + (1-\kappa) (1-\gbari)^{2k}\right)}{2k} \right].
    \end{align*}%
    To the limit $k \rightarrow \infty$, the \rhs equals $1$ because $(1-\gbari)^{2k}\rightarrow 0$, hence $|1-\gbari|=|1-\kappa \gbari|$, which implies the solution $\gbari[\infty]=\frac{2}{1+\kappa}$. \\
    %%%%%%%%%%%%%%%%%%%%%%%%%%
    \item \textbf{Case $\kappa < 0$.} By definition, $\gbari < 2$ for all integers $k\geq 1$ and condition $\TigLmu[k][\gbari[k]]=0$ implies:
    \begin{align*}
    (1-\gbari)^{-2k} 
        {}={}&
    1 + \frac{1 - (1-\kappa \gbari)^{-2k}}{-\kappa}
        {}\Leftrightarrow{} \\
    % %%%%%%%%%%%%%%%%%%%%%%%%%%%%%%%%%%%%%%%%%%%%%%%%%%%%%%%%%
    \left|1-\gbari\right| 
        {}={}&
    \exp \Bigg[ \frac{ -\ln \Big(1 + \frac{1 - (1-\kappa \gbari)^{-2k}}{-\kappa}\Big)}{2k} \Bigg].
\end{align*}
 To the limit $k \rightarrow \infty$, the \rhs equals $1$ because $(1-\kappa \gbari)^{-2k} \rightarrow 0$ due to $\kappa < 0$. Hence, $\left|1-\gbari\right|=1$, with the solution $\gbari[\infty]=2$.%
\item \textbf{Case $\kappa = 0$.} By definition, $\gbari[k][(0)] < 2$ and, similarly to the case $\kappa < 0$:%
 \begin{align*}
     \left|1-\gbari\right|
        {}={}&
     \left[ 1 + 2 k \gbari \right]^{\frac{-1}{2k}}
        {}={}
    \exp \Big[ \frac{-\ln \big( 1 + 2 k \gbari \big) }{2k} \Big].
 \end{align*}
  To the limit $k \rightarrow \infty$, the \rhs equals $1$, hence $\left|1-\gbari\right|=1$, with the solution $\gbari[\infty]=2$.
\end{enumerate}%
\item Condition $\gamma L \in [\gbari[k], \gbari[k+1])$ implies $\TigLmu[k][\gbari[k]] \geq 0$ and $\TigLmu[k+1][\gbari[k]] < 0$. Because $\delta_i(\gamma L, \kappa) < 0$ for all $i \geq k$, we have (i) $\TigLmu[i][\gbari[k][]] {}>{} 0$ for $i=1,\dots,k-1$ (since $k$ is uniquely defined) and (ii) $\TigLmu[i][\gbari[k][]] {}<{} 0$ for $i\geq k+1$. {}\qed{}%
\end{enumerate}
\end{proof}%%
\subsection{Proof of optimal constant stepsize for weakly convex functions}\label{appendix:optimal_constant_stepsize_nonconvex}%
\begin{proof}[\textbf{\cref{prop:gamma_star}}]
    Minimizing the upper bound from \cref{thm:wc_GM_hypo_ct_all} is equivalent to maximizing the denominator $\PN$ from \eqref{eq:P_den_full_rate_h_ct}.
    For a fixed curvature ratio $\kappa = \frac{\mu}{L}$, this reduces to maximizing the quantity $r(l) {}{}\coloneqq{}{} l \, p(l, \kappa l)$, where $l \in (0,2)$ is the normalized stepsize $\gamma L$ and $p(l,u)$ is defined in \eqref{eq:pi_hi_rate_GM}:%
    \begin{align*}
    \begin{aligned}%\label{eq:opt_problem_opt_step}    
        r(l)
            {}={}
    l
    \Big[
        2 - \tfrac{-\kappa l^2}{1-\kappa l - |1-l|}
    \Big]
            {}={}
         \left\{
        \arraycolsep=6pt
    \def\arraystretch{2}
        \begin{array}{ll}
            l \bigl( 2 - \tfrac{\kappa l^2}{l - \kappa l}\bigr) &  l \in (0,1]; \\
             l \tfrac{(2-l)(2-\kappa l)}{2-l- \kappa l} & l \in [1, \gbari[1]].
        \end{array}
    \right.
    \end{aligned}
\end{align*}
Its first branch is concave in $l$ and reaches its maximum for $l = 1$. Since this value is also feasible for the second branch, it means that the maximum belongs to the later expression and restrict the analysis to the interval $[1,2)$. One can check that $r(l)$ is strictly concave for $l \in [1,2)$. Moreover, $r'(1) =\frac{2}{(1-\kappa)^2}>0$ and $r'(2)=-2(1-\tfrac{1}{\kappa})< 0$, where%
    \begin{align}\label{eq:derivative_of_p_den}
    r'(l) {}={}
    \frac{2}{(2-l - \kappa l)^2} 
    \left[%
    -\kappa\big(1+\kappa\big) l^3 +
    \big[3\kappa + \left(1+\kappa\right)^2 \big] l^2
    -4\big(1+\kappa\big) l +
    4 \right].
    \end{align}%
Hence, there exists a unique maximizer $l^* \in [1,2)$, further on denoted by $\gstar[(\kappa)]$. The square bracket from \eqref{eq:derivative_of_p_den} is exactly the expression from equation \eqref{eq:equation_solve_for_h_star_one_step}. 

Moreover, the optimal stepsize $\gstar[(\kappa)]$ strictly increases with $\kappa$, from $\gstar[(-\infty)]=1$ to $\gstar[(0)] \nearrow 2$ and crosses the threshold $\gbari[1]$ at some specific $\bar{\kappa}$ (see \cref{fig:Full_picture_all_regimes}), which results by taking $\gstar[(\kappa)] {}={} \gbari[1]$ in \eqref{eq:equation_solve_for_h_star_one_step}. For longer stepsizes, some other transient and exponential terms in $\gamma$ are involved (see expression of denominator \eqref{eq:P_den_full_rate_h_ct}), therefore the solution $\gstar[(\kappa)] > \gbari[1]$ is only asymptotically optimal ($N \rightarrow \infty$).%
\qed%
\end{proof}%%
%
% \clearpage
\subsection{Scheduled stepsizes} \label{appendix:properties_variable_stepsizes}
\begin{proof}[\textbf{\cref{prop:properties_variable_stepsizes_sequence}}]\label{proof:prop:properties_variable_stepsizes_sequence}
\begin{enumerate}[label=(\roman*)]
    \item One can check that $s_0 = \gbari[1] > 1$, where $\gbari[1]$ is given explicitly in \cref{thm:wc_GM_hypo}. Further on, $s_k, s_{+} \in (1, \frac{2}{1+[\kappa]_{+}})$, for $k \geq 0$. By definition of $s_+$:
    \begin{align*}
    0 {}={}
        \frac{s_{k}}{2-s_{k}(1+\kappa)}
        - \frac{s_{+} }{2-s_{+}(1+\kappa)}
        + \frac{s_{+} (2-s_{+})(2-\kappa s_{+})}{2-s_{+}(1+\kappa)}.            
    \end{align*}
    Since $\frac{s_{+}(2-s_{+})(2-\kappa s_{+}) }{2-s_{+}(1+\kappa)} {}\geq{} 0$, it implies:
    \begin{align*}
        \frac{s_{+} }{2-s_{+}(1+\kappa)}
            {}\leq{}
        \frac{s_{k}}{2-s_{k}(1+\kappa)}
        {}\iff{}
            s_{k} \leq s_{+}.        
    \end{align*}%
    After simplifications in \eqref{eq:def_sequence_stepsizes}, $s_{+}$ are the roots of the function:
    \begin{align*}
        h(s_{+}) 
            {}={}
        -\kappa[2-s_{k}(1+\kappa)] s_{+}^3  +
        2(1+\kappa) [2-s_{k}(1+\kappa)] s_{+}^2  +
        2[-3 + s_{k}(1+\kappa)] s_{+}
        -2 s_{k}.
    \end{align*}%
The demonstration of the existence of a unique root $s_{k}^+$ is divided into scenarios of weak and strong convexity.
    \begin{enumerate}[align=right,itemindent=1em,labelsep=2pt,labelwidth=1em,leftmargin=0pt,nosep]
    \item \textbf{Case $\kappa \in (-\infty, 0]$.} Then $s_{k} \in [1, 2)$, $s_0 = \gbari[1][(0)]=\frac{3}{2}$ and the following holds:%
    \begin{enumerate}[align=right,itemindent=1em,labelsep=2pt,labelwidth=1em,leftmargin=0pt,nosep]
    \item $h(1) = -(1-\kappa)(2-\kappa s_{k}) < 0$;
    \item $h(3/2) = \frac{-s_{k}(1-3\kappa)}{2} + \frac{9}{8} [2-s_{k}(1+\kappa)] < 0$;
    \item $h(2) = 2(2-s_{k})$;
     \item $\frac{d h}{d s_{+}} (1) = (1+\kappa)(2-\kappa s_{k}) \geq 0$ for $\kappa \in [-1, 1)$;
    \item $\frac{d h}{d s_{+}} (\frac{3}{2}) = 2 + [2-s_{k}(1+\kappa)] ( 2 - \frac{3}{4}\kappa ) > 0$ for $\kappa \in (-\infty, 0]$;
    \item $\frac{d h}{d s_{+}} (2) = 2 + 4(1-\kappa)[2-s_{k}(1+\kappa)] > 0$ for $\kappa \in (-\infty, 0]$;
    \item $\frac{d^2 h(s_{+})}{d s_{+}^2}
            {}={}
        2[2-s_{k}(1+\kappa)] [2-\kappa -3\kappa (s_{+} - 1)] > 0$.
    \end{enumerate} %
    Therefore, $h$ is strictly increasing and possesses a unique root $s_{+} \in [\frac{3}{2},2)$.%
    \item \textbf{Strongly convex case.} We restrict to $\kappa \in [0,1)$ and $s_{k} \in [1, \frac{2}{1+\kappa})$. One can check:
\begin{enumerate}[align=right,itemindent=1em,labelsep=2pt,labelwidth=1em,leftmargin=0pt,nosep]    
    \item $h(1) = -(1-\kappa)(2-\kappa s_{k}) < 0$;
    \item $h(\frac{2}{1+\kappa}) = \frac{2(1-\kappa)^2[2-s_{k}(1+\kappa)]}{(1+\kappa)^3} > 0$;
    \item $\frac{d h}{d s_{+}} (1) = (1+\kappa)(2-\kappa s_{k}) > 0$;
    \item $\frac{d h}{d s_{+}} (\frac{2}{1+\kappa}) = 2(1+\kappa^2)[2-s_{k}(1+\kappa)] + s_{k} \kappa (1+\kappa) + (1-\kappa)^2 > 0$;
    \item $\frac{d^2 h(s_{+})}{d s_{+}^2} 
        {}={}
    2[2-s_{k}(1+\kappa)]
    \left[
    \frac{3\kappa [2-s_{k}^+(1+\kappa)] + 2 [ \kappa + (1-\kappa)^2 ] }{1+\kappa}
    \right] > 0$.
    \end{enumerate}
    Therefore, $h$ is strictly increasing and possesses a unique root $s_{+} \in [1, \frac{2}{1+\kappa})$.%
    \end{enumerate}%
%%%%%%%%%%%%%%%%%%%%%%%%%%%%%%%%%%%%%%%%%%%%%%%%%%%%%%%%%%%%%%%%%%%%%%%%%%%%%%%%%%%%%%%%%%%%%%%%%%%%    
\item Since $s_{k+1} = s_{+}$ is uniquely defined and $s_{k} \leq s_{+}$, the sequence is monotonically increasing.

\item Due to its monotonicity, the sequence $\left\{s_k(\kappa)\right\}_{k=-1}^{\infty}$ has a limit, denoted by $s_{\infty}$. We prove by contradiction that $s_{\infty} = \frac{2}{1+[\kappa]_{+}}$. Assume the contrary, i.e., $s_{\infty} \neq \frac{2}{1+[\kappa]_{+}}$, which is the upper bound of the interval of roots, so that $s_{\infty} \in [1, \frac{2}{1+[\kappa]_{+}})$. By replacing it in the definition of \eqref{eq:def_sequence_stepsizes}, it holds 
    $
        \frac{s_{\infty} (2-s_{\infty})(2-\kappa s_{\infty}) }{2-s_{\infty}(1+\kappa)}
            {}={}
        0,
    $
    which implies one of the solutions $s_{\infty} \in \{0,2,\frac{2}{\kappa}\}$, which are all outside of the interval $[1, \frac{2}{1+[\kappa]_{+}})$, hence the contradiction!

\item When $\kappa = 0$, we define $c_j \coloneqq 2 - s_j(0)$, $\forall j = \{-1,0,\dots\}$. 
From the recurrence it results $2-c_{j+1} - \frac{1}{c_{j+1}} + \frac{1}{c_j} = 0$. 
Since $\{s_j\}$ is monotonically increasing, with $s_0 = \frac{3}{2}$, we have that $c_0 = \frac{1}{2}$, $\{c_j\}$ is a decreasing sequence with $c_j \in (0, \frac{1}{2}]$, $\forall j \geq 0$, and $c_{\infty} = 0$. 
Thus, $\frac{1}{c_{j+1}} {}\geq{} \frac{1}{c_j} + \frac{3}{2}$, $\forall j \geq 0$. 
Telescoping for $j=0, \dots, k-2$ we get $\frac{1}{c_{k-1}} {}\geq{} \frac{1}{c_0} + \frac{3}{2} (k-1)$, hence $c_{k-1} {}\leq{} \frac{2}{3k+1}$ and $s_{k-1} {}\geq{} 2 - \frac{1}{\frac{3}{2} k}$.%

\item When $\kappa \in (0,1)$, we define $c_j \coloneqq \frac{2}{1+\kappa} - s_j(\kappa)$, $\forall j = \{-1,0,\dots\}$. 
Since $\{s_j\}$ is monotonically increasing, then $\{c_j\}$ is a monotonically decreasing sequence with $c_\infty = 0$.
Using $s_0 = \gbari[1] = \frac{3}{1+\kappa+\sqrt{1-\kappa+\kappa^2}}$, one can check that $c_0 \leq \min\{\frac{2}{1+\kappa} (\frac{1-\kappa}{1+\kappa})^2\,,\, \frac{1}{2} \}$, thus $c_j \in (0, \frac{1}{2}]$, $\forall j \geq 0$. 
We prove $\frac{c_{j+1}}{c_j} \leq (\frac{1-\kappa}{1+\kappa})^2$, $\forall j \geq 0$. By replacing in the recurrence of $s_j(\kappa)$ we get %
$$
        \frac{ (\frac{2}{1+\kappa}-c_{j+1}) (\frac{2\kappa}{1+\kappa}+c_{j+1})(\frac{2}{1+\kappa}+ \kappa c_{j+1}) }{c_{j+1}(1+\kappa)}
            -
        \frac{2}{c_{j+1}(1+\kappa)^2}
            +
        \frac{2}{c_j(1+\kappa)^2} = 0.%
$$%
After algebraic manipulations, this is equivalent with
$$
\begin{aligned}
\frac{c_{j+1}}{c_j} {}={}& 1 
        - \frac{1+\kappa}{2} (\frac{2}{1+\kappa}-c_{j+1}) (\frac{2\kappa}{1+\kappa}+c_{j+1})(\frac{2}{1+\kappa}+ \kappa c_{j+1}) \\
%     {}={}&
% (\frac{1-\kappa}{1+\kappa})^2 + 
% \frac{4\kappa}{(1+\kappa)^2} 
%         - \frac{1+\kappa}{2} (\frac{2}{1+\kappa}-c_{k+1}) (\frac{2\kappa}{1+\kappa}+c_{k+1})(\frac{2}{1+\kappa}+ \kappa c_{k+1}) \\
%%%%%%%%%%%%%%%%%%%%%%%%%%%%%%%%%%%%%%%%%%%%%
    % {}={}&
% (\frac{1-\kappa}{1+\kappa})^2 + 
% \frac{4\kappa}{(1+\kappa)^2}
% \Big[
%        1 - \frac{(1+\kappa)^3}{8\kappa} (\frac{2}{1+\kappa}-c_{k+1}) (\frac{2\kappa}{1+\kappa}+c_{k+1})(\frac{2}{1+\kappa}+ \kappa c_{k+1})
% \Big] \\
% %%%%%%%%%%%%%%%%%%%%%%%%%%%%%%%%%%%%%%%%%%%%%
%     {}={}&
% (\frac{1-\kappa}{1+\kappa})^2 + 
% \frac{4\kappa}{(1+\kappa)^2}
% \Big[
%        1 - (1-\frac{1+\kappa}{2}c_{k+1}) (\kappa + \frac{1+\kappa}{2} c_{k+1})( \frac{1}{\kappa} + \frac{1+\kappa}{2} c_{k+1})
% \Big] \\
% %%%%%%%%%%%%%%%%%%%%%%%%%%%%%%%%%%%%%%%%%%%%%
%     {}={}&
% (\frac{1-\kappa}{1+\kappa})^2 + 
% \frac{4\kappa}{(1+\kappa)^2}
% \frac{(1+\kappa)c_{k+1}}{8 \kappa}
% \Big[
%        \kappa c_{k+1}^2 (1+\kappa)^2 + 2 c_{k+1} (1+\kappa^3) - 4 \frac{1+\kappa^3}{1+\kappa}
% \Big] \\
% %%%%%%%%%%%%%%%%%%%%%%%%%%%%%%%%%%%%%%%%%%%%%
%     {}={}&
% (\frac{1-\kappa}{1+\kappa})^2 + 
% \frac{c_{k+1}}{2(1+\kappa)}
% \Big[
%        \kappa c_{k+1}^2 (1+\kappa)^2 + 2 c_{k+1} (1+\kappa^3) - 4 \frac{1+\kappa^3}{1+\kappa}
% \Big] \\
%%%%%%%%%%%%%%%%%%%%%%%%%%%%%%%%%%%%%%%%%%%%%
    {}={}&
(\frac{1-\kappa}{1+\kappa})^2 + 
\frac{\kappa (1+\kappa)c_{j+1}}{2}
\Big[
       c_{j+1}^2 + 2 c_{j+1} \frac{1+\kappa^3}{\kappa (1+\kappa)^2} - 4 \frac{1+\kappa^3}{\kappa (1+\kappa)^3}
\Big]. 
\end{aligned}
$$%
For the expression within the square bracket, by exploiting its monotonicity with respect to $c_{j+1}$ and evaluating it in $\frac{1}{2}$ we have that it is negative for any $c_{j+1} \in [0, \frac{1}{2}]$ and $\kappa \in (0,1)$. 
Therefore, we get $\frac{c_{j+1}}{c_j} 
    {}\leq{}
(\frac{1-\kappa}{1+\kappa})^2$. By telescoping for $j=0,\dots,k-2$, this implies $c_{k-1} \leq c_0 (\frac{1-\kappa}{1+\kappa})^{2(k-1)} \leq (\frac{1-\kappa}{1+\kappa})^{2k}$, hence $s_{k-1} {}\geq{} \frac{2}{1+\kappa} - (\frac{1-\kappa}{1+\kappa})^{2k}$. \qed%
\end{enumerate}%
\end{proof}%
%
% \begin{lemma}\label{lemma:var_sequence_kappa_m1}
%     If $\kappa = -1$, then $\lambda_{\kappa} (\gamma_{i-1} L) = 2 \cos \left( \frac{1}{3} \arccos( \frac{\gamma_{i-1} L}{2} ) \right)$.
% \end{lemma}%
% \begin{proof}
%     For simplicity of notation, further on we use $s_{k} = \gamma_{i-1} L$ and $s_{+} = \lambda_{-1} (s_{k})$. With $\kappa = -1$ in \eqref{eq:def_max_next_stepsize}, $s_{+}$ is the unique solution belonging to the interval $[1,2)$ of the third order equation:
%     \begin{align*}
%         s_{+} (3-s_{+}^2)
%             +
%         s_{k} = 0.
%     \end{align*}
%     One can check that $s_{+} = 2 \cos [ \frac{1}{3} \arccos \big( \frac{s_{k}}{2} \big) ]$ is one of the roots. It is well defined since, with $s_{k} \in (0,2)$, it holds $\frac{1}{3} \arccos( \frac{s_{k}}{2} ) \in (0, \frac{\pi}{6}]$, hence $s_{+} \in [\sqrt{3}, 2)$.
% \qed \end{proof}
%
%
%%
\subsection{Performance bounds for dynamic stepsizes}\label{proof:prop:tight_rate_all_kappas_increasing_sequence}
In this section we give the proof of \cref{prop:tight_rate_all_kappas_increasing_sequence}, within which we also show the results from 
\cref{thm:convergence_rates_strg_convex_variable_stepsizes} (for strongly convex function), \cref{corollary:convergence_rates_convex_variable_stepsizes} (for convex functions) and 
\cref{thm:GM_hypo_rate_sequence} (for weakly convex functions).%
\begin{proof}[\textbf{\cref{prop:tight_rate_all_kappas_increasing_sequence}}]
    By telescoping \eqref{eq:one_step_rate_h_geq_1} with stepsizes $\gamma_i L$:%
\begin{align}\label{eq:telescoped_N4SD}
\begin{aligned}
    f_0 - f_N
        {}\geq {}&
    \frac{\gamma_{N-1} L}{2-\gamma_{N-1} L -\gamma_{N-1} \mu} \frac{\|g_{N}\|^2}{2L} {}+{}
    \frac{\gamma_{0} L[(2-\gamma_{0} L)(2-\gamma_{0} \mu)-1]}{2-\gamma_{0} L -\gamma_{0} \mu}  \frac{\|g_{0}\|^2}{2L} {}+{}
    \\ & {}+{}
         \sum_{i=1}^{N-1} \frac{\|g_{i}\|^2}{2L}
     \Big[ 
        \frac{\gamma_{i-1} L}{2-\gamma_{i-1} L -\gamma_{i-1} \mu}  {}+{}
        \frac{\gamma_{i} L[(2-\gamma_i L)(2-\gamma_i \mu)-1]}{2-\gamma_{i} L -\gamma_{i} \mu}
    \Big].
\end{aligned}
\end{align}%
Given that all weights are nonnegative (the fact which we prove below), the bound \eqref{eq:tight_rate_all_kappas_increasing_sequence} is obtained by taking the minimum gradient norm:%
\begin{align*}
\begin{aligned}
    f_0 - f_N
        {}\geq {}&
    \min\limits_{0 \leq i \leq N} \left\{\frac{\|g_i\|^2}{2L}\right\}    
         \sum_{i=0}^{N-1} 
    \frac{\gamma_{i} L (2-\gamma_i L)(2-\gamma_i \mu)}{2-\gamma_{i} L -\gamma_{i} \mu}.
\end{aligned}
\end{align*}%
For $\kappa < 0$, we define $\tilde{N} \coloneqq \argmax_{0\leq i \leq N-1} \{ s_i \leq \gstar[(\kappa)] \}$. \\
\noindent\textbf{Case (i):} $\gamma_i L = s_{i}(\kappa)$ for $i=0,\dots,N-1$. 
This corresponds to the (strongly) convex case ($\kappa \in [0,1)$) or weakly convex ($\kappa < 0$) with $N \leq \tilde{N}$. By the definition of $s_{i}(\kappa)$, all weights but the one of $\|g_N\|$ vanish and the inequality \eqref{eq:telescoped_N4SD} reduces to:%
\begin{align}\label{eq:rate_nonconvex_dynamic_schedule_N_leq_Ntilde}
\begin{aligned}
    f_0 - f_N
        {}\geq {}&
    \frac{s_{N-1}}{2-s_{N-1}(1+\kappa)} \frac{\|g_{N}\|^2}{2L},
\end{aligned}
\end{align}%
which is of type \eqref{eq:I^k_N_general_ineq} and leads to the rate \eqref{eq:wc_rate_var_stepsizes_strg_convex} from \cref{thm:convergence_rates_strg_convex_variable_stepsizes}. \\
Moreover, in the convex case the expression of $s_k(0)$ can be computed analytically from \eqref{eq:def_sequence_stepsizes} as the root belonging to $[1,2)$ of:%
    $$
        \frac{s_{k}(0) (3 - 2s_{k}(0))}{2-s_{k}(0)}
            +
        \frac{s_{k-1}(0)}{2-s_{k-1}(0)} = 0.
$$% \\%
This solution is $s_k(0) = \frac{3-2s_{k-1}(0) + \sqrt{9 - 4 s_{k-1}(0)}}{2(2-s_{k-1}(0) )}$, which equivalently rewrites as in the hypothesis of \cref{corollary:convergence_rates_convex_variable_stepsizes}. %

\noindent\textbf{Case (ii):} $\gamma_i L = \min \{ s_i(\kappa) \,,\, \gstar[(\kappa)] \}$ for $i=0,\dots,N-1$. This case corresponds to $\kappa < 0$ and $N > \tilde{N}$. Then inequality \eqref{eq:telescoped_N4SD} becomes
\begin{align}\label{eq:dynamic_stepsizes_nonconvex}
\begin{aligned}
    f_0 - f_N
        {}\geq {}&
    \left[\frac{s_{\tilde{N}}}{2-s_{\tilde{N}}(1+\kappa)} {}+{}
    \frac{\gstar[] [(2-\gstar[])(2-\kappa \gstar[])-1]}{2-\gstar[](1+\kappa)} \right] \frac{\|g_{\tilde{N}+1}\|^2}{2L} {}+{}
    \\ &
    \sum_{i=\tilde{N}+2}^{N-1} 
     \Big[
        \frac{\gstar[](2-\gstar[])(2-\kappa \gstar[])}{2-\gstar[](1+\kappa)} \frac{\|g_{i}\|^2}{2L}
    \Big]  +
    \frac{\gstar[]}{2-\gstar[](1+\kappa)}  \frac{\|g_{N}\|^2}{2L}.
\end{aligned}
\end{align}%
It remains to prove the nonnegativity of $\|g_{\tilde{N}+1}\|^2$'s weight, namely:
\begin{align*}
    \frac{s_{\tilde{N}}}{2-s_{\tilde{N}}(1+\kappa)} {}+{}
    \frac{\gstar[] [(2-\gstar[])(2-\kappa \gstar[])-1]}{2-\gstar[](1+\kappa)} {}\geq{} 0.
\end{align*}
By definition of \eqref{eq:def_sequence_stepsizes}, it holds:
\begin{align*}
    \frac{s_{\tilde{N}}}{2-s_{\tilde{N}}(1+\kappa)} 
        {}={}
    -\frac{s_{\tilde{N}+1} [(2-s_{\tilde{N}+1})(2-\kappa s_{\tilde{N}+1})-1] }{2-s_{\tilde{N}+1}(1+\kappa)},
\end{align*}
with $s_{\tilde{N}+1} > \gstar[(\kappa)]$. It remains to demonstrate 
\begin{align*}
    -\frac{s_{\tilde{N}+1} [(2-s_{\tilde{N}+1})(2-\kappa s_{\tilde{N}+1})-1] }{2-s_{\tilde{N}+1}(1+\kappa)}
        {}\geq{}
    -\frac{\gstar[] [(2-\gstar[])(2-\kappa \gstar[])-1]}{2-\gstar[](1+\kappa)},
\end{align*}
which results from the monotone decreasing of the function $\sigma{:} [1, 2) {\rightarrow} \mathbb{R}$, $
    \sigma(t)
        {{}\coloneqq{}}
    \frac{t [(2{-}t)(2{-}\kappa t){-}1] }{2{-}t(1{+}\kappa)}.
$
This holds because 
$
    \frac{d \sigma (t)}{dt} 
        {}={}
    {- }\frac{2(1-\kappa t)(t-1)[1 + 2-t(1+\kappa)]}{[2-t(1+\kappa)]^2} \leq 0.
$ \\
%
%
%
% \vspace{2em}
By taking the minimum gradient norm in \eqref{eq:dynamic_stepsizes_nonconvex}, with $\pigLgmu$ defined in \eqref{eq:pi_hi_rate_GM}:%
\begin{align}\label{eq:rate_nonconvex_dynamic_schedule_N_geq_Ntilde}
\begin{aligned}
    f_0 - f_N
        {}\geq {}&
    \min\limits_{\tilde{N}+1 \leq i \leq N} \Big\{\frac{\|g_i\|^2}{2L}\Big\}
    \Big[
    \frac{s_{\tilde{N}}}{2-s_{\tilde{N}}(1+\kappa)} {}+{}
    p(\gstar[], \kappa \gstar[]) (N-\tilde{N}-1)    
    \Big].
\end{aligned}
\end{align}% \\%
Since sequence $\{s_k\}$ is increasing, the inequality $ \frac{s_{\tilde{N}}}{2-s_{\tilde{N}}(1+\kappa)} \geq \frac{s_{N-1}}{2-s_{N-1}(1+\kappa)} $ is equivalent to $\tilde{N} \geq N-1$. Consequently, the merged \rhs of \eqref{eq:rate_nonconvex_dynamic_schedule_N_leq_Ntilde} and \eqref{eq:rate_nonconvex_dynamic_schedule_N_geq_Ntilde} can be expressed as a maximum:%
\begin{align*}
  \max \Big\{ 
    \frac{s_{N{-}1}}{2{-}s_{N{-}1}(1+\kappa)}
        \,,\, 
    \frac{s_{\tilde{N}}}{2{-}s_{\tilde{N}}(1{+}\kappa)} {}+{}
    p(\gstar[], \kappa \gstar[]) (N{-}\tilde{N}{-}1)  
 \Big\} {}={} \\
 p(\gstar[], \kappa \gstar[]) N + 
 \max \limits_{0 \leq k \leq N} \Big\{ \frac{s_{k{-}1}}{2{-}s_{k{-}1}(1{+}\kappa)} - p(\gstar[], \kappa \gstar[]) k \Big\},
\end{align*}%
which is exactly as in \eqref{eq:GM_hypo_rate_sequence} from \cref{thm:GM_hypo_rate_sequence}.%
\end{proof}%
%%
%
% \clearpage
\section{Proofs on interpolability of smooth functions}%
\label{appendix:proof_interpolability}%
    %!TEX root = ../../main.tex
\subsection{Proof of interpolation conditions for smooth functions (\texorpdfstring{\cref{thm:interp_hypo}}{Theorem 3.2})%
}\label{subsec:proof_interp_conds_hypo}%
\begin{proof}[\textbf{\cref{thm:interp_hypo}}]%
    The proof follows the same steps outlined in \cite[Theorem 4]{taylor_smooth_2017}, with a minimal extension to accommodate $\mu < 0$. This extension impacts only the first step of the demonstration, known as minimal curvature subtraction, which involves shifting the lower curvature of the function to $0$, effectively converting it into a smooth convex function. This adjustment remains valid for $\mu < 0$, as showed further on. By applying the \textit{gradient inequality} for smooth convex functions: $h(x) {} \geq {} h(y) + \left \langle \nabla h(y), x-y \right \rangle, \, \forall x,y \in \mathbb{R}^d$ for $h = \frac{L}{2}\|\cdot\|^2-f$ and $h = f - \frac{\mu}{2}\|\cdot\|^2$, we obtain:%
    \begin{align}\label{eq:def_hypoconvex_smooth}
    \tfrac{\mu}{2} \|x-y\|^2 
        {} \leq {} 
    f(x) - f(y) - \langle \nabla f(y), x-y \rangle 
        {} \leq {} 
    \tfrac{L}{2} \|x-y\|^2 \quad \forall x, y \in \mathbb{R}^d.
\end{align}%        
    Let $\bar{f} \in \mathcal{F}_{0,L-\mu}$, $\bar{f}(x){}\coloneqq{}f(x)-\tfrac{\mu}{2}\|x\|^2$, with $\nabla \bar{f}(x)=\nabla f(x)-\mu x$; by replacing $f$ with $\bar{f}$ in the quadratic bounds inequality \eqref{eq:def_hypoconvex_smooth} it holds:%
$$
    % \frac{\mu}{2} \|x-y\|^2 {}\leq{}& f(x) - f(y) - \langle \nabla f(y), x-y \rangle {}\leq{} \frac{L}{2} \|x-y\|^2 \nonumber \\
    % & \frac{\mu}{2} \|x-y\|^2 {}\leq{} \big(f(x)-\frac{\mu}{2} \|x\|^2\big) - \big(f(y)-\frac{\mu}{2} \|y\|^2\big) - \langle \nabla \big(f(y)-\mu y\big), x-y \rangle {}\leq{} \frac{L}{2} \|x-y\|^2 \nonumber \\
    % \frac{\mu}{2} \|x-y\|^2 {}\leq{}& \bar{f}(x) - \bar{f}(y)- \langle \nabla \bar{f}(y), x-y \rangle + \frac{\mu}{2} \|x-y\|^2 {}\leq{} \frac{L}{2} \|x-y\|^2 \nonumber \\
    0 {}\leq{} \bar{f}(x) - \bar{f}(y)- \langle \nabla \bar{f}(y), x-y \rangle {}\leq{} \tfrac{L-\mu}{2} \|x-y\|^2.
$$
Hence, the set $\{(x_i, g_i, f_i)\}_{i \in \mathcal{I}}$ is $\mathcal{F}_{\mu,L}$-interpolable if and only if the set $\{(x_i, g_i-\mu x_i, f_i-\tfrac{\mu}{2}\|x_i\|^2)\}_{i \in \mathcal{I}}$ is $\mathcal{F}_{0, L-\mu}$-interpolable. 
Subsequently, the proof continues identically to \cite[Theorem 4]{taylor_smooth_2017}.
\qed\end{proof}%%
\subsection{Proof of optimal point's characterization for smooth functions} \label{appendix:proof_charact_min}

\begin{proof}[\textbf{\cref{lemma:characterization_min}}]\label{lemma_proof:characterization_min}%
We extend \cite[Theorem 7]{drori2021complexity} to smooth functions with any lower curvature $\mu$, generalizing from the original convex ($\mu=0$) and non-necessarily ($\mu=-L$) cases. Consider the function%
\begin{align*}
    \begin{aligned}
        Z(y) {}{}{\coloneqq}{}{} \min_{\alpha \in \Delta_{\mathcal{I}}} 
        &\big\{
            \tfrac{L{-}\mu}{2} 
            \big\|
            y {-} \sum_{i \in \mathcal{I}} \alpha_i
             \big[x_i {-} \tfrac{1}{L{-}\mu} \big(g_i {-} \mu x_i \big) \big] 
            \big\|^2 {+} 
            % \right. \\ & \left. 
            % \\ & \qquad
            \sum_{i \in \mathcal{I}} \alpha_i
             \big( f_i {-} \tfrac{\mu}{2} \|x_i\|^2 {-} \tfrac{1}{2(L{-}\mu)} \|g_i{-}\mu x_i\|^2 \big)
        \big\}
    \end{aligned}
\end{align*}%
where $\Delta_{\mathcal{I}}$ is the $|\mathcal{I}|$-dimensional unit simplex:%
\begin{align*}%
    \Delta_{{\mathcal{I}}} {}\coloneqq{} \big\{ {\alpha} \in \mathbb{R}^n: \sum_{i\in\mathcal{I}} \alpha_i = 1, \alpha_i {}\geq{} 0, \forall i \in \mathcal{I} \big\}.%
\end{align*}%
The function $Z$ is the primal interpolation function $W_{\mathcal{T}}^{C}$ as defined in \cite[Definition 2.1]{drori2018_primal_interp_func}. One can check this by replacing in the definition $C \leftarrow \{0\}$, $L \leftarrow L-\mu$ and %
\begin{align*}
    \mathcal{T} \leftarrow \big\{\left(x_i, g_i-\mu x_i, f_i - \tfrac{\mu}{2} \|x_i\|^2\right) \big\}_{i \in \mathcal{I}}.
\end{align*}
The set $\mathcal{T}$ is $\mathcal{F}_{0, L-\mu}$-interpolable. Hence, from \cite[Theorem 1]{drori2018_primal_interp_func} it follows that $Z$ is convex, with upper curvature $L-\mu$, and satisfies
\begin{align*}
    \begin{aligned}
        Z(x_i) = f_i - \frac{\mu}{2} \|x_i\|^2 \quad \text { and } \quad
        \nabla Z(x_i) = g_i - \mu x_i.
    \end{aligned}
\end{align*}
Consider the function
\begin{align*}
    \hat{W}(y) {}\coloneqq{} Z(y) + \tfrac{\mu}{2} \|y\|^2,
\end{align*}
which by curvature shifting belongs to the function class $\mathcal{F}_{\mu,L}$ and satisfies 
\begin{align*}
    \begin{aligned}
        \hat{W}(x_i) = f_i \quad \text { and } \quad
        \nabla \hat{W}(x_i) = g_i.
    \end{aligned}
\end{align*}
Using algebraic manipulations, $\hat{W}$ can be expressed as
\begin{align}\label{eq:hat_W_equiv}
\begin{aligned}
    \hat{W}(y) {}={} \min_{\alpha \in \Delta_{\mathcal{I}}} 
            &\big\{
                \tfrac{L}{2}
                \big\|
                y - \sum\limits_{i \in \mathcal{I}} \alpha_i
                 \big(x_i - \tfrac{1}{L}g_i \big) 
                \big\|^2 + 
                \tfrac{L}{2}\tfrac{\mu}{L -\mu} \big\|\sum\limits_{i \in \mathcal{I}} \alpha_i
                 \big(x_i - \tfrac{1}{L}g_i \big) \big\|^2 +
                \\ & \qquad 
                \sum\limits_{i \in \mathcal{I}} \alpha_i
                 \big( f_i - \tfrac{1}{2L} \|g_i\|^2 - \tfrac{L}{2} \tfrac{\mu}{L-\mu} \big\|x_i - \tfrac{1}{L} g_i \big\|^2 \big)
            \big\}.
\end{aligned}            
\end{align}
Further, we \textit{lower bound} the first squared norm by $0$:%
\begin{align*}%\label{eq:W_hat_upper_bound}
    \begin{aligned}
        \hat{W}(y) {}{}{\geq}{}{} & 
            \min_{\alpha \in \Delta_{\mathcal{I}}} 
            \big\{
                \tfrac{L}{2}\tfrac{\mu}{L{-}\mu} \big\|\sum\limits_{i \in \mathcal{I}} \alpha_i
                 \big(x_i {-} \tfrac{1}{L}g_i \big) \big\|^2 +
                \sum\limits_{i \in \mathcal{I}} \alpha_i
                 \big( f_i {-} \tfrac{1}{2L} \|g_i\|^2 {-} \tfrac{L}{2} \tfrac{\mu}{L{-}\mu} \big\|x_i {-} \tfrac{1}{L} g_i \big\|^2 \big)
            \big\}.
    \end{aligned}
\end{align*}
Applying Jensen inequality for the squared norm, which is a convex function, we get:%
\begin{align*}
    \big\|\sum\limits_{i \in \mathcal{I}} \alpha_i
                 \big(x_i - \tfrac{1}{L}g_i \big) \big\|^2 
                {}\leq{}
    \sum\limits_{i \in \mathcal{I}} \alpha_i
                 \big\|x_i - \tfrac{1}{L} g_i \big\|^2.
\end{align*}
Then for $\mu \leq 0$ we can further lower bound the above inequality and obtain:%
\begin{align*}%\label{eq:W_hat_upper_bound}
    \begin{aligned}
        \hat{W}(y)
            {}\geq{} & 
            % \min_{\alpha \in \Delta_{\mathcal{I}}} 
            % \big\{
            %     \tfrac{L}{2}\tfrac{\mu}{L-\mu} \sum\limits_{i \in \mathcal{I}} \alpha_i\big\|
            %      x_i - \tfrac{1}{L}g_i \big\|^2 +
            %     \sum\limits_{i \in \mathcal{I}} \alpha_i
            %      \big( f_i - \tfrac{1}{2L} \|g_i\|^2 - \tfrac{L}{2} \tfrac{\mu}{L-\mu} \big\|x_i - \tfrac{1}{L} g_i \big\|^2 \big)
            % \big\} \\
            % {}={}&
            \min_{\alpha \in \Delta_{\mathcal{I}}} 
            \big\{
                \sum\limits_{i \in \mathcal{I}} \alpha_i
                 \big( f_i - \tfrac{1}{2L} \|g_i\|^2 \big)
            \big\} %\\
            {}={}
            \min_{i \in \mathcal{I}}\,\,\,\,
                 \big\{ f_i - \tfrac{1}{2L} \|g_i\|^2 \big\} %\\
            {}={}
            f_{i_*} - \tfrac{1}{2L} \|g_{i_*}\|^2.%
    \end{aligned}
\end{align*}
For the \textit{upper bound}, in \eqref{eq:hat_W_equiv} we take $y{}\coloneqq{}x_{i_*}-\tfrac{1}{L}g_{i_*}$ and ${\alpha}={e}_{i_*}$ (the ${i_*}$-th unit vector):%
\begin{align*}%\label{eq:hat_W_lower_bound}
\begin{aligned}
    \hat{W}\big(x_{i_*}-\tfrac{1}{L}g_{i_*}\big) 
                {}&\leq{}
                \tfrac{L}{2}\tfrac{\mu}{L-\mu} \big\|x_{i_*} - \tfrac{1}{L}g_{i_*} \big\|^2 +
                 f_{i_*} - \tfrac{1}{2L} \|g_{i_*}\|^2 - \tfrac{L}{2} \tfrac{\mu}{L-\mu} \big\|x_{i_*} - \tfrac{1}{L} g_{i_*} \big\|^2 
                 \\ &
                 {}={}
                 f_{i_*} - \tfrac{1}{2L} \|g_{i_*}\|^2.
\end{aligned}            
\end{align*}%
Therefore, from both lower and upper bounds it follows
$\hat{W}\big(x_{i_*}-\tfrac{1}{L}g_{i_*}\big) = f_{i_*} - \tfrac{1}{2L} \|g_{i_*}\|^2$.%
\qed
\end{proof}
%
%
% \clearpage %
\section{Tightness proofs (details)}\label{appendix:tightness_proofs}%
% \subsection{A worst-case function in the convex case}
% \begin{proof}[Proof of \cref{corollary:cvx_rate}]\label{proof:cvx_1D_wc_func}
% The lower bounds result by setting $\mu=0$ in \cref{prop:tightness_short_steps}; then for $\gamma_i L \in (0, \frac{3}{2}]$ a worst-case example function is linear between iterates and extended outside by quadratics of curvature $L$; with $U^2 = \tfrac{f_0-f_N}{\sum_{i=0}^{N-1} \gamma_i}$ it reads:%
% \begin{align}\label{eq:worst-case_func_example_cvx}    
%     f(x)
%         {}={} 
%     \left\{
%     \arraycolsep=8pt
%     \def\arraystretch{1.25}
%     \begin{array}{ll}
%         \tfrac{L}{2} (x+\tfrac{U}{L})^2 - \tfrac{U^2}{2L}
%         & 
%         \text{ if }x \in \left(-\infty, x_N\right]; \\
%         U \left(x - x_{i+1}   \right) + f_{i+1} & 
%         \text{ if }x \in [x_{i+1}, x_i]; \\
%         \tfrac{L}{2} \left(x-x_0\right)^2 + U \left(x - x_{0} \right) + f_{0}
%         &
%         \text{ if }x \in \left[x_{0}, \infty\right).
%     \end{array}
%     \right.
% \end{align}%
% %
% For constant step-sizes $\gamma L \in (\frac{3}{2}, 2)$, a worst-case function is Huber-like, namely the maximum of: (i) the linear function \eqref{eq:worst-case_func_example_cvx} when $\gamma L \in (\tfrac{3}{2}, \gbari[N][(0)])$ and (ii) the quadratic function $\tfrac{L}{2}\|x\|^2$  when $\gamma L \in (\gbari[N][(0)], 2)$ (see \cref{prop:tightness_linear_regime}).
% \qed\end{proof}
%
%
%!TEX root = ../../main.tex
% \clearpage
\subsection{Sublinear regimes with stepsizes \texorpdfstring{$\gamma_i L \in (1, \gbari[1][]]$}{in first interval}} 
\begin{proof}[\textbf{\cref{prop:nec_eq_conds_variable_gamma}}]\label{proof:nec_eq_conds_variable_gamma_mid_steps}
    The equality case from \cref{thm:wc_GM_hypo}'s proof implies all gradient norms equal with each other as defined in \eqref{eq:defs_U_tightness_long_h_ct}. By replacing this expression in the equality conditions from \cref{lemma:sufficient_decrease_h_geq_1}, identity  \eqref{eq:def_ci} results. Identity \eqref{eq:cond_fi_medium_steps} is implied by the equality case in \eqref{eq:one_step_rate_h_geq_1}:
$$
    f_i - f_{i+1} 
        {}={}
    \frac{\gamma_i L(2-\gamma_i \mu)(\text{\oldtext{\color{gray}{\sout{$1$}}}%
    }{\shownewtext{2}}-\gamma_{i} L)}{2-\gamma_i L - \gamma_i \mu} \frac{U^2}{2L}.
$$
The normalized inner product after one iteration is $\cgLgmu[\gamma_i] \in (-1,1)$, therefore $\thetagLgmu$ is well defined in \eqref{eq:def_theta_i}. Moreover, recurrence \eqref{eq:def_reccurence_gi_medium_steps} satisfies the necessary condition \eqref{eq:def_ci} since
$$
    \langle g_i , g_{i+1} \rangle
        {}={}
    U^2 \cos \big( (-1)^i \thetagLgmu[\gamma_i] \big) 
        {}={}
    \shownewtext{U^2} \, \cgLgmu[\gamma_i].
$$
\eqref{eq:cond_*_mid_steps} holds because all gradient norms are equal and the function values decrease.\qed
\end{proof}

\begin{proof}[\textbf{\cref{prop:tightness_long_ct_steps} (cont.)}]\label{appendix:proof:tightness_long_steps_ct} 
We check that all interpolation conditions \eqref{eq:distance-N_interpolation_conditions_I_ij} and \eqref{eq:distance-N_interpolation_conditions_I_ji} hold $\forall (i,j) \in \{0,1,\dots,N\}$: (without loss of generality $j>i$)%
% { \fontsize{8.85pt}{11pt}\selectfont
$$
\begin{aligned}
    &\Iij:
    f_{i}{-}f_{j} {-} \gamma \big\langle g_{j} \,,\, \sum_{k=i}^{j{-}1} g_k\big \rangle 
        {}{\geq}{}
    \tfrac{1}{2(L{-}\mu)}\Big(\big\|g_{i}-g_{j}\big\|^{2} + 
    \gamma L\, \gamma \mu \big\| \sum_{k=i}^{j-1} g_k\big\|^2-
    2 \gamma \mu \big\langle g_{i}{-}g_{j}, \sum_{k=i}^{j-1} g_k\big\rangle\Big); \\
    %%%%%%%%%%%%%%%%%%%%%%%%%%%%%%%%%%%%%%%%%%%%%%%%%%%%%%%%%%%%
    &\Iji:
    f_{j}{-}f_{i} {+} \gamma \big\langle g_{i} \,,\, \sum_{k=i}^{j{-}1} g_k\big \rangle 
        {}{\geq}{}
    \tfrac{1}{2(L{-}\mu)}\Big(\big\|g_{i}-g_{j}\big\|^{2} + 
    \gamma L\, \gamma \mu \big\| \sum_{k=i}^{j-1} g_k\big\|^2-
    2 \gamma \mu \big\langle g_{i}{-}g_{j}, \sum_{k=i}^{j-1} g_k\big\rangle\Big).
\end{aligned} %  
$$%
Note the same expressions of the \rhs in both inequalities. Let $\Delta{}\coloneqq{}f_0-f_N$ and recall $\cgLgmu=\frac{1+(1-\gamma L)(1-\gamma \mu)}{2-\gamma L - \gamma \mu}$. Then  $1+c = \frac{(2-\gamma L)(2-\gamma \mu)}{2-\gamma L - \gamma \mu}$ and 
\[%\label{eq:def_U_gamma=ct}
    U^2 
        {}{}\coloneqq{}{}
    \frac{\Delta}{\frac{1}{2L}\frac{\gamma L(2-\gamma L)(2-\gamma \mu)}{2-\gamma(L+\mu)} N}
        {}={} \frac{2L \Delta}{(1+c) N}.
\]%
\begin{itemize}[align=right,itemindent=1em,labelsep=2pt,labelwidth=1em,leftmargin=0pt,nosep]%
\item {\textbf{Case $j-i$ even.}} Then $g_i = g_j$ and we use the identity:%
$$
    \big\langle g_{i} \,,\, \sum_{k=i}^{j-1}g_k \big \rangle 
        {}={}
    \frac{U^2(1+c)(j-i)}{2} %
$$%
to rewrite the interpolation inequalities as follows:%
$$
    \begin{aligned}
        \Iij: \quad &
    \frac{\Delta (j-i)}{N} - \frac{\gamma U^2(1+c) (j-i)}{2}
        {}\geq {}
    \frac{\gamma L \gamma \mu}{2(L-\mu)}
    \big\| \sum_{k=i}^{j-1} g_k\big\|^2; \\
    %%%%%%%%%%%%%%%%%%%%%%%%%%%%%%%%%%%%%%%%%%%%%%%%%%%%%%%%%%%%
    \Iji: \quad &
    \frac{\Delta (i-j)}{N} + \frac{\gamma U^2(1+c) (j-i)}{2}
        {}\geq {}
    \frac{\gamma L \gamma \mu}{2(L-\mu)}
    \big\| \sum_{k=i}^{j-1} g_k\big\|^2.
    \end{aligned}
$$%
By replacing the definition of $U^2$, the \lhs in both inequalities equals zero, while the \rhs is non-positive because $\mu \leq 0$; therefore, both interpolation inequalities are valid for $i+j$ even.%
%%%%%
\item {\textbf{Case $j-i$ odd.}} Then $g_i \neq g_j$ and $\langle g_i, g_j \rangle = c U^2$, and we use the following identities:%
$$ \big \langle g_{i} \,,\, \sum_{k=i}^{j-1} g_k \big \rangle {}={}
        \frac{U^2(j-i)(1+c)}{2} {} + {} %
        \frac{U^2(1-c)}{2}; $$%
$$\big \langle g_{j} \,,\, \sum_{k=i}^{j-1} g_k \big \rangle {}={}
        \frac{U^2(j-i)(1+c)}{2} {} - {} %
        \frac{U^2(1-c)}{2};$$%
$$\big\langle g_i - g_j , \sum_{k=i}^{j-1} g_k \big\rangle 
            {}={}U^2 (1-c);$$%
$$\big\| \sum_{k=i}^{j-1}  g_k \big\|^2 
            {}={}
        \frac{U^2}{2} \big[(1+c)(j-i)^2 + (1-c) \big];$$%
$$\big\| g_i - g_j \big\|^2 {}={} 2(1-c) U^2.$$%
After substitutions, the \rhs and \lhs of both inequalities $\Iij$ and $\Iji$ are the same:%
\[
    \begin{aligned}
    % R {} = {} & 
    & \text{\rhs}_{\Iij}
        {}={}
    \text{\rhs}_{\Iji}
        {}={}
    % \\ & \qquad \qquad 
    \frac{U^2}{4(L-\mu)} 
      \big[
        \gamma L \gamma \mu (1+c) (j-i)^2 + (1-c) \big(\gamma L \gamma \mu - 4 \gamma \mu + 4\big)
      \big]
    \end{aligned}
\]%
and, respectively:%
% { \fontsize{8.5pt}{10pt}\selectfont
$$
    \begin{aligned}
    \text{\lhs}_{\Iij}: \,\, &
        \frac{\Delta(j-i)}{N}-\gamma \big\langle g_{j}, \sum_{k=i}^{j-1} g_k \big \rangle
            {}={}
        \frac{\gamma U^2(1-c)}{2} ;
        \\
        %%%%%%%%%%%%%%%%%%%%%%%%%%%%%%%
    \text{\lhs}_{\Iji}: \,\, &
        \frac{\Delta(i-j)}{N}+\gamma \big\langle g_{i}, \sum_{k=i}^{j-1} g_k \big \rangle  
            {}={}
        \frac{\gamma U^2(1-c)}{2} .
    \end{aligned}%
$$%}%
It remains to prove (with $j-i = 1,3,\dots$):%
$$
    \begin{aligned}
    \frac{\gamma U^2(1-c)}{2} 
        {}\geq {}&
   \frac{U^2}{4(L-\mu)} 
      \Big[
        \gamma L \gamma \mu (1+c) (j-i)^2 + (1-c) \big[(\gamma L)(\gamma \mu) - 4 (\gamma \mu) + 4\big]
      \Big],
    \end{aligned}%
$$%
equivalent to%
$$
    0  
        {}\geq {}
    \gamma L \gamma \mu (1+c) (j-i)^2 + (1-c) \big[\gamma L \gamma \mu - 2 \gamma(L+\mu) + 4 \big].
$$%
By using identity%
$$
    (1-c) \big[\gamma L \gamma \mu - 2 \gamma(L+\mu) + 4 \big] = -\gamma L \gamma \mu (1+c),
$$%
the inequality to prove becomes:%
$$
    \begin{aligned}
    0 {}\geq {} \gamma L \gamma \mu (1+c) \big[(j-i)^2 - 1\big].
    \end{aligned}
$$%
This is true since $j-i \geq 1$, $\mu \leq 0$ and $\cgLgmu \in [-1,1]$; the equality case holds for $j=i+1$, i.e., the necessary conditions from equality in the distance-1 interpolation inequalities used to derive the upper bound.%
\end{itemize}%

We have shown that the proposed set of triplets satisfy the interpolation conditions for all pairs $(i,j)$, hence, following \cref{thm:interp_hypo}, it is $\mathcal{F}_{\mu,L}$-interpolable. {}\qed{}
\end{proof}%%
\subsection{Proof of the three-dimensional worst-case example from \texorpdfstring{\cref{conj:3D_triplets}}{Proposition 5.13} }\label{sec:proof:3D_wc_nonconvex}%
\textbf{Preliminaries.} 
We denote 
$\m \coloneqq 1-\gamma \mu \in (1, \infty)$, 
$\n \coloneqq 1-\gamma L \in (-1, 0)$, $\tilde{U} \coloneqq \sqrt{\frac{(1{-}\n^2)(\m^2{-}1)}{\m^2{-}\n^2}}$, and rewrite the expressions defining the triplets as:%
 \begin{align}
    & g_i 
        {}={}
    U \tilde{U}
    \left[
    \begin{array}{c}
       \tfrac{\n^{i-\bar{N}}}{\sqrt{1-\n^2}}  \\ 
       \tfrac{\m^{i-\bar{N}}}{\sqrt{\m^2-1}}  \\
       0
    \end{array}
    \right]
    {}\quad{} \forall i=0,1,\dots,\bar{N}+1; \nonumber \\
    %%%%%%%%%%%%%%%%%%%%%%%%%%%%%%%%%%%%%%%
    & g_{i+1}
        {}={}
    \n g_{i} + (\m-\n) \cos(\theta_{i+1})
    \Big[ \begin{array}{cc}  0  &  \\ &  R(\theta_{i+1}) \end{array} \Big] g_{i} {}\quad{} \forall i=\bar{N}, \dots, N. \label{eq:recurrence_gi} \\ % \\
    %%%%%%%%%%%%%%%%%%%%%%%%%%%%%%%%%%%%%%%
    % & q_{i} {}={} \tan(\theta_i) {}={} (-1)^{i-(\bar{N}+1)} \sqrt{(\m^2-1) \tfrac{1-\n^{2(i-(\bar{N}+1))}}{1-\n^2}}; \nonumber \\ 
    %%%%%%%%%%%%%%%%%%%%%%%%%%%%%%    
        & f_i {}={}
        \Delta
    \tfrac{ N-i + \frac{(1-\n)(1-\m)}{\m-\n} \sum_{j=1}^{\bar{N}-i}\Tinm[j] }
    {N + \frac{(1-\n)(1-\m)}{\m-\n} \sum_{j=1}^{\bar{N}}\Tinm[j]}  {}\quad{} i=0,1,\dots,N. \label{eq:def_fi_3D}
\end{align} 
The recurrence \eqref{eq:recurrence_gi} also works to define $g_{\bar{N}+1}$ due to $q_{\bar{N}+1} = 0$ and $\theta_{\bar{N}+1} = 0$. We consider 
\begin{align}\label{eq:def_Sk}
    S_{k} \coloneqq \sum_{p=0}^{k-(\bar{N}+2)} \n^{2p}
        {}={} \frac{1 - \n^{2(k-(\bar{N}+1))}}{1-\n^2},
\end{align}
satisfying the properties $S_{k+1} = 1 + \n^2 S_{k}$ and $S_{\infty} \coloneqq \lim_{k \rightarrow \infty} S_k = \frac{1}{1-\n^2}$. For each $k = \bar{N}+1, \bar{N}+2, \dots, N$ we have that
\begin{align}\label{eq:exp_qk}
q_{k} {}={} \tan(\theta_k)
   % {}={} (-1)^{k-(\bar{N}+1)} \sqrt{ (\m^2-1) \frac{1-\n^{k-(\bar{N}+1)}}{1-\n^2} }
    {}={} (-1)^{k-(\bar{N}+1)} \sqrt{ (\m^2-1) S_{k} }.
\end{align}%
\paragraph*{1. Proof of reaching the performance bound. {}} 
~We demonstrate that $\min_{0 \leq i \leq \bar{N}} \|g_i\| {}={} U$. \\
Note that:
\begin{align}
    \|g_i\|^2 
        {}={}  
    U^2 \tilde{U}^2 \big( \tfrac{\m^{2(i-\bar{N})}}{\m^2-1} + \tfrac{\n^{2(i-\bar{N})}}{1-\n^2} \big) {}\quad{} i = 0,\dots,\bar{N}+1; \label{eq:grads_nrm_eq_thm} 
\end{align}
in particular, $\|g_{\bar{N}}\|^2 {}={} \|g_{\bar{N}+1}\|^2 {}={} U^2$. For the later, one can use the identity $\frac{\m^2}{\m^2-1} - \frac{\n^2}{1-\n^2} {}={} \frac{1}{\m^2-1} - \frac{1}{1-\n^2}$, obtained by adding $1$ to both ratios. 
We show that all gradients with indices $i \geq \bar{N}$ have norm $U$, namely $\|g_i\|^2 = U^2$ for any $i = \bar{N}, \dots, N$. \\
    %%%%%%%%%%%%%%%%%%%%%%%%%%%%%%%%%%%
Recurrence \eqref{eq:recurrence_gi} implies that the first component of $g_i$ is always $U \tilde{U} \tfrac{\n^{i-\bar{N}}}{\sqrt{1-\n^2}}$, being a scaling with $\n$ of the first component of $g_{i-1}$. 
Let $v_i \coloneqq \begin{bmatrix}
    0 & 1 & 0 \\
    0 & 0 & 1
\end{bmatrix} g_i$ be the vector of the last two entries in $g_i$. Thus showing that $\|g_i\|^2 = U^2$ for $i \geq \bar{N}$ reduces to prove that 
\begin{align}\label{eq:norm_of_vk}
    \|v_i\|^2
    {}={} U^2 \Big(1 - \tilde{U}^2 \frac{\n^{2(i-\bar{N})}}{1-\n^2} \Big)    
    , \,\, \forall i=\bar{N}, \bar{N}+1, \dots, N\shownewtextfinal{.}
\end{align}
We show this by induction. Firstly, it rewrites as
\begin{align}\label{eq:norm_of_vk_equiv_written}
    \|v_i\|^2
    {}={} U^2 \tilde{U}^2 \frac{1 + q_{i+1}^2}{\m^2-1} 
    % {}={} \frac{U^2 \tilde{U}^2}{(\m^2-1) \cos(\theta_{i+1})}
    , \,\, \forall i=\bar{N}, \bar{N}+1, \dots, N.
\end{align}
Cases $i=\{\bar{N}, \bar{N}+1\}$ result from the expressions of $g_{\bar{N}}$ and $g_{\bar{N}+1}$. 
Assuming that \eqref{eq:norm_of_vk} holds for some index $k$, we show its validity for $k+1$. 
The recurrence of the gradients implies
\begin{align}\label{eq:recurrence_vk}
    v_{k+1} {}={} \n v_{k} + (\m-\n) \cos(\theta_{k+1}) R(\theta_{k+1}) v_k.
\end{align}
Taking the norm, expanding the squares, using the properties that the rotation matrix preserves the vector's norm and $\langle v_k, R(\theta_{k+1}) v_k \rangle = \|v_k\|^2 \cos(\theta_{k+1})$ yields    
\begin{align*}
    \|v_{k+1}\|^2 
        {}={}& 
    % \|\n v_{k} + (\m-\n) \cos(\theta_{k+1}) R(\theta_{k+1}) v_k\|^2 \\
    %     {}={}&
    % \n^2 \|v_{k}\|^2 + 
    % (\m-\n)^2 \cos^2(\theta_{k+1}) \| R(\theta_{k+1}) v_k\|^2
    % + 2 \n(\m-\n) \cos(\theta_{k+1}) \langle  v_{k}, R(\theta_{k+1}) v_k \rangle \\
    %     {}={}&
    % \n^2 \|v_{k}\|^2 + 
    % (\m-\n)^2 \cos^2(\theta_{k+1}) \|v_{k}\|^2
    % + 2 \n(\m-\n) \cos^2(\theta_{k+1}) \|v_{k}\|^2 \\
    %     {}={}&
    % \|v_{k}\|^2 
    % {[} \n^2 + 
    % (\m-\n)^2 \cos^2(\theta_{k+1})
    % + 2 \n(\m-\n) \cos^2(\theta_{k+1})
    % {]} \\
    %     {}={}&
    \|v_{k}\|^2 
    {[} \n^2 + 
    (\m^2-\n^2) \cos^2(\theta_{k+1})
    {]}.
\end{align*}%
Using $\cos^2(\theta_{k+1}) = \frac{1}{1+q_{k+1}^2}$ and the identity $\m^2 + \n^2 q_{k+1}^2 {}={} 1 + q_{k+2}^2$ we get%
\begin{align*}
    \|v_{k+1}\|^2
        {}={}
    \|v_{k}\|^2
    \frac{ 1 + q_{k+2}^2 }{1 + q_{k+1}^2},
\end{align*}%
% Replacing the expressions of $\|v_k\|^2$ we get
% \begin{align}\label{eq:norm_of_vk}
%     \|v_{k+1}\|^2
%     {}={} 
%     U^2 \left(1 - \frac{(\m^2-1)(1-\n^2)}{\m^2-\n^2} \frac{\n^{2(k-\bar{N})}}{(1-\n^2)} \right) 
%     \frac{ 1 + q_{k+2}^2 }{1 + q_{k+1}^2}
% \end{align}
% 
where replacing $\|v_k\|$ using \eqref{eq:norm_of_vk_equiv_written} yields the expression in \eqref{eq:norm_of_vk} written for $k+1$.

Since $U^2 = \frac{2L \, \Delta}{ \gamma L \PN }$ and $f_0 - f_N = \Delta$, the triplets reach the performance bound with respect to $f(x_0) - f(x_N)$ from \cref{conj:3D_triplets}. To prove the bound with respect to  $f(x_0) - f_*$, it remains to show that the required assumption in \cref{lemma:decouple_f*} holds, namely
\begin{align}\label{eq:condition_tight_decouple}
    N{}\in{}\argmin_{0 \leq i \leq N} \big\{f_i - \tfrac{1}{2L}\|g_i\|^2\big\}.%
\end{align}%
After replacing the expressions of gradient norms in \eqref{eq:def_fi_3D} and several simplifications we get:
$$
    f_i - \frac{\|g_i\|^2}{2L}
        {}={}
    f_N - \frac{U^2}{2L} + 
    \frac{U^2}{2L} \frac{(1-\n)(1+\n)(1+\m)}{\m+\n}
    \Big[ N - \bar{N} + \sum_{j=1}^{\bar{N}-i} \m^{-2j} \Big], \forall i = 0, \dots, N,
$$
where, with an abuse of notation, the sum is zero for $i \geq \bar{N}$. 
Its minimum is reached for $i=N$ and therefore identity \eqref{eq:condition_tight_decouple} holds.%
%%%%%%%%%%%%%%%%%%%%%%%%%%%%%%%%%%%%%%%%%%%%%%%%%%%%%%%%%%%%%%%%%%%%%%%%%%%%%%%%%%%%%%%%%%%%%%%%%%%%%%%%%%%%%%%%
%%%%%%%%%%%%%%%%%%%%%%%%%%%%%%%%%%%%%%%%%%%%%%%%%%%%%%%%%%%%%%%%%%%%%%%%%%%%%%%%%%%%%%%%%%%%%%%%%%%%%%%%%%%%%%%%
\paragraph*{2. The set of triplets $\mathcal{T}_{\mathcal{I}}$ is $\mathcal{F}_{\mu,L}$-interpolable.{}}
~We show that all interpolation conditions \eqref{eq:Interp_hypoconvex} are satisfied, namely for any $i,j$ with $0\leq i,j \leq N$, $i \neq j$, it holds that $Q_{i,j} \geq 0$, where 
\begin{align*}
    Q_{i,j} 
        {}\coloneqq{}
    f_{i}-f_{j}-\langle g_{j}, x_{i}-x_{j} \rangle 
        {}-{}
    \frac{1}{2L} \|g_{i}-g_{j}\|^{2} {}-{} \frac{\mu}{2L(L-\mu)} \|g_{i}-g_{j} - L (x_i-x_j)\|^2.
\end{align*}%
We divide the analysis into the cases: $i>j$ and $i<j$, showing that $Q_{i,j}$ rewrites as a sum of nonnegative contributions in each of them. The central challenge lies in finding an exact expression for the inner product $\left\langle g_{k+1} - \m g_{k} \,,\, g_{i+1} - \n g_i \right\rangle$ for all $i,k \geq \bar{N}$, which we obtain below in \eqref{eq:eq_mixed_grads_eval_to_use}. \\
\textbf{Preliminaries.} 
For simplicity of notation and without loss of generality, we consider $U=1$ since all gradients are multiplied by $U = \min_{0 \leq i \leq N} \|\nabla f(x_i)\|$. 

Following \eqref{eq:def_Sk}, observe that $q_{\bar{N}+1} = 0$, $q_{\bar{N}+2} = -\sqrt{\m^2-1}$ and $|q_k|$ is monotonically increasing towards $|q_{\infty}| {}={} \sqrt{\frac{\m^2-1}{1-\n^2}}$. Consequently, $|\theta_k| {}={} \arctan(|q_{k}|) \in (0, \frac{\pi}{2})$, $|\theta_k|$ monotonically increases towards $\arctan(|q_{\infty}|)$ and $\sgn( \theta_k ) = (-1)^{k-(\bar{N}+1)}$. \\
For all $k=0,\dots,\bar{N}+1$, the gradients read as $g_k {}={} \tilde{U} [\frac{\n^{k-\bar{N}}}{\sqrt{1-\n^2}}, \frac{\m^{k-\bar{N}}}{\sqrt{\m^2-1}}, 0]^T$ and we have
\begin{align}\label{eq:zero_gradients_mix}
    \left\langle g_{k+1} - \m g_{k} \,,\, g_{i+1} - \n g_i \right\rangle 
        {}={} \langle [*,0,0]^T \,,\, [0,*,*]^T \rangle {}={} 0, \,\, \forall k \leq \bar{N}, \, i = 0,\dots,N-1,
\end{align}
where $*$ denotes for the entries \shownewtextfinal{of the} three dimensional vectors multiplied by zero.  
Only the first component $g_{k+1} - \eta g_{k}$ is nonzero (for $k \leq \bar{N}$), while the first component in $g_{i+1} - \n g_{i}$ is always zero for any $i=0,\dots,N-1$. \\
\medskip \\
\textbf{Calculation of $\left\langle g_{k+1} - \m g_{k} \,,\, g_{i+1} - \n g_i \right\rangle$ for any $i,k \geq \bar{N}$.} 
Identity \eqref{eq:norm_of_vk_equiv_written} implies for all $k=\bar{N}, \bar{N}+1, \dots, N-1$ that
$
    \|v_k\| \cos(\theta_{k+1})
    {}={} \frac{\tilde{U}}{\sqrt{\m^2-1}}
$ and thus:%
\begin{align}
    %%%%%%%%%%%%%%%%%%%%%%%%%%%%%%%%%%%%%%%%%%%%%%%%%%%%%%%%%%%%%%%%%%%%%%%%%%%%%%%%%%%%
    v_{k+1} - \n v_{k} 
        % &{}={} (\m-\n) \cos(\theta_{k+1}) R(\theta_{k+1}) v_k \nonumber \\
        &{}={} (\m-\n) \frac{\tilde{U}}{\sqrt{\m^2-1}}  R(\theta_{k+1}) \frac{v_k}{\|v_k\|}; \label{eq:vip_minus_nvi} \\
    %%%%%%%%%%%%%%%%%%%%%%%%%%%%%%%%%%%%%%%%%%%%%%%%%%%%%%%%%%%%%%%%%%%%%%%%%%%%%%%%%%%%
    v_{k+1} - \m v_{k} 
        &{}={} -(\m-\n) \big[I_2 - \cos(\theta_{k+1}) R(\theta_{k+1})\big] v_k \nonumber \\
        &{}\stackrel{\eqref{eq:norm_of_vk_equiv_written}}{=}{} -(\m-\n) \frac{\tilde{U}}{\sqrt{\m^2-1}} q_{k+1} R\big({-}\frac{\pi}{2}\big) R(\theta_{k+1}) \frac{v_k}{\|v_k\|}. \label{eq:vip_minus_mvi}
\end{align}
Then we have
\begin{align}\label{eq:intermediate_eq_mixed_grads_eval}
    \left\langle g_{k+1} - \m g_{k} \,,\, g_{i+1} - \n g_i \right\rangle
        {}={}
    \frac{-(\m-\n)^2 \tilde{U}^2}{\m^2-1}
    q_{k+1} \cos(\beta_{k,i}),
\end{align}
where $\beta_{k,i} {}\coloneqq{} \angle(R(-\frac{\pi}{2}) R(\theta_{k+1}) \frac{v_k}{\|v_k\|} \,,\, R(\theta_{i+1}) \frac{v_i}{\|v_i\|})$. We normalize all vectors to transform the problem to maneuvering angles on the unit circle. 
Let $\phi(y)$ be the oriented angle of a two-dimensional vector $y$. We express the angles as follows:%
\shownewtextfinal{%
\[%
    \begin{aligned}
        \beta_{k,i} & {}={} \phi\Big(R(\theta_{i+1}) \frac{v_i}{\|v_i\|}\Big) - \phi\Big( R(-\frac{\pi}{2}) R(\theta_{k+1})\Big); \\
        %%%%%%%%%%%%%%%%%%%%%%%%%%%%%%%%%%%%%%%%%%%%%%%%%%%%%%%%%%%%%%%%%%%%%
        -\frac{\pi}{2} &{}={} \phi\Big( R(-\frac{\pi}{2}) R(\theta_{k+1}) \frac{v_k}{\|v_k\|}\Big) - \phi\Big(R(\theta_{k+1}) \frac{v_k}{\|v_k\|}\Big); \\
        %%%%%%%%%%%%%%%%%%%%%%%%%%%%%%%%%%%%%%%%%%%%%%%%%%%%%%%%%%%%%%%%%%%%%
        \theta_{k+1} &{}={} \phi\Big(R(\theta_{k+1}) \frac{v_k}{\|v_k\|}\Big) - \phi\Big(\frac{v_k}{\|v_k\|}\Big); \\
        %%%%%%%%%%%%%%%%%%%%%%%%%%%%%%%%%%%%%%%%%%%%%%%%%%%%%%%%%%%%%%%%%%%%%
        \theta_{i+1} &{}={} \phi\Big(R(\theta_{i+1}) \frac{v_i}{\|v_i\|}\Big) - \phi\Big(\frac{v_i}{\|v_i\|}\Big).
    \end{aligned}
\]%
}%
Using these identities yields
$$
\beta_{k,i} 
% {}={} \phi(R(\theta_{i+1}) \frac{v_i}{\|v_i\|}) - \phi( R(-\frac{\pi}{2}) R(\theta_{k+1}) \frac{v_k}{\|v_k\|}) 
    {}={} 
\frac{\pi}{2} + \theta_{i+1} - \theta_{k+1} + \phi\big(\frac{v_i}{\|v_i\|}\big) - \phi\big(\frac{v_k}{\|v_k\|}\big).
$$%
For any $l \geq \bar{N}$, let $\delta_{l,l+1} \coloneqq \phi\big(\frac{v_{l+1}}{\|v_{l+1}\|}\big) - \phi\big(\frac{v_l}{\|v_l\|}\big)$ be the oriented angle from $\frac{v_l}{\|v_l\|}$ to $\frac{v_{l+1}}{\|v_{l+1}\|}$ and $\epsilon_{l,l+1} \coloneqq \theta_{l+2} - \theta_{l+1} + \delta_{l,l+1}$. 
We rewrite $\beta_{k,i}$ as
$$
\beta_{k,i} 
    {}={} 
\frac{\pi}{2} 
+ \sum_{l=\bar{N}}^{i-1} \epsilon_{l,l+1} 
- \sum_{l=\bar{N}}^{k-1} \epsilon_{l,l+1}
    {}={}
\begin{cases}
    \frac{\pi}{2} + \sum_{l=k}^{i-1} \epsilon_{l,l+1}, & \text{ if $i \geq k$}; \\
    \frac{\pi}{2} - \sum_{l=i}^{k-1} \epsilon_{l,l+1}, & \text{ if $i < k$}.
\end{cases}
$$
We compute a closed form of $\sin(\epsilon_{l,l+1})$ with $l \geq \bar{N}$:
\begin{align*}
    \sin(\epsilon_{l,l+1})
                       {}={}& \sin(\delta_{l,l+1}) \cos( \theta_{l+2} - \theta_{l+1} ) + 
                              \cos(\delta_{l,l+1}) \sin( \theta_{l+2} - \theta_{l+1} ).
\end{align*} 
Using that $q_l = \tan(\theta_l)$, we get:
\begin{align*}
    \begin{aligned}
        % & 
        %%%%%%%%%%%%%%%%%%%%%%%%%%%%%%%%%%%%%%%%%%%%%%%%%%%%%%%%%%%%%%%%%%%%%%%%%%%%%%%%%%%
        %%%%%%%%%%%%%%%%%%%%%%%%%%%%%%%%%%%%%%%%%%%%%%%%%%%%%%%%%%%%%%%%%%%%%%%%%%%%%%%%%%%
        & \cos( \theta_{l+2} - \theta_{l+1} ) 
            {}={} 
        \frac{1+q_{l+1}q_{l+2}}{\sqrt{(1+q_{l+2}^2)(1+q_{l+1}^2)}}; \\      
        %%%%%%%%%%%%%%%%%%%%%%%%%%%%%%%%%%%%%%%%%%%%%%%%%%%%%%%%%%%%%%%%%%%%%%%%%%%%%%%%%%%
        & \sin( \theta_{l+2} - \theta_{l+1} ) 
            {}={} 
        \frac{q_{l+2}-q_{l+1}}{\sqrt{(1+q_{l+2}^2)(1+q_{l+1}^2)}}.
    \end{aligned}
\end{align*}
Since
$\cos(\delta_{l,l+1})  {}={} \langle \frac{v_l}{\|v_l\|}, \frac{v_{l+1}}{\|v_{l+1}\|} \rangle$ and
$\sin(\delta_{l,l+1})  {}={} \frac{v_l}{\|v_l\|} \times \frac{v_{l+1}}{\|v_{l+1}\|} $, 
we obtain
\begin{align*}
    \begin{aligned}
        % & 
        %%%%%%%%%%%%%%%%%%%%%%%%%%%%%%%%%%%%%%%%%%%%%%%%%%%%%%%%%%%%%%%%%%%%%%%%%%%%%%%%%%%
        %%%%%%%%%%%%%%%%%%%%%%%%%%%%%%%%%%%%%%%%%%%%%%%%%%%%%%%%%%%%%%%%%%%%%%%%%%%%%%%%%%%        
        \cos(\delta_{l,l+1}) 
            {}={} \frac{\m + \n q_{l+1}^2}{\sqrt{(1+q_{l+2}^2)(1+q_{l+1}^2)}}; \,\,
        %%%%%%%%%%%%%%%%%%%%%%%%%%%%%%%%%%%%%%%%%%%%%%%%%%%%%%%%%%%%%%%%%%%%%%%%%%%%%%%%%%%
        \sin(\delta_{l,l+1}) 
        %%%%%%%%%%%%%%%%%%%%%%%%%%%%%%%%%%%%%%%%%%%%%%%%%%%%%%%%%%%%%%%%%%%%%%%%%%%%%%%%%%%
            {}={} \frac{(\m-\n) q_{l+1}}{\sqrt{(1+q_{l+2}^2)(1+q_{l+1}^2)}}.
    \end{aligned}
\end{align*}
Replacing in the expansion of $\sin(\epsilon_{l,l+1})$, after simplifications we get:
\begin{align*}
    \sin( \epsilon_{l,l+1} )
    %     {}={}&
    % \frac{(\m-\n) q_{l+1} (1+q_{l+1}q_{l+2}) + (\m + \n q_{l+1}^2)(q_{l+2}-q_{l+1}) }{(1+q_{l+2}^2)(1+q_{l+1}^2)} \\
%%%%%%%%%%%%%%%%%%%%%%%%%%%%%%%%%%%%%%%%%%%%%%%%%%%%%%%%%%%%%%%%%%%%%%%%%%%%%%%%%%%
        {}={}&
    % \frac{(\m q_{l+2} - \n q_{l+1}) (1+q_{l+1}^2)}{(1+q_{l+2}^2)(1+q_{l+1}^2)}
    %     {}={} 
    \frac{\m q_{l+2} - \n q_{l+1}}{1+q_{l+2}^2}.
\end{align*}
Substituting $q_{l+1}$ and $q_{l+2}$ and amplifying by the conjugate square root yields:
\begin{align*}
    \sin( \epsilon_{l,l+1} )
    %     {}={}&
    % \sqrt{\m^2-1}
    % \frac{\m (-1)^{l+2-(\bar{N}+1)} \sqrt{S_{l+2}} - \n (-1)^{l+1-(\bar{N}+1)} \sqrt{S_{l+1}}}{1 + (\m^2-1) S_{l+2}} \\
    %     {}={}& 
    % (-1)^{l-(\bar{N}+1)}\sqrt{\m^2-1}
    % \frac{\m \sqrt{S_{l+2}} - \sqrt{S_{l+2}-1}}{1 + (\m^2-1) S_{l+2}} \\
        {}={}&     
    (-1)^{l-(\bar{N}+1)} \frac{\sqrt{\m^2-1}}{\m \sqrt{S_{l+2}} - \n \sqrt{S_{l+1}}}.
\end{align*}
For simplicity of notation, we define the angle $\alpha_{l}$ as
\begin{align}\label{eq:def_angle_alpha}%\tag{$\alpha_l$}
    \alpha_{l} 
        \coloneqq 
    % \arcsin \left( \frac{\sqrt{\m^2-1}}{\m \sqrt{S_{l+2}} + \sqrt{S_{l+2}-1}} \right)
    %     {}={} 
    \arcsin \left( \frac{\sqrt{\m^2-1}}{\m \sqrt{S_{l+2}} - \n \sqrt{S_{l+1}}} \right), \,\, \forall l \geq \bar{N}.
\end{align}
Since $\rho < 0$, $\eta > 1$ and $S_{l}$ is monotonically increasing with $l$, we have that $\alpha_l$ is monotonically decreasing with $l$, hence $\alpha_l \in [\alpha_{\infty}, \alpha_{\bar{N}}] \subset (0, \frac{\pi}{2})$, where $\alpha_{\bar{N}} = \arcsin\big(\sqrt{1 - \frac{1}{\m^2}} \big)$ and $\alpha_{\infty} = \arcsin\big(\frac{\sqrt{(\m^2-1)(1-\n^2)}}{\m-\n}\big)$. \\%
Using $q_{k+1} = \sqrt{(\m^2-1) S_{k+1}} (-1)^{k-\bar{N}}$, after substituting everything into \eqref{eq:intermediate_eq_mixed_grads_eval} we get%
\begin{align}
\begin{aligned}
    & \left\langle g_{k+1} {-} \m g_{k} \,,\, g_{i+1} {-} \n g_i \right\rangle
        {}={} \\
    & \qquad \frac{{-} (\m{-}\n)^2 \tilde{U}^2}{\sqrt{\m^2{-}1}}
    \sqrt{S_{k+1}} \cdot
    \begin{cases}
        \sin( \sum_{l=k}^{i-1} ({-1})^{l{-}k} \alpha_l), & \text{ if } \bar{N} \leq k \leq i \leq N-1; \\
        \sin( \sum_{l=i}^{k-1} ({-1})^{k{-}1{-}l} \alpha_l), & \text{ if } \bar{N} \leq i < k \leq N-1.
    \end{cases} \label{eq:eq_mixed_grads_eval_to_use}%
\end{aligned}
\end{align}%
All distance-1 interpolation inequalities $Q_{i,j}$ with $|i-j|=1$ are involved in the proofs and thus hold with equality. This fact can also be checked by direct substitutions of our triplets and using that the inner product between consecutive gradients is known in closed form; for $i \geq \bar{N}$ we have \(\langle g_{\bar{N}}, g_{\bar{N}+1} \rangle = c U^2\), with $c = \frac{1+\eta \rho}{\m + \n}$ (as in \eqref{eq:def_ci}). %
\smallskip \\
\textbf{Case $i>j$.} Let us consider a fixed $j \in \{0, 1, \dots, N-1\}$. For each $i \in \{j+1, \dots, N\}$, we define %
\begin{align*}%\label{def:Delta_Q_ij_i>j}
    \Delta Q_{i} \coloneqq Q_{i+1,j} - Q_{i,j} - Q_{i+1,i}.
\end{align*}%
The distance-1 interpolation inequalities hold exactly, hence $Q_{i+1,i} = 0$ and $Q_{j+1,j} = 0$. 
Therefore to prove that $Q_{i,j} \geq 0$ for all $i>j$ it is sufficient to show that $\Delta Q_{i} \geq 0$, as $Q_{i,j}$ is obtained as a sum of positive increments. \\
Substituting the expressions of $Q_{i+1,j}$, $Q_{i,j}$ and $Q_{i+1,i}$ we subsequently obtain:%
\begin{align*}%\label{def:Delta_Q_ij_i>j}
    \Delta Q_{i} 
    {}={}& \frac{1}{L-\mu} 
        \left\langle g_i - g_j - \mu(x_i - x_j) \,,\, g_i - g_{i+1} - L(x_i - x_{i+1}) \right\rangle \\
    %%%%%%%%%%%%%%%%%%%%%%%%%%%%%%%%%%%%%%%%%%%%%%%%%%%%%%%%%%%%%%%%%%%%%%%%%%%%%%%%%%%%%%%%%%%
    % {}={}& \frac{1}{L-\mu} 
    %     \left\langle g_i - g_j + \gamma \mu \sum_{k=j}^{i-1} g_k \,,\, (1-\gamma L) g_i - g_{i+1} \right\rangle \\
    %%%%%%%%%%%%%%%%%%%%%%%%%%%%%%%%%%%%%%%%%%%%%%%%%%%%%%%%%%%%%%%%%%%%%%%%%%%%%%%%%%%%%%%%%%%
    {}={}& \frac{1}{L-\mu} 
        \langle g_j - g_i - \gamma \mu \sum_{k=j}^{i-1} g_k \,,\, g_{i+1} - (1-\gamma L) g_i \rangle \\
    %%%%%%%%%%%%%%%%%%%%%%%%%%%%%%%%%%%%%%%%%%%%%%%%%%%%%%%%%%%%%%%%%%%%%%%%%%%%%%%%%%%%%%%%
    %%%%%%%%%%%%%%%%%%%%%%%%%%%%%%%%%%%%%%%%%%%%%%%%%%%%%%%%%%%%%%%%%%%%%%%%%%%%%%%%%%%%%%%%
    % {}={}& \frac{-1}{L-\mu} \sum_{k=j}^{i-1}
    %     \left\langle g_{k+1} - (1-\gamma \mu) g_{k} \,,\, g_{i+1} - (1-\gamma L) g_i \right\rangle \\
    {}={}& \frac{-\gamma }{\m-\n} \sum_{k=j}^{i-1}
        \left\langle g_{k+1} - \m g_{k} \,,\, g_{i+1} - \n g_i \right\rangle,
\end{align*}
where in the second step we use the gradient descent iteration. 
For any $k \leq \bar{N}$, from \eqref{eq:zero_gradients_mix} we have that $\left\langle g_{k+1} - \m g_{k} \,,\, g_{i+1} - \n g_i \right\rangle = 0$. 
Hence it is sufficient to focus on the range $j \geq \bar{N}+1$ (thus $i \geq \bar{N}+2$). 
Using \eqref{eq:eq_mixed_grads_eval_to_use} in the case $i \geq k$, the expression of $\Delta Q_{i}$ can be rewritten as
\begin{align*}
    \Delta Q_{i} 
    {}={}& \frac{\gamma (\m-\n) \tilde{U}^2}{\sqrt{\m^2-1}}
    \sum_{k=j}^{i-1} \sqrt{S_{k+1}} \sin \Big( \sum_{l=k}^{i-1} (-1)^{l-k} \alpha_l\Big).
\end{align*}
Because $\alpha_l$ decreases monotonically on $(0, \frac{\pi}{2})$ and the first term of the sum is always positive, then the sine of the sum is also positive, implying $\Delta Q_i \geq 0$. \\
%
%
%
%
%
%
% \vspace{2em} \\%
\smallskip \\
\textbf{Case $i<j$.} Let us consider a fixed $j\in \{1,2,\dots,N\}$. For each $i \in \{0,1,\dots,j-1\}$, we define 
\begin{align*}%\label{def:Delta_Q_ij_i<j}
    \Delta Q_{i} \coloneqq Q_{i,j} - Q_{i+1,j} - Q_{i,i+1}.
\end{align*}
The distance-1 interpolation inequalities hold exactly, hence $Q_{i+1,i} = 0$ and $Q_{j+1,j} = 0$. 
Therefore to prove $Q_{i,j} \geq 0$ for all $i<j$ it is sufficient to show that $\Delta Q_{i} \geq 0$, since $Q_{i,j}$ is obtained as a sum of positive increments. \\
Similarly to the previous case, substituting the expressions of $Q_{i,j}$, $Q_{i+1,j}$ and $Q_{i,i+1}$, we subsequently obtain:%
\begin{align*}%\label{def:Delta_Q_ij_i>j}
    \Delta Q_{i} 
    {}={}& \frac{1}{L-\mu} 
        \left\langle g_{i+1} - g_j - \mu(x_{i+1} - x_j) \,,\, g_{i+1} - g_{i} - L(x_{i+1} - x_{i}) \right\rangle \\
    % %%%%%%%%%%%%%%%%%%%%%%%%%%%%%%%%%%%%%%%%%%%%%%%%%%%%%%%%%%%%%%%%%%%%%%%%%%%%%%%%%%%%%%%%%%%
    % {}={}& \frac{1}{L-\mu} 
    %     \left\langle g_{i+1} - g_j - \gamma \mu \sum_{k=i+1}^{j-1} g_k \,,\, g_{i+1} - (1-\gamma L) g_i \right\rangle \\
    %%%%%%%%%%%%%%%%%%%%%%%%%%%%%%%%%%%%%%%%%%%%%%%%%%%%%%%%%%%%%%%%%%%%%%%%%%%%%%%%%%%%%%%%%%%
    {}={}& \frac{-\gamma }{\m-\n} \sum_{k=i+1}^{j-1}
        \left\langle g_{k+1} - \m g_{k} \,,\, g_{i+1} - \n g_i \right\rangle.
\end{align*}
From \eqref{eq:zero_gradients_mix} we have that $\Delta Q_{i} = 0$ if $j \leq \bar{N}+1$. Therefore, further on we assume $j \geq \bar{N}+2$. 
Moreover, we have the identity
$g_{k+1} - \n g_{k} = \m^{\bar{N}-k} (g_{\bar{N}+1} - \n g_{\bar{N}})$ for $k = 0,1,\dots, \bar{N}$, 
% and 
% $g_{k+1} - \m g_{k} = \n^{\bar{N}-k} (g_{\bar{N}+1} - \m g_{\bar{N}})$
hence, using \eqref{eq:zero_gradients_mix}, $\Delta Q_{i}$ rewrites as  %
\begin{align*}%\label{def:Delta_Q_ij_i>j}
    \Delta Q_{i}
    %%%%%%%%%%%%%%%%%%%%%%%%%%%%%%%%%%%%%%%%%%%%%%%%%%%%%%%%%%%%%%%%%%%%%%%%%%%%%%%%%%%%%%%%%%%
    {}={}& \frac{-\gamma }{\m-\n} \m^{\max\{0,\bar{N}-i\}} \sum_{k=\max\{i,\bar{N}\}+1}^{j-1}
        \left\langle g_{k+1} - \m g_{k} \,,\, g_{\max\{i,\bar{N}\}+1} - \n g_{\max\{i,\bar{N}\}} \right\rangle.
\end{align*}
From \eqref{eq:eq_mixed_grads_eval_to_use} in the case $i<k$ we get
\begin{align*}%\label{def:Delta_Q_ij_i>j}
    \Delta Q_{i}
    %%%%%%%%%%%%%%%%%%%%%%%%%%%%%%%%%%%%%%%%%%%%%%%%%%%%%%%%%%%%%%%%%%%%%%%%%%%%%%%%%%%%%%%%%%%
    {}={}    
    \frac{\gamma(\m-\n) \tilde{U}^2 \m^{\max\{0,\bar{N}-i\}}}{\sqrt{\m^2-1}}
    \sum_{k=\max\{i,\bar{N}\}+1}^{j-1}    
    \sqrt{S_{k+1}} 
        \sin \left( \sum_{l=\max\{i,\bar{N}\}}^{k-1} (-1)^{k-1-l} \alpha_l \right).
\end{align*}
It is sufficient to check the positivity of the following expression for any $\bar{N} \leq i < j$
\begin{align*}
    V \coloneqq 
    \sum_{k=i+1}^{j-1}    
    \sqrt{S_{k+1}} 
        \sin \Big( \sum_{l=i}^{k-1} (-1)^{l-(k-1)} \alpha_l\Big).
\end{align*}%
Let $\Sigma_{k} \coloneqq \sum_{l=i}^{k-1} (-1)^{l-(k-1)} \alpha_l$ and $Z_{k} {}\coloneqq{} \sqrt{S_{k+1}} \sin \left( \Sigma_{k} \right)$. 
%%%%%%%%%%%%%%%%%%%%%%%%%%%%%%%%%%%%%%%%%%%%%%%%%%%%%%%%%%%%%%%%%%%%%%%%
For any integer $p\geq 0$, one can check the following properties:%
\shownewtextfinal{%
\begin{enumerate}
    \item $\Sigma_{(i+1)+2p} + \Sigma_{(i+1)+2p+1} {}={} \alpha_{(i+1)+2p}$;
    \item $\Sigma_{(i+1)+2p}$ is a decreasing sequence;
    \item $\Sigma_{(i+1)+2p} \geq \alpha_{(i+1)+2p}$;
    \item $\Sigma_{(i+1)+2p} \in (0, \frac{\pi}{2})$. 
\end{enumerate}%
}%
%%%%%%%%%%%%%%%%%%%%%%%%
Then we have%
\[%
\begin{aligned}
%%%%%%%%%%%%%%%%%%%%%%%%%%%%%%%%%%%%%%%%%%%%%%%%%%%%%%%%%%%%%%%%%%%%%%%%%%%%%%%%%%%%%%%%%%%%%%%%%%%%%%%%%%%%%%%%
    Z_{(i+1)+2p} &{}={} \sqrt{S_{(i+1)+2p+1}} \sin ( \Sigma_{(i+1)+2p} ) > 0; \\
%%%%%%%%%%%%%%%%%%%%%%%%%%%%%%%%%%%%%%%%%%%%%%%%%%%%%%%%%%%%%%%%%%%%%%%%%%%%%%%%%%%%%%%%%%%%%%%%%%%%%%%%%
    Z_{(i+1)+2p+1} &{}={} -\sqrt{S_{(i+1)+2p+2}} \sin ( \Sigma_{(i+1)+2p} - \alpha_{(i+1)+2p} ) < 0.
\end{aligned}
\]%
A sufficient condition for the positivity of $V$ is $Z_{(i+1)+2p}  + Z_{(i+1)+2p+1} \geq 0$ for any $p \geq 0$. This condition is equivalent to 
$\sqrt{\frac{S_{(i+1)+2p+1}}{S_{(i+1)+2p+2}}} {}>{}
     \frac{\sin ( \Sigma_{(i+1)+2p} - \alpha_{(i+1)+2p} )}{\sin ( \Sigma_{(i+1)+2p} )}$. \\
Using the monotonic decrease of the cotangent function on $(0, \frac{\pi}{2})$ and the property $\Sigma_{(i+1)+2p} {}\leq \Sigma_{(i+1)} {}={} \alpha_{i}$, we have that
 \begin{align*}
     \frac{\sin ( \Sigma_{(i+1)+2p} - \alpha_{(i+1)+2p} )}{\sin ( \Sigma_{(i+1)+2p} )}
        {}={} &
    \sin (\alpha_{(i+1)+2p}) ( \cot (\alpha_{(i+1)+2p}) - \cot (\Sigma_{(i+1)+2p}) ) \\
        {}\leq{} & 
    % \sin (\alpha_{(i+1)+2p}) ( \cot (\alpha_{(i+1)+2p}) - \cot (\Sigma_{(i+1)}) ) \\
    %     {}={} & 
    \sin (\alpha_{(i+1)+2p}) ( \cot (\alpha_{(i+1)+2p}) - \cot (\alpha_{i}) ).
 \end{align*}
From \eqref{eq:def_angle_alpha} we have
$
     \cot( \alpha_l ) 
        {}={} 
     % \frac{\m \sqrt{S_{l+2}-1} + \sqrt{S_{l+2}}}{\sqrt{\m^2-1}}
     %    {}={}
     \frac{-\m \n \sqrt{S_{l+1}} + \sqrt{S_{l+2}}}{\sqrt{\m^2-1}},
$
for any $l \geq \bar{N}$. Then
\begin{align*}
    \sin (\alpha_{(i+1)+2p}) & ( \cot (\alpha_{(i+1)+2p}) - \cot \alpha_{i} ) \\
    %     {}={} \\
    % \frac{\sqrt{\m^2-1}}{\m \sqrt{S_{(i+1)+2p+2}} - \n \sqrt{S_{(i+1)+2p+1}}} \\
    % \Big( 
    %     \frac{-\m \n ( \sqrt{S_{(i+1)+2p+1}} - \sqrt{S_{i+1}} )  + ( \sqrt{S_{(i+1)+2p+2}} - \sqrt{S_{i+2}} )}{\sqrt{\m^2-1}}
    % \Big)  \\            
    %       &{}={} 
    % \frac{-\m \n ( \sqrt{S_{(i+1)+2p+1}} - \sqrt{S_{i+1}} )  + ( \sqrt{S_{(i+1)+2p+2}} - \sqrt{S_{i+2}} )}{\m \sqrt{S_{(i+1)+2p+2}} - \n \sqrt{S_{(i+1)+2p+1}}} \\
    %%%%%%%%%%%%%%%%%%%%%%%%%%%%%%%%%%%%%%%%%%%%%
          &{}={} 
    \frac{
    -\m \n \left( \sqrt{\frac{S_{(i+1)+2p+1}}{S_{(i+1)+2p+2}}} - \sqrt{\frac{S_{i+1}}{S_{(i+1)+2p+2}}} \right) + 
    \left( 1 - \sqrt{\frac{S_{i+2}}{S_{(i+1)+2p+2}}} \right)
    }{\m - \n \sqrt{\frac{S_{(i+1)+2p+1}} {S_{(i+1)+2p+2} } } },
\end{align*}
where we amplified by $\frac{1}{\sqrt{S_{(i+1)+2p+2}}}$. 
Thus it is left to show that for each $p \geq 0$ it holds
\begin{align*}
    \sqrt{\frac{S_{(i+1)+2p+1}} {S_{(i+1)+2p+2} } }
        {}\geq{}
    \frac{-\m \n \left( \sqrt{\frac{S_{(i+1)+2p+1}}{S_{(i+1)+2p+2}}} - \sqrt{\frac{S_{i+1}}{S_{(i+1)+2p+2}}} \right) + \left( 1 - \sqrt{\frac{S_{i+2}}{S_{(i+1)+2p+2}}} \right)}{\m - \n \sqrt{\frac{S_{(i+1)+2p+1}} {S_{(i+1)+2p+2} } } },
\end{align*}
which is equivalent to:
\begin{align*}
%     \m - \n \sqrt{\frac{S_{(i+1)+2p+1}} {S_{(i+1)+2p+2} } } 
%         {}\geq{} &
%     \frac{-\m \n ( \sqrt{\frac{S_{(i+1)+2p+1}}{S_{(i+1)+2p+2}}} - \sqrt{\frac{S_{i+1}}{S_{(i+1)+2p+2}}} ) + ( 1 - \sqrt{\frac{S_{i+2}}{S_{(i+1)+2p+2}}} )}{ \sqrt{\frac{S_{(i+1)+2p+1}} {S_{(i+1)+2p+2} } } } \\
% %%%%%%%%%%%%%%%%%%%%%%%%%%%%%%%%%%%%%%%%%%%%%%%%%%%%%%%%%%%%%%%%%%%%%%%%%%%%%%%%%%%%%%%%%
%     \big( \m(1+\n) - \n \sqrt{\frac{S_{(i+1)+2p+1}} {S_{(i+1)+2p+2} } } \big) \sqrt{\frac{S_{(i+1)+2p+1}} {S_{(i+1)+2p+2} } } 
%         {}\geq{} &
%     -\m \n ( - \sqrt{\frac{S_{i+1}}{S_{(i+1)+2p+2}}} ) + ( 1 - \sqrt{\frac{S_{i+2}}{S_{(i+1)+2p+2}}} ) \\
% %%%%%%%%%%%%%%%%%%%%%%%%%%%%%%%%%%%%%%%%%%%%%%%%%%%%%%%%%%%%%%%%%%%%%%%%%%%%%%%%%%%%%%%%%%
%     \m(1+\n) - \n \sqrt{\frac{S_{(i+1)+2p+1}} {S_{(i+1)+2p+2} } } 
%         {}\geq{} &
%     \m \n \sqrt{\frac{S_{i+1}}{S_{(i+1)+2p+1}}} + \sqrt{\frac{S_{(i+1)+2p+2} } {S_{(i+1)+2p+1}}  } - \sqrt{\frac{S_{i+2}}{S_{(i+1)+2p+1}}} \\
%%%%%%%%%%%%%%%%%%%%%%%%%%%%%%%%%%%%%%%%%%%%%%%%%%%%%%%%%%%%%%%%%%%%%%%%%%%%%%%%%%%%%%%%%%
    \m(1+\n) - \left(\n + \frac{S_{(i+1)+2p+2} } {S_{(i+1)+2p+1}}\right) \sqrt{\frac{S_{(i+1)+2p+1}} {S_{(i+1)+2p+2} } }
        {}\geq{} &
    \m \n \sqrt{\frac{S_{i+1}}{S_{(i+1)+2p+1}}} - \sqrt{\frac{S_{i+2}}{S_{(i+1)+2p+1}}}.
\end{align*}
Using $\frac{S_{(i+1)+2p+2}}{S_{(i+1)+2p+1}} = \n^2 + \frac{1}{S_{(i+1)+2p+1}}$, the inequality equivalently rewrites as:
\begin{align*}
%%%%%%%%%%%%%%%%%%%%%%%%%%%%%%%%%%%%%%%%%%%%%%%%%%%%%%%%%%%%%%%%%%%%%%%%%%%%%%%%%%%%%%%%%%
    % \m(1+\n) - (\n + \n^2 + \frac{1} {S_{(i+1)+2p+1}}) \sqrt{\frac{S_{(i+1)+2p+1}} {S_{(i+1)+2p+2} } }
    %     {}\geq{} 
    % \m \n \sqrt{\frac{S_{i+1}}{S_{(i+1)+2p+1}}} - \sqrt{\frac{S_{i+2}}{S_{(i+1)+2p+1}}} \\
%%%%%%%%%%%%%%%%%%%%%%%%%%%%%%%%%%%%%%%%%%%%%%%%%%%%%%%%%%%%%%%%%%%%%%%%%%%%%%%%%%%%%%%%%%
    (1+\n) 
    \left(\m-\n \sqrt{\frac{S_{(i+1)+2p+1}} {S_{(i+1)+2p+2} } } \right)
    + \frac{\frac{-1} { \sqrt{ S_{(i+1)+2p+2} } } 
    % - \sqrt{\frac{1} { S_{(i+1)+2p+1} S_{(i+1)+2p+2} } }
    + \sqrt{S_{i+2}}
    - \m \n \sqrt{S_{i+1}}} { \sqrt{ S_{(i+1)+2p+1} } }
        {}\geq{} 
    0.
%%%%%%%%%%%%%%%%%%%%%%%%%%%%%%%%%%%%%%%%%%%%%%%%%%%%%%%%%%%%%%%%%%%%%%%%%%%%%%%%%%%%%%%%%%
\end{align*}
The first term is positive because $\n \in (-1,0)$, whereas the second one due to $S_{i+2} \geq 1$ and $\frac{1}{\sqrt{S_{(i+1)+2p+2}}} < 1$.%
\qed%%
%
}%
%\end{appendices}
% BibTeX users please use one of
%\bibliographystyle{spbasic}      % basic style, author-year citations
% \bibliographystyle{spmpsci}      % mathematics and physical sciences
% \bibliographystyle{spphys}       % APS-like style for physics
\bibliographystyle{alpha} % my own (FG) preference
\bibliography{ms.bib}   % name your BibTeX data base
%% EXAMPLE TEMPLATE
%\section{Introduction}
%\label{intro}
%Your text comes here. Separate text sections with
%\section{Section title}
%\label{sec:1}
%Text with citations \cite{RefB} and \cite{RefJ}.
%\subsection{Subsection title}
%\label{sec:2}
%as required. Don't forget to give each section
%and subsection a unique label (see Sect.~\ref{sec:1}).
%\paragraph{Paragraph headings} Use paragraph headings as needed.
%\begin{equation}
%a^2+b^2=c^2
%\end{equation}
%
%% For one-column wide figures use
%\begin{figure}
%% Use the relevant command to insert your figure file.
%% For example, with the graphicx package use
%  \includegraphics{example.eps}
%% figure caption is below the figure
%\caption{Please write your figure caption here}
%\label{fig:1}       % Give a unique label
%\end{figure}
%%
%% For two-column wide figures use
%\begin{figure*}
%% Use the relevant command to insert your figure file.
%% For example, with the graphicx package use
%  \includegraphics[width=0.75\textwidth]{example.eps}
%% figure caption is below the figure
%\caption{Please write your figure caption here}
%\label{fig:2}       % Give a unique label
%\end{figure*}
%%
%% For tables use
%\begin{table}
%% table caption is above the table
%\caption{Please write your table caption here}
%\label{tab:1}       % Give a unique label
%% For LaTeX tables use
%\begin{tabular}{lll}
%\hline\noalign{\smallskip}
%first & second & third  \\
%\noalign{\smallskip}\hline\noalign{\smallskip}
%number & number & number \\
%number & number & number \\
%\noalign{\smallskip}\hline
%\end{tabular}
%\end{table}
\end{document}